# The existence of instanton solutions to the $\mathbb{R}$-invariant Kapustin-Witten equations on $(0, \infty) \times \mathbb{R}^2 \times \mathbb{R}$

Clifford Henry Taubes[†]

Department of Mathematics
Harvard University
Cambridge, MA 02138

ABSTRACT: A non-negative integer labeled set of model solutions to the $\mathbb{R}$-invariant Kapustin-Witten equations on $(0,\infty) \times \mathbb{R}^2 \times \mathbb{R}$ plays a central role in Edward Witten's program to interpret the colored Jones polynomials of a knot in the context of SU(2) gauge theory. This paper explains why there are $\mathbb{R}$-invariant solutions to these equations on $(0,\infty) \times \mathbb{R}^2 \times \mathbb{R}$ that interpolate between two model solutions as the $(0, \infty)$ parameter increases from 0 to $\infty$ while respecting the $\mathbb{R}^2$ factor asymptotics. The only constraint on the limiting pair of model solutions is this: Letting $m_0$ and $m_\infty$ denote their respective non-negative integer labels, then $m_0 - m_\infty$ must be a positive, even integer. (As explained in the paper, there is a $\mathbb{C}^{(m_0 - m_\infty - 2)/2} \times \mathbb{C}^*$ moduli space of these interpolating solutions.)

[†] Supported in part by the NSF under the grant DMS-2002771

# 1. Introduction

The purpose of this paper is to first state and then prove a theorem asserting the existence of what I will call *instanton* solutions to the ($\mathbb{R}$-invariant) Kapustin-Witten equations on $(0,\infty) \times \mathbb{R}^2 \times \mathbb{R}$ for the group SU(2). These are solutions that interpolate in a suitable sense between distinct versions of the model solutions found by Edward Witten [W1] (see also [GT] and [MW]. A precise characterization of an instanton solution is given momentarily in Section 1c. Section 1a which follows directly describes the Kapustin-Witten equations and Section 1b describes Witten's model solutions. The existence theorem for these solutions is also in Section 1c. The remaining sections 2–7 of this paper contain the proof of the existence theorem.

It is not clear at present if or how these instanton solutions fit in to Witten's program (see [W1-4]) to obtain the colored Jones polynomial via an algebraic count of solutions to the Kapustin-Witten equations.

## a) The $S^1$-invariant Kapustin-Witten equations

To set the notation: The space $(0,\infty) \times \mathbb{R}^2 \times \mathbb{R}$ is given the product metric. The notation has t denoting the Euclidean coordinate for the $(0,\infty)$ factor. Let P denote a principal SU(2) bundle over $(0,\infty) \times \mathbb{R}^2 \times \mathbb{R}$ and let A denote a given connection on P. The notation uses $\nabla_{At}$ to denote the A-covariant derivative along the $(0,\infty)$ factor of $(0,\infty) \times \mathbb{R}^2 \times \mathbb{R}$. It also uses $d_A$ to denote the A-exterior derivative along the $\mathbb{R}^2 \times \mathbb{R}$ factor. The Hodge star along this same factor is denoted by $*$. Finally the curvature of A is written below as $dt \wedge E_A + *B_A$ with $E_A$ and $B_A$ being sections of the bundle $\text{ad}(P) \otimes T^*(\mathbb{R}^2 \times \mathbb{R})$. (The bundle $\text{ad}(P)$ is the associated vector bundle to P using the adjoint representation of SU(2) on its Lie algebra.)

For the purpose of this paper, a solution to the Kapustin-Witten equations on $(0,\infty) \times \mathbb{R}^2 \times \mathbb{R}$ consists of a principle SU(2) bundle over this space (denoted by P), and a pair $(A, \mathfrak{a})$ of a connection on P and section of $\text{ad}(P) \otimes T^*(\mathbb{R}^2 \times \mathbb{R})$ obeying the following:

- $E_A = *d_A \mathfrak{a}$ .
- $\nabla_{At} \mathfrak{a} = B_A - *(\mathfrak{a} \wedge \mathfrak{a})$ .
- $d_A * \mathfrak{a} = 0$.

(1.1)

A pair $(A, \mathfrak{a})$ of connection on P and section of $\text{ad}(P) \times T^*(\mathbb{R}^2 \times \mathbb{R})$ is deemed to be $\mathbb{R}$ invariant when A's curvature annihilates tangent vectors to the $\mathbb{R}$ factor of $(0,\infty) \times \mathbb{R}^2 \times \mathbb{R}$ and when A's directional covariant derivative along this same $\mathbb{R}$ factor annihilates $\mathfrak{a}$.

It is convenient to rewrite the equations in (1.1) for an $\mathbb{R}$ invariant pair $(A, \mathfrak{a})$ so as to distinguish the invariant $\mathbb{R}$ direction. To this end: Let $(z_1, z_2)$ denote Euclidean coordinates for the $\mathbb{R}^2$ factor. (When a complex coordinate is needed, it is $z = z_1 + i z_2$.)
The Euclidean coordinate along the $\mathbb{R}$ factor of $(0,\infty) \times \mathbb{R}^2 \times \mathbb{R}$ is denoted by $x_3$.
The components of $E_A$ and $B_A$ with respect to the basis $\{dz_1, dz_2, dx_3\}$ for $T^*((0,\infty) \times \mathbb{R}^2 \times \mathbb{R})$ are written as $(E_{A1}, E_{A2}, E_{A3})$ and $(B_{A1}, B_{A2}, B_{A3})$; and likewise the components of $\mathfrak{a}$ are written as $(\mathfrak{a}_1, \mathfrak{a}_2, \mathfrak{a}_3)$. Meanwhile, A's directional covariant derivative along the dual basis of vector fields as $\nabla_{A1}, \nabla_{A2}$ and $\nabla_{A3}$. To say that $(A,\mathfrak{a})$ is $\mathbb{R}$-invariant is to say that

$$E_{A3} = B_{A1} = B_{A2} = 0 \quad \text{and} \quad \nabla_{A3}\mathfrak{a} = 0.$$

(1.2)

Now introduce $\varphi \equiv \mathfrak{a}_1 - i\mathfrak{a}_2$ and $\varphi^* = \mathfrak{a}_1 + i\mathfrak{a}_2$. Given (1.2), then the equations in (1.1) can be written as follows:

- $E_{A1} = \nabla_{A2}\mathfrak{a}_3$ and $E_{A2} = -\nabla_{A1}\mathfrak{a}_3$.
- $(\nabla_{A1} + i\nabla_{A2})\varphi = 0$.
- $\nabla_{At}\varphi = i[\mathfrak{a}_3, \varphi]$.
- $\nabla_{At}\mathfrak{a}_3 = B_{A3} + \frac{i}{2}[\varphi, \varphi^*]$.

(1.3)

These are the equations of interest for what is to come.

**b) Witten's model solutions**

Edward Witten introduced in [W1] the model solutions to (1.3) that are depicted below in (1.5) as pairs of connection on the product principal SU(2) bundle and section of $T^*((0,\infty) \times \mathbb{R}^2 \times \mathbb{R})$ with values in the product Lie algebra bundle.

To set the notation: The Lie algebra of SU(2), denoted by $\mathfrak{su}(2)$, is viewed as the space of $2 \times 2$ anti-Hermitian, traceless $\mathbb{C}$-valued matrices. A basis $\{\sigma_1, \sigma_2, \sigma_3\}$ for $\mathfrak{su}(2)$ is hereby chosen to obey the following algebraic relations:

- $\sigma_1^2 = \sigma_2^2 = \sigma_3^2 = -1$
- $\sigma_1\sigma_2 = -\sigma_3$ and $\sigma_2\sigma_3 = -\sigma_1$ and $\sigma_3\sigma_1 = -\sigma_2$.

(1.4)

Note the - sign in the second bullet. This is not the convention used by most people. The bi-invariant inner product on $\mathfrak{su}(2)$ is defined by declaring this basis to be orthonormal. Letting $\theta_0$ denote the product connection on the product principal SU(2) bundle, the basis $\{\sigma_1, \sigma_2, \sigma_3\}$ for $\mathfrak{su}(2)$ is viewed as a $\theta_0$-covariantly constant, orthonormal basis for the associated Lie algebra bundle (which is also a product bundle).

More notation: A real valued function $\Theta$ on $(0,\infty)\times(\mathbb{R}^2-0)$ is defined by setting $\sinh(\Theta) = \frac{t}{|z|}$. Meanwhile, $x$ denotes the function on $(0,\infty)\times\mathbb{R}^2\times\mathbb{R}$ given by the rule $(t,z) \to (t^2+|z|^2)^{1/2}$.

Witten's model solutions (from [W1]) are indexed by a non-negative integer which is denoted by $m$ in what follows. The integer $m$ version is this:

- *The section* $\mathfrak{a}$:
  a) $\mathfrak{a}_3 = -\frac{1}{2t} \frac{(m+1)\sinh(\Theta)}{\sinh((m+1)\Theta)} \frac{\cosh((m+1)\Theta)}{\cosh(\Theta)} \sigma_3$.
  b) $\varphi = \mathfrak{a}_1 - i\mathfrak{a}_2 = -\frac{1}{2t} \frac{(m+1)\sinh(\Theta)}{\sinh((m+1)\Theta)} (\frac{z}{|z|})^m (\sigma_1 - i\sigma_2)$.
- *Connection* $A$ *and its curvature*:
  a) $A = \theta_0 + \frac{(m+1)}{2}(1 - \frac{\sinh(\Theta)}{\cosh(\Theta)}\frac{\cosh((m+1)\Theta)}{\sinh((m+1)\Theta)}) \frac{1}{|z|^2}(z_1 dz_2 - z_2 dz_1)\sigma_3$.
  b) $B_A = \frac{(m+1)}{2x^2} \frac{\sinh(\Theta)}{\cosh(\Theta)} \frac{\cosh((m+1)\Theta)}{\sinh((m+1)\Theta)} (1 - \frac{(m+1)\sinh(\Theta)\cosh(\Theta)}{\sinh((m+1)\Theta)\cosh((m+1)\Theta)}) \sigma_3\, dx_3$
  c) $E_A = -\frac{(m+1)}{2x^2} \frac{\cosh((m+1)\Theta)}{\sinh((m+1)\Theta)} (1 - \frac{(m+1)\sinh(\Theta)\cosh(\Theta)}{\sinh((m+1)\Theta)\cosh((m+1)\Theta)}) \sigma_3 \frac{1}{x}(z_1 dz_2 - z_2 dz_1)$.

(1.5)

By way of terminology: The $m = 0$ version of (1.2) is called the Nahm pole solution; it is the solution with $A = \theta_0$ and $\{\mathfrak{a}_a = -\frac{1}{2t}\sigma_a\}_{a=1,2,3}$.

These solutions have special properties which are listed next and in (1.7) and (1.8).

- *The Lie-algebra element* $\sigma_3$ *is* $A$-*covariantly constant; so the connection* $A$ *is Abelian*.
- *The component* $\mathfrak{a}_3$ *is proportional to* $\sigma_3$; *it will be written as* $\alpha\sigma_3$. *Meanwhile,* $\alpha$ *is negative with values between* $-\frac{1}{2t}$ *and* $-\frac{(m+1)}{2t}$.
- *The* $\mathfrak{su}(2)\times_{\mathbb{R}}\mathbb{C} = \mathfrak{sl}(2;\mathbb{C})$ *valued function* $\varphi$ *is pointwise orthogonal to* $\sigma_3$. *Moreover*
  a) $\mathrm{trace}(\varphi^2) \equiv 0$. *(Thus* $\mathfrak{a}_1$ *and* $\mathfrak{a}_2$ *are pointwise orthogonal and have the same norm.)*
  b) $[\frac{i}{2}\sigma_3, \varphi] = \varphi$.
  c) $|\varphi| \leq \frac{1}{\sqrt{2}t}$ *with equality only in the case when* $m = 0$.
- *The norms of the* $\mathfrak{su}(2)$ *valued 1-forms* $E_{A1}, E_{A2}$ *and* $B_{A3}$ *and* $\nabla_{A1}\mathfrak{a}$ *and* $\nabla_{A2}\mathfrak{a}$ *are bounded by a constant multiple of* $\frac{t}{x^3}$. *Meanwhile, the norm of* $\nabla_t\mathfrak{a}$ *is bounded by a constant multiple of* $\frac{1}{t^2}$.
- *The function* $\langle\sigma B_{A3}\rangle$ *is positive (it is the norm of* $B_{A3}$*); and its integral over each constant* $t \in (0,\infty)$ *version of* $\{t\}\times\mathbb{R}^2\times\{0\}$ *is finite and equal to* $2\pi m$.

(1.6)

The next property concerns the family of coordinate rescaling diffeomorphisms of the $(0,\infty)\times\mathbb{R}^2$ factor in $(0,\infty)\times\mathbb{R}^2\times\mathbb{R}$; the family is parametrized by $(0,\infty)$ and its action is defined by the rule whereby the diffeomorphism that is parametrized by a given positive number $\lambda$ sends $(t, z_1, z_2, x_3) \to (\lambda t, \lambda z_1, \lambda z_2, \lambda x_3)$.

*The model solutions from (1.5) are fixed by every coordinate rescaling diffeomorphism from this family: Each model solution is the same as its pull-back by any of them.*

(1.7)

A final property: Write $\mathfrak{a}_3$ as $\alpha \sigma_3$ with $\alpha \equiv -\frac{1}{2t} \frac{(m+1)\sinh(\Theta)}{\sinh((m+1)\Theta)} \frac{\cosh((m+1)\Theta)}{\cosh(\Theta)}$. Then the equations in (1.2) when written using $(\alpha, \varphi)$ and $(E_{A1}, E_{A2}, B_{A3})$ say this:

- $\nabla_t \varphi - 2\alpha \varphi = 0$  *and*  $(\nabla_1 + i\nabla_2)\varphi = 0$.
- $E_{A1} = \frac{\partial \alpha}{\partial z_2} \sigma_3$  *and*  $E_{A2} = -\frac{\partial \alpha}{\partial z_1} \sigma_3$.
- $B_{A3} = (\frac{\partial \alpha}{\partial t} - |\varphi|^2) \sigma_3$.

(1.8)

The integer $m$ model solution from (1.5) is denoted henceforth as $(A^{(m)}, \mathfrak{a}^{(m)})$.

### c) Kapustin-Witten instantons

Kapustin-Witten instantons are $\mathbb{R}$-invariant solutions to (1.3) that obey the following constraints (this is CONSTRAINT SET 2 in [T1]):

- *The function* $t|\mathfrak{a}|$ *is uniformly bounded on* $(0, \infty) \times \mathbb{R}^2 \times \mathbb{R}$.
- *There exist positive numbers* $t_0$, $\varepsilon$ *and* $r > 0$ *such that* $t|\mathfrak{a}| > \varepsilon$ *where both* $t < t_0$ *and* $|z| \geq rt$ *in* $(0, \infty) \times \mathbb{R}^2 \times \mathbb{R}$.
- *The function* $|B_{A3}|^2 + |E_{A1}|^2 + |E_{A2}|^2$ *has finite integral on any given* $R > 0$ *version of the domain in* $(0, \infty) \times \mathbb{R}^2 \times \{0\}$ *where* $|z| > R$; *and that integral is bounded by an* R-*independent multiple of* $\frac{1}{R}$.

(1.9)

These constraints are satisfied by the model solutions. Of course, if P is a principal SU(2) bundle over $(0, \infty) \times \mathbb{R}^2 \times \mathbb{R}$, then P is isomorphic to $P_0$ and the pull-back of any model solution is also a solution (for the most part, no distinction will be made between a model solution on $P_0$ and its pull-back by a principle bundle isomorphism).

According to Theorem 2 in [T1], if a solution to (1.2) and (1.3) for a given principal bundle P obeys the constraints in (1.9), and if it is not the pull-back of some model solutions from (1.5) via a principle bundle isomorphism, then it none-the-less looks like a model solution where $t \ll 1$ and where $t \gg 1$ and where $\frac{|z|}{t} \gg 1$. To be precise, Theorem 2 in [T1] says that there exist a non-negative integer $m$, a positive integer $p$ and isomorphisms $g_\infty$ and $g_0$ from P to $P_0$ defined on the respective $t > 1$ and $t < 1$ parts of $(0, \infty) \times \mathbb{R}^2 \times \mathbb{R}$, and an isomorphism $g_>$ from $P/\{\pm 1\}$ to $P_0/\{\pm 1\}$ that is defined in the $|z| > t$ part; and these are such that

- $\lim_{t \to 0} (t|g_0^* A - A^{(m+2p)}| + t|g_0^* \mathfrak{a} - \mathfrak{a}^{(m+2p)}|) = 0$,

- $\lim_{t\to\infty} (t|g_\infty{}^*A - A^{(m)}| + t|g_\infty{}^*\mathfrak{a} - \mathfrak{a}^{(m)}| + |\ln(\frac{|\varphi|}{|\varphi^{(m)}|})|) = 0$ ,

- $\lim_{|z|/t\to\infty} (t|g_>{}^*A - A^{(0)}| + t|g_>{}^*\mathfrak{a} - \mathfrak{a}^{(0)}|) = 0$ .

(1.10)

(The limits in the top two bullets are uniform with respect to the $\mathbb{R}^2$ coordinate z; and the limit in the third bullet is uniform with respect to $x = (t^2+|z|^2)^{1/2}$.) As a parenthetical remark: The $t \to \infty$ limit convergence is stronger than the $t \to 0$ limit convergence by virtue of the convergence of the norm of the logarithm of the ratio of $|\varphi|$ to $|\varphi^{(m)}|$. It is a consequence of what is said in Proposition 3.1 of [T1] that the logarithm of the ratio of $|\varphi|$ to $|\varphi^{(m+2p)}|$ diverges on the whole of the $z = 0$ locus.

The analysis in [T1] did not determine whether solutions obeying (1.10) exist. The theorem below says they do.

**Theorem 1**: *Fix a non-negative integer* m *and a positive integer* p. *Then there exist solutions to (1.3) that are described by (1.9) and (1.10). Moreover, there is a moduli space of solutions to (1.3) obeying (1.9) and (1.10) which is parametrized by* $\mathbb{C}^{p-1} \times (\mathbb{C}-0)$ *with the property that respective solutions labeled by distinct points in* $\mathbb{C}^{p-1} \times (\mathbb{C}-0)$ *are not pull-backs of each other by any principle bundle isomorphism (which is to say that they are not gauge equivalent).*

Much more is said about the behavior of these solutions in [T1]. As noted, the proof of Theorem 1 occupies the rest of this paper; the last section (Section 7) summarizes the arguments.

By way of a look ahead, the proof of Theorem 1 exploits an interpretation of the the equations in (1.3) which is due to Gaiotto and Witten [GW], the top three equations being the integrability equations for a vector bundle complex structure (a sort of $\bar{\partial}$-operator). The fourth equation in (1.3) is then viewed as a constraint on the data that defines the $\bar{\partial}$ operator. The complex interpretation of the equations in (1.3) was used to great effect by Mazzeo and Witten (see Section 2 of [MW]) and in the fundamental papers of He and Mazzeo in [HM1-HM3].

Ideas introduced by Donaldson in [D] in the context of the Hermitian-Yang-Mills equations and anti-self-dual Yang-Mills equations are antecedents of this complex structure interpretation of (1.3). Of particular note is that Donaldson solved his constraint equation (the Hermitian Yang-Mills analog of the fourth bullet in (1.3)) by deforming an initial choice of data for the vector bundle $\bar{\partial}$-operator along a 1-parameter family of complex equivalent $\bar{\partial}$-operators using notions from Hermitian geometry. The idea of solving the fourth bullet equation via a 1-parameter family of complex equivalences is borrowed for use here, but it is not done here using Hermitian geometric tools. (One could probably use a Hermitian geometric approach here also, but that seemed daunting to me. He and Mazzeo [HM1-HM3] do brave the Hermitian geometric approach along the lines used by Donaldson.)

TABLE OF CONTENTS FOR THE REST OF THIS ARTICLE

By way of an acknowledgement, I thank Weifeng Sun for his thoughts about the complex interpretation of the equations in (1.3).

**d) Notation and conventions**

With regards to notation: The notation has $c_0$ denoting a number that is greater than 1 and independent of what ever is relevant to a given inequality. For example, if an inequality concerns a value or values of t from $[0,\infty)$ and/or points in $\mathbb{R}^2$, then $c_0$ will be independent of t and the points in $\mathbb{R}^2$. In general, it should be clear from the context what does and doesn't determine an upper bound for $c_0$. It is always the case that $c_0$ increases between successive appearances.

With regards to SU(2) conventions: The Lie algebra of SU(2) (denoted by $\mathfrak{su}(2)$) is identified implicitly with the vector space of $2 \times 2$, anti-hermitian matrices. If $\sigma$ is any given $2\times 2$, $\mathbb{C}$-valued matrix, the notation $\langle\sigma\rangle$ means this:

$$\langle\sigma\rangle = -\tfrac{1}{2}\,\mathrm{trace}(\sigma)\,.$$
(1.11)

The inner product of two elements, say $\sigma$ and $\tau$, in $\mathfrak{su}(2)$ is then $\langle\sigma\tau\rangle$. This inner product on $\mathfrak{su}(2)$ defines a fiberwise inner product on the bundle ad(P). The inner product between elements $\sigma$ and $\tau$ in the fiber of ad(P) over any point is denoted by $\langle\sigma\tau\rangle$ also.

The complexification of $\mathfrak{su}(2)$ is the vector space of $2 \times 2$ complex matrices with zero trace which is the Lie algebra of the group $SL(2;\mathbb{C})$. The vector bundle associated to P via the adjoint representation of SU(2) on this complexification is denoted by $\mathrm{ad}(P)\otimes_{\mathbb{R}}\mathbb{C}$. If $\eta$ is an element of the complexification of $\mathfrak{su}(2)$ or of $\mathrm{ad}(P)\otimes_{\mathbb{R}}\mathbb{C}$, then $\eta^*$ is used to denote -1 times the Hermitian conjugate of η (thus, -η†). With this notation understood, then the Hermitian inner product between two elements (call the first one η and the second one λ) in the complexified Lie algebra or in $ad(P)\otimes_{\mathbb{R}}\mathbb{C}$ is $\langle \eta^*\lambda \rangle$,

Inner products on vectors, differential forms and higher rank tensors on $(0,\infty)\times\mathbb{R}^2\times\mathbb{R}$ are defined using the Euclidean metric that is defined so that dt, $dx_1$, $dx_2$ and $dx_3$ are orthonormal. The inner product between tensor $\mathfrak{t}$ and $\mathfrak{h}$ of any given rank is denoted by $\langle \mathfrak{t},\mathfrak{h} \rangle$. This Euclidean inner product and the previously defined fiber inner product for ad(P) gives an inner product on ad(P) valued tensors of any rank. If $\mathfrak{t}$ and $\mathfrak{h}$ are any two such tensors, their inner product is also denoted by $\langle \mathfrak{t},\mathfrak{h} \rangle$.

With regards to covariant derivatives: Supposing that $\mathfrak{t}$ is an tensor, its covariant derivative (using the Euclidean metric's flat connection) is denoted by $\nabla\mathfrak{t}$. Supposing that $\mathfrak{t}$ is an ad(P) valued tensor, the covariant derivative using the connnection A and the Euclidean metric's connection is denoted by $\nabla_A$. (It has components $\nabla_{At}$, $\nabla_{A1}$, $\nabla_{A2}$ and $\nabla_{A3}$ with respect to the basis dt, $dx_1$, $dx_2$ and $dx_3$. But $\nabla_{A3}$ will act as zero on every tensor and ad(P) valued tensor of interest.)

Another convention is with regards to 'cut-off functions' and/or 'bump functions': All such functions are constructed from a basic model function on $\mathbb{R}$ to be denoted by χ. This function χ is a smooth, non-increasing function that equals 1 on $(-\infty, \frac{1}{4}]$ and equals zero on $[\frac{3}{4},\infty)$. Cut-off functions in the context of $(0,\infty)\times\mathbb{R}^2)$ are obtained from χ by composing the latter with a suitably chosen map to $\mathbb{R}$. The advantage to constructing all such functions from χ is that the norms of the derivatives of any such function has a priori bounds given bounds for the norms of the derivatives of the chosen map to $\mathbb{R}$.

One last convention: For the most part, the $\mathbb{R}$ factor in $\mathbb{R}^2\times\mathbb{R}$ will be ignored in what follows because its only role is to supply the $x_3$ components of $\mathfrak{a}$ and $B_A$.

## 2. Rewriting (A, $\mathfrak{a}$)

Let (A, $\mathfrak{a}$) denote a given pair of connection on P and ad(P)-valued 1-form along the $\mathbb{R}^2\times\mathbb{R}$ factor of $(0,\infty)\times\mathbb{R}^2\times\mathbb{R}$ that obeys (1.2) and the first three bullets of (1.3). It is also assumed in what follows that the zero locus of $\varphi \equiv \mathfrak{a}_1 - i\mathfrak{a}_2$ is the z = 0 locus. This section rewrites (A, $\mathfrak{a}$) in terms a new sort of data set. Looking ahead: This rewriting of solutions to the first three bullets of (1.3) is reversed in Section 3 to construct pairs that obey the first three bullets of (1.3) and are described by (1.9) and (1.10), thus accounting for all but one of Theorem 1's requirements (missing only the fourth bullet of (1.3)). The last subsection here explains how the rewriting done here ties into the complex vector bundle picture of Giaotto and Witten [GW], [MW] and [HM1-3]. (This rewriting of the pair (A, $\mathfrak{a}$) summarizes part of what is done in Section 3 of [T1].)

### a) The section σ

Some background: A length 1 section of the bundle ad(P) can be used to decompose $ad(P) \otimes_\mathbb{R} \mathbb{C}$ as a direct sum of eigenbundles for the section's commutator action as an endomorphism of ad(P): Letting σ denote the section in question, the decomposition is

$$ad(P) \otimes_\mathbb{R} \mathbb{C} = \mathcal{L}^+ \oplus \mathbb{C}\sigma \oplus \mathcal{L}^-$$
(2.1)

with $\mathcal{L}^+$ and $\mathcal{L}^-$ denoting the respective +1 and -1 eigenbundles of the endomorphism $[\frac{i}{2}\sigma, \cdot]$. These are complex line bundles over the domain of σ (which is implicitly assumed to be an open set in $(0,\infty) \times \mathbb{R}^2 \times \mathbb{R}$). Two observations are used throughout, the first being that if η is in, $\mathcal{L}^+$ then $\langle \eta^2 \rangle = 0$ and $[\eta, \eta] = 0$. The second observation is that $\eta^*$ is in $\mathcal{L}^-$ and $[\eta, \eta^*] = -2i|\eta|^2 \sigma$. (Remember that $\eta^*$ is -1 times η's Hermitian conjugate.)

Now suppose that $(A, \mathfrak{a})$ obeys (1.2) and the top three bullets of (1.3) with φ not identically zero. Assume in addition that $t|\mathfrak{a}|$ is bounded on $(0,\infty) \times \mathbb{R}^2 \times \mathbb{R}$. Granted these assumptions, Proposition 3.1 in [T1] describes a unit length section of the bundle ad(P), denoted by σ, that obeys

$$\langle \sigma\varphi \rangle = 0 \quad and \quad [\tfrac{i}{2}\sigma, \varphi] = \varphi .$$
(2.2)

These imply that φ is a section of the corresponding line bundle $\mathcal{L}^+$. (By way of an example: When $(A, \mathfrak{a})$ is a model solution from (1.5), then the corresponding σ is $\sigma_3$.)

### b) The function α, the connection Â and the sections β and 𝔟

The unit length section σ can be used to write $\mathfrak{a}_3$ as

$$\mathfrak{a}_3 = \alpha\sigma + \beta + \beta^* ,$$
(2.3)

with α being a real valued function on $(0,\infty) \times \mathbb{R}^2 \times \mathbb{R}$ and with β being the $\mathcal{L}^+$ part of $\mathfrak{a}_3$. (The section β vanishes for the model solutions from (1.5).)

As explained directly, this section β of $\mathcal{L}^+$ also appears as part of the connection A. To elaborate, introduce a connection to be denoted by Â by writing

$$A = \hat{A} + \tfrac{1}{4}[\sigma, \nabla_A \sigma] .$$
(2.4)

The connection Â is defined so that $\nabla_{\hat{A}} \sigma = 0$. This writes A as $\hat{A} + \mathfrak{b}$ with 𝔟 denoting an ad(P)-valued 1-form obeying $\langle \sigma \mathfrak{b} \rangle = 0$. This decomposition of A is useful because $\nabla_{\hat{A}}$ preserves the summands in (2.1). (The 𝔟 part of A is zero when A is a model solution).

If (A, 𝔞) obeys (1.2) and the first three bullets of (1.3) and φ isn't identically zero, then the dt and $dz_1$ and $dz_2$ components of $\frac{1}{4}[\sigma, \nabla_A \sigma]$ must behave as follows:

- $\frac{1}{4}[\sigma, \nabla_{At}\sigma] = -i(\beta - \beta^*)$ .
- $\flat \equiv \frac{1}{8}([\sigma, \nabla_{A1}\sigma] + [\sigma, i\nabla_{A2}\sigma])$ is a section of $\mathcal{L}^+$ .

(2.5)

This behavior is explained in Section 3c of [T1]. Section 3c of [T1] also explains why $[\sigma, \nabla_A \sigma]$ has no $dx_3$ component.

**b) Rewriting the equations in (1.3)**

The equations in (1.2) and (1.3) when written using (φ, α, Â, β, $\flat$) are equivalent to the following system of equations:

- *The top bullet in (1.3)*:
  a) $\langle \sigma E_{\hat{A}1} \rangle = \partial_2 \alpha$ and $\langle \sigma E_{\hat{A}2} \rangle = -\partial_1 \alpha$ .
  b) $(\nabla_{\hat{A}t} - 2\alpha) \flat = -i(\nabla_{\hat{A}1} + i\nabla_{\hat{A}2})\beta$.
- *The second and third bullets in (1.3)*:
  a) $(\nabla_{\hat{A}1} + i\nabla_{\hat{A}2})\varphi = 0$ .
  b) $\nabla_{\hat{A}t}\varphi - 2\alpha\varphi = 0$ .
- *The fourth bullet in (1.3)*:
  a) $\partial_t \alpha = \langle \sigma B_{\hat{A}3} \rangle + |\varphi|^2 + 4(|\beta|^2 + |\flat|^2)$ .
  b) $(\nabla_{\hat{A}t} + 2\alpha)\beta = -i(\nabla_{\hat{A}1} - i\nabla_{\hat{A}2})\flat$.

(2.6)

The equation in Item a) of the second bullet of (2.6) says in effect that α can be written where φ is non-zero as

$$\alpha = \partial_t w .$$

(2.7)

with w defined where $\varphi \neq 0$ by writing |φ| as $e^{2w}$. With regards to the second bullet of (2.6): It's two equations can be solved simultaneously only in the event that the differential operators $\nabla_{\hat{A}t} - 2\alpha$ and $\nabla_{\hat{A}1} + i\nabla_{\hat{A}2}$ commute on sections of $\mathcal{L}^+$. The condition in Item a) of the top bullet of (2.6) is necessary and sufficient for that to happen. The second bullet in (2.6), this says in effect that there is a principle bundle isomorphism on the complement of the φ = 0 locus in $(0, \infty) \times \mathbb{R}^2 \times \mathbb{R}$ that pulls back the connection Â so as to have the form

$$\hat{A} = \theta_{\ddagger} + (\partial_2 w\, dz_1 - \partial_1 w\, dz_2)\sigma .$$

(2.8)

with $\theta_\ddagger$ denoting the unique flat connection on the complement of the $\varphi = 0$ locus that annhilates $\sigma$ and $\frac{\varphi}{|\varphi|}$. It then follows as a consequence of this depiction of $\hat{A}$ that $B_{\hat{A}3}$ can be written where $\varphi \neq 0$ as

$$B_{\hat{A}3} = -(\partial_1^2 + \partial_2^2)\, w\, \sigma.$$
(2.9)

Granted this depiction of $B_{\hat{A}3}$ and granted (2.7), then Item a) of the third bullet of (2.6) becomes a second order differential equation for w:

$$-(\partial_t^2 + \partial_1^2 + \partial_2^2)\, w + e^{4w} + 4(|\beta|^2 + |\mathfrak{b}|^2) = 0.$$
(2.10)

With regards to $\beta$ and $\mathfrak{b}$: If the operator $\nabla_{\hat{A}t} - 2\alpha$ and $\nabla_{\hat{A}1} + i\nabla_{\hat{A}2}$ commute when acting on sections of $\mathcal{L}^+$, then the equation in Item b) of the top bullet in (2.6) will hold on any open set in $(0, \infty) \times \mathbb{R}^2$ where $\beta$ and $\mathfrak{b}$ can be written using a section $q$ of $\mathcal{L}^+$ as

$$\beta = (\nabla_{\hat{A}t} - 2\alpha)\, q \quad and \quad \mathfrak{b} = -i(\nabla_{\hat{A}1} + i\nabla_{\hat{A}2})\, q.$$
(2.11)

In this event, Item b) of the third bullet in (2.6) when written in terms of $q$ (with (2.6) used to write $E_{\hat{A}1}$, $E_{\hat{A}2}$ and $B_{\hat{A}3}$) says that

$$-(\nabla_{\hat{A}t}^2 + \nabla_{\hat{A}1}^2 + \nabla_{\hat{A}2}^2)\, q + (4\alpha^2 + 2|\varphi|^2 + 8|\beta|^2 + 8|\mathfrak{b}|^2)\, q = 0.$$
(2.12)

What with (2.7), (2.8) and (2.11), the equations in (2.10) and (2.12) form a closed system for the pair (w, $q$).

The logic in the preceding section will be turned around in Section 3 to use an $\mathbb{R}$-valued function and a $\mathbb{C}$-valued function (both defined on a given open set in $(0,\infty) \times \mathbb{R}^2 \times \{0\}$) to obtain a solution to (2.6) on that same open set. To elaborate: The real valued function will be called w and the $\mathbb{C}$-valued function will be called $q_\ddagger$. A solution to (2.4) can be constructed from w and $q_\ddagger$ if these obey

- $-(\partial_t^2 + \partial_1^2 + \partial_2^2)\, w + e^{4w} + 4 e^{4w}(|\partial_t(e^{-2w} q_\ddagger)|^2 + |(\partial_1 + i\partial_2)(e^{-2w} q_\ddagger)|^2) = 0$,
- $(\partial_t(e^{4w}\partial_t(e^{-2w} q_\ddagger))) + (\partial_1 - i\partial_2)(e^{4w}(\partial_1 + i\partial_2)(e^{-2w} q_\ddagger)) = 0$.

(2.13)

To this end, fix a unit length section of $ad(P_0)$ over the open set to play the role of $\sigma$. Having done that, fix a length $\sqrt{2}$ section of the corresponding version of $\mathcal{L}^+$ to play the role of $\frac{\varphi}{|\varphi|}$. Then set $\varphi$ to be $e^{2w}\frac{\varphi}{|\varphi|}$. Meanwhile, define the function $\alpha$ and the connection $\hat{A}$ by taking $\alpha = \partial_t w$ and using (2.8) for $\hat{A}$. Finally, define $\beta$ and $b$ via (2.11) with $q \equiv q_{\pm} \frac{\varphi}{|\varphi|}$.

To elaborate on (2.13): When written using $\hat{A}, \alpha, \beta$ and $b$, the left hand side of the top bullet equation is the equation in Item a) of the third bullet in (2.6). Meanwhile, the left hand side of the lower bullet equation in (2.13) when written using $\hat{A}, \alpha, \beta$ and $b$ is the equation in Item b) of the third bullet in (2.6). The other equations in (2.6) follow automatically from the definition of $\hat{A}, \alpha, \beta$ and $b$ in terms of w and $q$.

A parenthetical remark: An equivalent point of view is to introduce the $\mathbb{C}$ valued function $Q = e^{-2w} q_{\pm}$ and then write $\beta$ and $b$ in terms of w and Q:

- $\beta = e^{2w} \partial_t Q$ and $b = e^{2w}(\partial_1 + i\partial_2)Q$.
- $-(\partial_t^2 + \partial_1^2 + \partial_2^2)w + e^{4w}(1 + 4(|\partial_t Q|^2 + |(\partial_1 + i\partial_2)Q|^2)) = 0$.
- $\partial_t(e^{4w}\partial_t Q) + (\partial_1 - i\partial_2)(e^{4w}(\partial_1 + i\partial_2)Q) = 0$.

(2.14)

This is equivalent to writing $q$ in (2.11) and (2.12) as $Q\varphi$.

An important point now concerns an ambiguity for the section $q$ that appears in (2.11) and (2.12): Supposing that w and $q$ are defined on some open set $U \subset (0, \infty) \times \mathbb{R}^2$, view $\mathbb{R}^2$ as $\mathbb{C}$ and let $h$ denote the restriction to U of a t-independent, holomorphic function of the complex coordinate $z = z_1 + iz_2$ that is defined for all z with $(t,z) \in U$. Then, the change

$$q \to q + h\varphi$$

(2.15)

does not change $\beta$ nor does it change $b$. (This is equivalent to adding $h$ to the function Q that appears in (2.14).) This ambiguity is important because the versions of $q$ for the solutions in Theorem 1 will have a meromorphic pole along the $z = 0$ locus. Even so, this pole will be removable by just such a transformation. In effect, $h$ defines a Čech cocycle for the intersection of two open sets in $(0, \infty) \times \mathbb{C}$ with one being the set where z is non-zero and the other being a neighborhood of the $z = 0$ locus.

### c) Sl(2;$\mathbb{C}$) point of view

As noted in the introduction, Giaotto-Witten [GT], [MW] and [HM1-3] rewrote the top three equations (1.3) as integrability equations for a complex vector bundle structure. As explained below, the rewriting of (1.3) in the previous sections is a disguised form of the complex vector bundle rewriting in [GT], [MW] and [HM1-3].

To start the explanation, introduce by way of notation $\lambda$ to denote $\sigma - 4q$. The identities in (2.5) and (2.11) say in effect that $\lambda$ obeys the identities

$$\nabla_{At}\lambda - i[\mathfrak{a}_3, \lambda] = 0 \quad and \quad (\nabla_{A1} + i\nabla_{A2})\lambda = 0 .$$

(2.16)

Thus, both $\lambda$ and $\varphi$ are annihilated by the two operators $\nabla_{At} - i[\mathfrak{a}_3, \cdot]$ and $\nabla_{A1} + i\nabla_{A2}$. Since $\lambda$ and $\varphi$ are linearly independent in $\text{ad}(P_0) \otimes_{\mathbb{R}} \mathbb{C}$ where $z \neq 0$, this pair uniquely determine A and $\mathfrak{a}_3$. To see this explicity, note that the pair $(\varphi, \lambda)$ also obey the algebraic identities

$$\langle \lambda^2 \rangle = 1 \quad and \quad \langle \varphi^2 \rangle = 0 \quad and \quad \tfrac{i}{2}[\lambda, \varphi] = \varphi$$

(2.17)

at each point. Meanwhile, the pair $(\sigma, \eta \equiv \tfrac{\varphi}{|\varphi|}\hat{\imath})$ obey the analogous three algebraic constraints, $\langle \sigma^2 \rangle = 1$, $\langle \eta^2 \rangle = 0$ and $\tfrac{i}{2}[\sigma, \eta] = 0$. Since any pair in $\text{sl}(2; \mathbb{C})$ obeying these three constraints is conjugate to $(\sigma, \eta)$ via the action of an $\text{Sl}(2; \mathbb{C})$ matrix, it follows that there exists a unique element in $P \times_{\text{Ad}(SU(2))} \text{Sl}(2; \mathbb{C})$ over any given $z \neq 0$ point (call it $\mathfrak{g}$ for now) such that

$$\mathfrak{g} \sigma \mathfrak{g}^{-1} = \lambda \quad and \quad \mathfrak{g} \eta \mathfrak{g}^{-1} = \varphi .$$

(2.18)

A computation finds that $\mathfrak{g}$ in the case at hand can be written using the data $\sigma$, w and $q$ as

$$\mathfrak{g} = \cosh w + i \sinh w \, \sigma + 2i e^{-w} q .$$

(2.19)

As a parenthetical remark for now but a crucial input for Section 4 and subsequently: Because $\nabla_{\theta_\ddagger}$ annihilates both $\sigma$ and $\eta$, the connection A can be obtained by conjugating $\nabla_{\theta_\ddagger}$ with $\mathfrak{g}$. To do this, note first that the operator $\mathfrak{g}^{-1} \nabla_{\theta_\ddagger}(\mathfrak{g}(\cdot))$ when acting on sections of $\text{ad}(P_0) \otimes_{\mathbb{R}} \mathbb{C}$ will annihilate both $\lambda$ and $\varphi$ because of (2.18). Since (2.16) is obeyed, and since the operator $\nabla_A$ is $\nabla_{\theta_\ddagger} + [\mathfrak{b}, \cdot]$ whereas the operator $\mathfrak{g}^{-1} \nabla_{\theta_\ddagger}(\mathfrak{g}(\cdot))$ is $\nabla_{\theta_\ddagger} + \mathfrak{g}^{-1} \nabla_{\theta_\ddagger} \mathfrak{g}$, it follows that when A is written as $\theta_\ddagger + \mathrm{A}$, then $\mathfrak{a}_3$ and the $\text{ad}(P_0)$ valued 1-form $\mathrm{A}$ are given in terms of $\mathfrak{g}$ using the following rules:

- $A_t - i \mathfrak{a}_3 = -(\nabla_{\theta_\ddagger t} \mathfrak{g}) \mathfrak{g}^{-1}$ .
- $A_1 + i A_2 = -((\nabla_{\theta_\ddagger 1} \mathfrak{g}) \mathfrak{g}^{-1} + i (\nabla_{\theta_\ddagger 2} \mathfrak{g}) \mathfrak{g}^{-1})$ .

(2.20)

This gives another view regarding the fact that $\sigma$, w and $q$ determine $\mathfrak{a}_3$ and A.

Crucially for what is to comes after Section 3, there is a generalization of (2.20) along the following lines: Suppose that $(A_\diamond, \varphi_\diamond, \mathfrak{a}_{\diamond 3})$ is any given data set that obeys the first three bullets in (1.3), and suppose that $\mathfrak{g}$ is any given section of $P \times_{Ad(SU(2))} SL(2;\mathbb{C})$. Then the data set $(A, \varphi, \mathfrak{a}_3)$ that is obtained from $(A_\diamond, \varphi_\diamond, \mathfrak{a}_{\diamond 3})$ and $\mathfrak{g}$ by the rules that follow in (2.21) also obeys the first three bullets of (1.3).

- $(A - A_\diamond)_t - i(\mathfrak{a}_3 - \mathfrak{a}_{\diamond 3}) = -(\nabla_{A_{\diamond t}} \mathfrak{g}) \mathfrak{g}^{-1} + i[\mathfrak{a}_{\diamond 3}, \mathfrak{g}]\mathfrak{g}^{-1}$.
- $(A - A_\diamond)_1 + i(A - A_\diamond)_2 = -(\nabla_{A_{\diamond 1}} + i\nabla_{A_{\diamond 2}}) \mathfrak{g} \mathfrak{g}^{-1}$.
- $\varphi = \mathfrak{g} \varphi_\diamond \mathfrak{g}^{-1}$.

(2.21)

This observation and (2.21) are recapitulated at the start of Section 4 because they are the basis for what is done there and in Sections 5 and 6.

### 3. An initial choice for (A, $\mathfrak{a}$)

The purpose of this subsection is to present a family of $\mathbb{R}$-invariant pairs $(A, \mathfrak{a})$ that obey the asymptotic conditions in (1.9) and (1.10) for a given non-negative integer $m$ and positive integer p with each member of this family obeying the top three bullet equations in (1.3). The subsequent sections will explain how to deform each member of this family so as to obtain a solution to all of (1.3) plus (1.9) and (1.10).

Looking ahead, the family is parametrized by a set of p complex numbers $(a_1, \ldots, a_p)$ with $a_p \neq 0$ and a nuisance parameter $\delta \in (0,1)$. Having chosen this data, the subsections that follow describe the corresponding version of $(A, \mathfrak{a})$. By way of notation for these subsections, the integer $m+2p$ is denoted below by $k$. (Thus, $p = \frac{1}{2}(k-m)$.)

### a) The initial choice of (A, $\mathfrak{a}$) where t very small and $|z| < t$

The definition of A and $\mathfrak{a}_3$ for t small and $|z| \leq t$ is straightforward enough:

$$A = A^{(k)} \quad and \quad \mathfrak{a}_3 = \mathfrak{a}^{(k)}_3.$$

(3.1)

This definition guarantees a priori the top bullet equation in (1.3). The definition of $\varphi$ is less straightforward because the second third bullets in (1.3) will be obeyed plus the second bullet in (1.10). In addition, $\varphi$ will vanish only on the $z = 0$ locus and it will vanish there with degree $m$. In this regard: Although $\varphi = z^{m-k}\varphi^{(k)}$ obeys the second and third bullets in (1.3) and vanishes only at $z = 0$ and with degree $m$, it is not viable because it runs afoul of the second bullet in (1.10).

Some background for defining $\varphi$: A length 1 section of the bundle ad(P) can be used to decompose $ad(P) \otimes_\mathbb{R} \mathbb{C}$ as a direct sum of eigenbundles for the section's commutator action: Letting $\sigma$ denote the section in question, the decomposition is given in (2.1). The initial choice of (A, 𝔞) that is described momentarily will have a corresponding unit length section σ for which (2.2) holds; so that φ is a section of the corresponding $\mathcal{L}^+$. And, as such, φ will have an order $m$ zero on the $z = 0$ locus with $m$ being the integer that labels the $t \to \infty$ model solution that appears in (1.10). (The $z = 0$ locus is its only zero.)

Preliminary definitions are needed to obtain an appropriate version of σ and φ. To this end, suppose henceforth that $(a_1, \ldots, a_p) \in \mathbb{C}^{m+p}$ has been chosen (with $a_p \neq 0$). Introduce $\mathfrak{P}$ to denote the polynomal in the coordinate z given by

$$\mathfrak{P} = \sum_{j=1}^{p} a_j z^{k-j}$$

(3.2)

Let ς denote the function

$$\varsigma \equiv \tfrac{1}{2} \frac{(k+1)\sinh^{k+1}(\Theta)}{\sinh((k+1)\Theta)} \frac{z^k}{t^{k+1}},$$

(3.3)

This function ς appeared before in (1.5)'s definition of $\varphi^{(k)}$; the latter being $-\varsigma(\sigma_1 - i\sigma_2)$.

With $\mathfrak{P}$ and ς in hand, define γ by the rule

$$\gamma = \frac{\mathfrak{P}}{\varsigma}$$

(3.4)

A note-worthy feature of γ is its order p pole at $z = 0$:

$$\gamma \sim \frac{2}{k+1} t^{k+1} \sum_{j=1}^{p} \frac{a_j}{z^j} \quad \textit{for z very near } 0.$$

(3.5)

The function γ is used next to define a length 1 section of the bundle $\text{ad}(P_0)$ over the complement of the $z = 0$ locus in $(0, \infty) \times \mathbb{R}^2 \times \mathbb{R}$. The section is denoted by σ and here is its definition:

$$\sigma = \tfrac{1}{1+|\gamma|^2} (\overline{\gamma}(\sigma_1 - i\sigma_2) + (1 - |\gamma|^2)\sigma_3 + \gamma(\sigma_1 + i\sigma_2)) .$$

(3.6)

The polar behavior in (3.5) notwithstanding, this section σ extends over the $z = 0$ locus to define a smooth section of $\text{ad}(P_0)$ on the whole of $(0, \infty) \times \mathbb{R}^2 \times \mathbb{R}$. The corresponding version of the line bundle $\mathcal{L}^+$ is the span (where $z \neq 0$) of its section

$$\hat{\imath} \equiv \tfrac{1}{1+|\gamma|^2}((\sigma_1 - i\sigma_2) - 2\gamma\sigma_3 - \gamma^2(\sigma_1 + i\sigma_2)) .$$

(3.7)

Although its norm is equal to $\sqrt{2}$, it doesn't extend continuously across the $z = 0$ locus.

Granted the preceding, set $\varphi$ to be $-\varsigma(1+|\gamma|^2)\hat{\imath}$, thus

$$\varphi = -\varsigma((\sigma_1 - i\sigma_2) - 2\gamma\sigma_3 - \gamma^2(\sigma_1 + i\sigma_2)) \ . \tag{3.8}$$

Although $\gamma$ has a pole at $z = 0$, the section $\varphi$ is smooth across the $z = 0$ locus because of the behavior of $\varsigma$ near $z = 0$. In particular, where both $|z|$ is much less than $t$ and $|\gamma|$ is much greater than 1, the section $\varphi$ looks like

$$\varphi = \tfrac{2}{k+1} a_p^2 t^{k+1} z^m (\sigma_1 + i\sigma_2) + \mathcal{O}(|z|^{m+1}) \tag{3.9}$$

which has the same order of vanishing as $\varphi^{(m)}$. On the other hand, where $|\gamma|$ is much less than 1, the section $\varphi$ differs little from $\varphi^{(k)}$.

The next lemma makes a formal assertion with regards to these limits, and it asserts that the second and third bullets of (1.3) are obeyed.

**Lemma 3.1:** *Having fixed a p-tuple of complex numbers $(a_1, \ldots, a_p)$ with $a_p \neq 0$, there exists $\kappa > 1$ with the following significance: Define $\gamma$ via (2.1)–(2.3) and then $\varphi$ via (2.7).*

- *If $t < \tfrac{1}{\kappa}$ and if $|z| \leq 4t$, then the function $\gamma$ obeys*

  a) $\gamma = \tfrac{2}{k+1}(1+\mathfrak{e}_1)t^{k+1} \sum_{j=1}^{p} \tfrac{a_j}{z^j}$ *with the norm of $\mathfrak{e}_1$ obeying* $|\mathfrak{e}_1| < \kappa \tfrac{|z|^2}{t^2}$.

  b) $\kappa^{-1} \tfrac{t^{k+1}}{|z|^p} < |\gamma| < \kappa \tfrac{t^{k+1}}{|z|^p}$ .

- *If $t < \tfrac{1}{\kappa}$ and if $|z| < 4t$, then $\varphi$ obeys*

  a) $\varphi = \varphi^{(k)} + \mathfrak{e}_2$ *with* $|\mathfrak{e}_2| \leq \kappa |\varphi^{(k)}| \tfrac{t^{k+1}}{|z|^p}$ .

  b) $\varphi = \tfrac{2}{k+1} a_p^2 t^{k+1} z^m (\sigma_1 + i\sigma_2) + \mathfrak{e}_3$ *with* $|\mathfrak{e}_3| \leq \kappa(|z|^{m+p} + \tfrac{|z|^k}{t^{k+1}})$

- $\nabla_{At}\varphi - [i\mathfrak{a}_3, \varphi] = 0$ *and* $(\nabla_{A1} + i\nabla_{A2})\varphi = 0$.

*Proof of Lemma 3.1*: If $|z| \leq c_0^{-1} t$, then $\varsigma$ can be written as

$$\varsigma = (1+\mathfrak{r}_1)\tfrac{(k+1)}{2}\tfrac{z^k}{t^{k+1}} \tag{3.10}$$

with $\mathfrak{r}_1$'s norm obeying $|\mathfrak{r}_1| \leq c_0 \tfrac{|z|^2}{t^2}$. In general, the norm of $\varsigma$ obeys

$$c_0^{-1} \tfrac{1}{t} \tfrac{|z|^k}{(t^k + |z|^k)} \leq |\varsigma| \leq \tfrac{1}{\sqrt{2}\,t} \tag{3.11}$$

These two bounds imply what is asserted by Item a) of the top bullet of the lemma. Meanwhile, if $t < c_0^{-1}$, then $\mathfrak{P}$ can be written as

$$\mathfrak{P} = (1 + \mathfrak{r}_3) \, a_p z^{k-p}$$

(3.12)

with $\mathfrak{r}_3$'s norm obeying $|\mathfrak{r}_3| \le c_0 |z|$. (This is because $a_p \ne 0$.) This last bound with (3.10) and (3.11) lead to the upper and lower bounds in Item b) of the top bullet of the lemma. The assertions made by the second bullet of the lemma follow directly from (3.11) using the top bullet of the lemma.

The third bullet of the lemma follows from three observations. The first is that $\varphi$ can be written as

$$\varphi = \varphi^{(k)} - 2\mathfrak{P}\sigma_3 - \mathfrak{P}^2 \frac{\varphi^{(k)*}}{|\varphi^{(k)}|^2} \ .$$

(3.13)

The second observation is that $\mathfrak{P}$ is t-independent and annihilated by $\partial_1 + i\partial_2$. The third is that $A = A^{(k)}$ and $\mathfrak{a}_3 = \mathfrak{a}_3^{(k)}$, and that the $A = A^{(k)}$ and $\mathfrak{a}_3 = \mathfrak{a}^{(k)}{}_3$ versions of the operators $\nabla_{At} - [i\mathfrak{a}_3, \cdot]$ and $\nabla_{A1} + i\nabla_{A2}$ annihilate both $\varphi^{(k)}$ and $\frac{\varphi^{(k)*}}{|\varphi^{(k)}|^2}$.

**b) The description of $(A, \mathfrak{a})$ where t is small and $|z| \ge t$**

Introduce by way of notation $\chi_*$ to denote the function on $(0, \infty) \times \mathbb{R}^2$ given by the rule $\chi_*(t, z) = \chi(\frac{|z|}{t} - 1))$, it being zero where $|z| > 2t$ and equal to one where $|z| < t$.

Given the set $(a_1, \ldots, a_p)$ of complex numbers subject to the constraint that $a_p \ne 0$, reintroduce the complex function $\gamma$ from (3.4) and introduce, by way of notation, $\gamma_*$ to denote $\chi_* \gamma$. Granted this notation, henceforth define $\sigma$ and $\hat{\imath}$ by replacing $\gamma$ with $\gamma_*$ in (3.6) and (3.7). This is to say that from now on,

- $\sigma \equiv \frac{1}{1+|\gamma_*|^2} (\bar{\gamma}_*(\sigma_1 - i\sigma_2) + (1 - |\gamma_*|^2)\sigma_3 + \gamma_*(\sigma_1 + i\sigma_2))$ .
- $\hat{\imath} = \frac{1}{1+|\gamma_*|^2} ((\sigma_1 - i\sigma_2) - 2\gamma_*\sigma_3 - \gamma_*^2(\sigma_1 + i\sigma_2))$ .

(3.14)

Letting $\kappa_*$ denote the version of $\kappa$ from the $\{a_1, \ldots, a_p\}$ version of Lemma 3.1, define $\varphi$ where $|z| > t$ and $t < \frac{1}{\kappa_*}$ to be $-\varsigma(1 + |\gamma_*|^2)\hat{\imath}$, thus

$$\varphi \equiv -\varsigma\big((\sigma_1 - i\sigma_2) - 2\gamma_*\sigma_3 - \gamma_*^2(\sigma_1 + i\sigma_2)\big)$$

(3.15)

which is the formula in (3.8) with $\gamma$ replaced by $\chi_*\gamma$.

The first step towards defining $A$ and $\mathfrak{a}_3$ where $|z| \ge t$ and $t < \frac{1}{\kappa_*}$ is to depict the latter using $(\alpha, \beta, \hat{A}, \hat{b})$ as per the discussion in Section 2b. To do that, first write the norm of $|\varphi|$ from (3.15) as $e^{2w}$. Then introduce $\theta_{\ddot{+}}$ to denote the unique flat connection on $(0, \infty) \times (\mathbb{R}^2 - \{0\})$ whose covariant derivative annihilates both $\sigma$ and $\frac{\varsigma}{|\varsigma|}\hat{\imath}$. The corresponding function $\alpha$ and the connection $\hat{A}$ are then given where $t < \frac{1}{\kappa_*}$ in terms of $w$ and $\theta_{\ddot{+}}$ by (2.7) and (2.8). To obtain $\beta$ and $\mathfrak{b}$: These are defined where $|z| > 0$ and $t < \frac{1}{\kappa_*}$ via (2.11) using the following section $q$ of $\mathcal{L}^+$ that is defined where $z \neq 0$ by the rule

$$q = \tfrac{1}{4} \overline{\gamma}_* \hat{\imath} \; .$$
(3.16)

A computation using specific properties of $\varsigma$ leads to formulas for $\beta$ and $\mathfrak{b}$ which are given below (the notation has $w^{(k)}$ denoting $\tfrac{1}{2} \ln |\varphi^{(k)}|$.):

- $\beta = -\tfrac{1}{4} \tfrac{1}{1+|\gamma_*|^2} (4\partial_t w^{(k)} \overline{\gamma}_* - \partial_t \chi_* (1 - |\gamma_*|^2) \overline{\gamma}) \hat{\imath}$.
- $\mathfrak{b} = -\tfrac{i}{4} \tfrac{1}{1+|\gamma_*|^2} ((\partial_1 + i\partial_2)\overline{\gamma}_* - 2((\partial_1 + i\partial_2) w^{(k)}) \overline{\gamma}_* - ((\partial_1 + i\partial_2)\chi_*) \overline{\gamma} |\gamma_*|^2) \hat{\imath}$.

(3.17)

With regards to $q$ (for later reference): It can be written near the $z = 0$ locus (where $\gamma_* = \gamma$) as

$$q = \tfrac{1}{4} \tfrac{1}{\mathfrak{P}} \varphi - \tfrac{1}{4} \tfrac{1}{\overline{\gamma}} \hat{\imath} \; .$$
(3.18)

with $\mathfrak{P}$ being the holomorphic polynomial from (3.2). Notwithstanding the fact that $\tfrac{1}{\mathfrak{P}} \varphi$ is singular along the $z = 0$ locus (its norm where $|z| \ll t$ is a positive multiple of $\tfrac{1}{|z|^p}$), the $\tfrac{1}{\mathfrak{P}} \varphi$ term in (3.18) makes no contribution to either $\beta$ or $\mathfrak{b}$ because $\mathfrak{P}$ is holomorphic and because $\varphi$ is annihilated by both $\nabla_{\hat{A}t} - 2\alpha$ and $\nabla_{\hat{A}1} + i\nabla_{\hat{A}2}$. (This is a manifestation of what is said regarding (2.15).) Meanwhile, the term $\tfrac{1}{\overline{\gamma}}\hat{\imath}$ in (3.18) extends smoothly across the $z = 0$ locus which is why $\beta$ and $\mathfrak{b}$ are actually smooth across this locus. (The norm of $\tfrac{1}{\overline{\gamma}}\hat{\imath}$ near $z = 0$ is bounded by $c_0 \tfrac{|z|^p}{t^{k+1}}$.)

The following lemma summarizes the salient properties of this version of $(A, \mathfrak{a})$.

**Lemma 3.2**: *Having fixed complex numbers* $(a_1, \ldots, a_p)$ *with* $a_p \neq 0$, *there exists* $\kappa > \kappa_*$ *with the following significance: Define* $(A, \mathfrak{a})$ *where* $t < \tfrac{1}{\kappa}$ *as done above.*
- *Where* $|z| < t$, *the pair* $(A, \mathfrak{a}_3)$ *is equal to* $(A^{(k)}, \mathfrak{a}^{(k)}{}_3)$, *and $\varphi$ is given by* (3.8)
- *Where* $t \le |z| < 4t$, *the pair* $(A, \mathfrak{a})$ *obeys*
  a) $|A - A^{(k)}| + |\mathfrak{a} - \mathfrak{a}^{(k)}| \le \kappa t^{k-p}$.
  b) $|\nabla_{A^{(k)}} (A - A^{(k)})| + |\nabla_{A^{(k)}} (\mathfrak{a} - \mathfrak{a}^{(k)})| \le \kappa t^{k-p-1}$
- *Where* $|z| > 4t$, *the pair* $(A, \mathfrak{a})$ *is equal to* $(A^{(k)}, \mathfrak{a}^{(k)})$.

*Proof of Lemma 3.2*:  There are three parts to the proof, one for each bullet of the lemma. The convention in this proof is that $c_0$ can depend on the data $(a_1, \ldots, a_p)$.

*Part 1*:  The key observation with regards to the top bullet is that $\gamma_* = \gamma$ where $|z| < t$.  As a direct consequence, the formulas for $\varphi$ in (3.8) and (3.15) are identical.  With regards to $\mathfrak{a}_3$:  Since $|\varphi| = |\varphi^{(k)}|(1 + |\gamma|^2)$, the function w is $w^{(k)} + \frac{1}{2}\ln(1 + |\gamma|^2)$.  Since the t-derivative of $\varsigma$ is $2\alpha^{(k)}\varsigma$, the t-derivative of $\gamma$ is $-2\alpha^{(k)}\gamma$.  This implies that

$$\alpha = \alpha^{(k)} \frac{1 - |\gamma|^2}{1 + |\gamma|^2} \ .$$

(3.19)

Because $\beta = -\frac{1}{2}\alpha^{(k)} \frac{\bar{\gamma}}{1+|\gamma|^2} \hat{\imath}$, it follows from (2.3) that $\alpha\sigma + \beta + \beta^*$ is $\alpha^{(k)}\sigma_3$ which is $\mathfrak{a}^{(k)}_3$. The proof that $A = A^{(k)}$ is along similar lines and is left to the reader.

*Part 2*:  With regards to the second bullet of the lemma:  The first key point here is the bound from Item b) of the first bullet of Lemma 3.1 to the effect that $|\gamma| \leq c_0 \, t^{k-p+1}$ where $|z|$ is between $t$ and $4t$ if $t < c_0^{-1}$  This implies that $|\sigma - \sigma_3| \leq c_0 t^{k-p+1}$ and likewise that $|\hat{\imath} - (\sigma_1 - i\sigma_2)| \leq c_0 t^{k-p+1}$ where $|z|$ is between $t$ and $4t$.  The latter lead directly to the bound $|\varphi - \varphi^{(k)}| \leq c_0 t^{k-p}$ where $|z|$ is between $t$ and $4t$.  The second key point is that $|\nabla\gamma| \leq c_0 t^{k-p}$ where $|z|$ is between $t$ and $4t$.  This implies that $|\alpha - \alpha^{(k)}| \leq c_0 t^{k-p}$ where $|z|$ is between $t$ and $4t$ because w is equal to $w^{(k)} + \frac{1}{2}\ln(1 + |\gamma|^2)$.  It also implies directly that both $|\beta|$ and $|\mathfrak{b}|$ are bounded by $c_0 t^{k-p}$ where $|z|$ is between $t$ and $4t$.  The bounds for the norms of $\alpha - \alpha^{(k)}$ and $\beta$ and $\sigma - \sigma_3$ imply that $|\mathfrak{a}_3 - \mathfrak{a}^{(k)}_3| \leq c_0 t^{k-p}$ where $|z|$ is between $t$ and $4t$.  Meanwhile, the bounds for $|\mathfrak{b}|$ and the formula for w imply that $|A - A^{(k)}|$ is also bounded by $c_0 t^{k-p}$ where $|z|$ is between $t$ and $4t$.  Similar arguments lead to the asserted bounds on the $A^{(k)}$-covariant derivatives of $A - A^{(k)}$ and $\mathfrak{a}_3 - \mathfrak{a}_3^{(k)}$.

*Part 3*:  With regards to the third bullet:  This follows directly from the fact that $\gamma_*$ is zero where $|z| > 2t$.

### c) The description of $(A, \mathfrak{a})$ where t is bounded away from zero

Let $\kappa_\diamond$ denote the version of $\kappa$ from Lemma 3.2.  This subsection describes the pair $(A, \mathfrak{a})$ where t is bounded away from zero.  This task has four parts. By way of notation:  The complex functions $\varsigma$ and $\gamma$ are still defined via (3.3); and $\gamma_*$ is still $\chi_*\gamma$.  Likewise, $\sigma$ and $\hat{\imath}$ are still given by the formulas in (3.14).  One addition is needed to these definitions which requires choosing a positive number $\delta$ with upper bound $\frac{1}{100\kappa_\diamond}$.  This number is used to define the purely t-dependent cut-off function $\varpi_\delta$ by the rule whereby $\varpi_\delta(t) = \chi(\frac{2t}{\delta} - 1)$. This function is equal to 1 where $t < \frac{1}{2}\delta$ and it is equal to 0 where $t > \delta$.

*Part 1*: The section $\varphi$ where t is bounded away from zero is written as

$$\varphi = \tfrac{1}{\sqrt{2}} e^{2w} \tfrac{\varsigma}{|\varsigma|} \hat{\imath}$$

(3.20)

with w given by the rule

$$w = \tfrac{1}{2}\varpi_\delta \ln\big((1+|\gamma_*|^2)\tfrac{(k+1)\sinh(\Theta)}{\sinh((k+1)\Theta)}\big) + \tfrac{1}{2}(1-\varpi_\delta)\ln\big(\tfrac{(m+1)\sinh(\Theta)}{\sinh((m+1)\Theta)}\big) - \tfrac{1}{2}\ln(\sqrt{2}\,t)\,.$$

(3.21)

(The norm of $\varphi$ is $e^{2w}$.) This definition has $\varphi$ equal to its namesake from Section 3b where $t < \tfrac{1}{2}\delta$. The definition also has the norm of $\varphi$ where $t > \delta$ being equal to the norm of $\varphi^{(m)}$. (But $\varphi$ is not equal to $\varphi^{(m)}$ because $\hat{\imath}$ is not $(\sigma_1 - i\sigma_2)$.) Another point to note: The section $\varphi$ is smooth across the $z = 0$ locus; it looks very near this locus like a non-zero multiple of $z^m(\sigma_1 + i\sigma_2) + \mathcal{O}(|z|^{m+1})$. (In this regard: The line bundle $\mathcal{L}^+$ is the span of $\sigma_1 + i\sigma_2$ on the $z = 0$ locus because the section $\sigma$ on the $z = 0$ locus is equal to $-\sigma^3$.)

*Part 2*: It remains now to define A and $\mathfrak{a}_3$ where $t \geq \tfrac{1}{2}\delta$. This is done by writing A and $\mathfrak{a}_3$ in terms of a data set $(\hat{A}, \alpha, \beta, \mathfrak{b})$ as described in Section 2. In particular, the function $\alpha$ and the connection $\hat{A}$ are given in (2.7) and (2.8) in terms of w which is depicted in (3.21). Note in this regard that the formula for w where $t > \delta$ has $\varpi_\delta \equiv 0$ and so $w = w^{(m)}$. Meanwhile, $\theta_\ddagger$ where $t > \delta$ is the unique flat connection on $P_0$ whose covariant derivative annihilates both $\sigma$ and $\tfrac{\varphi}{|\varphi|}$ (the latter being $\tfrac{\varsigma}{|\varsigma|}\hat{\imath}$). Therefore, the definitions in (2.7) and (2.8) say in effect that the $t > \delta$ versions of $\hat{A}$ and $\alpha$ are

- $\hat{A} = \theta_\ddagger + \tfrac{1}{2}\big(1 - (m+1)\tfrac{\sinh(\Theta)}{\cosh(\Theta)}\tfrac{\cosh((m+1)\Theta)}{\sinh((m+1)\Theta)}\big)\tfrac{1}{|z|^2}(z_1 dz_2 - z_2 dz_1)\,\sigma$.
- $\alpha = -\tfrac{1}{2t}\tfrac{(m+1)\sinh(\Theta)}{\sinh((m+1)\Theta)}\tfrac{\cosh((m+1)\Theta)}{\cosh(\Theta)}$.

(3.22)

A key obsrvation: This $t > \delta$ definition of the pair $(\hat{A}, \alpha\sigma)$ is the pull-back of the model pair $(A^{(m)}, \mathfrak{a}^{(m)})$ on the $t > \delta$ part of $(0, \infty) \times \mathbb{R}^2 \times \mathbb{R}$ by an automorphism of $P_0$.

*Part 3*: This part defines the $\beta$ and $\mathfrak{b}$ parts of the connection A (remember that $\beta$ is part of $\mathfrak{a}_3$ also) where $t \geq \tfrac{1}{2}\delta$ on $(0, \infty) \times \mathbb{R}^2 \times \mathbb{R}$. The constraints on $\beta$ and $\mathfrak{b}$ are as follows: They are given by (3.17) where $t < \tfrac{1}{2}\delta$, they have to be smooth where $z = 0$, they have to limit to zero as $\tfrac{|z|}{t} \to \infty$, and they have to vanish for $t \to \infty$. In addition, $\beta$ and $\flat$ must obey the equation that is depicted in Item b) of the top bullet of (2.6). That precludes extending the versions of $\beta$ and $\flat$ in (3.17) by multiplying the right hand sides of the latter formulae by $\varpi_\delta$ or some other cut-off function.

In order to guarantee the identity from Item b) of the top bullet in (2.6), the pair $(\beta,\flat)$ will be defined on the whole of $(0,\infty)\times\mathbb{R}^2$ via (2.11) using an appropriate extension of the section $q$ that is depicted in (3.16). But note that taking $q$ to be $\frac{1}{4}\varpi_\delta\,\overline{\gamma}_*\hat{\imath}$ <u>does not</u> work because the resulting version of $\beta$ will be singular at points along the $z = 0$ locus: Its norm will be greater than a positive multiple of $|\partial_t\varpi_\delta|\frac{1}{|z|^p}$ near $z = 0$ because of the $\frac{1}{\mathfrak{P}}\varphi$ term in (3.18).

To define $q$, first write the $\frac{1}{\mathfrak{P}}$ from the $\frac{1}{\mathfrak{P}}\varphi$ term in (3.18) near $z = 0$ as

$$\tfrac{1}{\mathfrak{P}} = \tfrac{1}{\mathfrak{a}_p}\left(\tfrac{1}{z^{m+p}} + \mu_{p-1}\tfrac{1}{z^{m+p-1}} + \cdots + \mu_1\tfrac{1}{z^{m+1}}\right) + \mathfrak{R}$$

(3.23)

with $(\mu_1, \ldots, \mu_{p-1})$ being complex numbers and with $\mathfrak{R}$ being a Laurent series in $z$ with norm bounded near $z = 0$ by $\frac{1}{|z|^m}$. Now write the right hand side of (3.16) as

$$\tfrac{1}{4}\overline{\gamma}_R\hat{\imath} = \tfrac{1}{4}\tfrac{1}{\mathfrak{a}_p}\left(\tfrac{1}{z^{m+p}} + \mu_{p-1}\tfrac{1}{z^{m+p-1}} + \cdots + \mu_1\tfrac{1}{z^{m+1}}\right)\varphi + q\,.$$

(3.24)

Note in particular that what is denoted above by $q$ is bounded across the $z = 0$ locus because $\varphi$ near $z = 0$ has the from $\varphi = x(t)z^m(\sigma_1 + i\sigma_2) + \mathcal{O}(|z|^{m+1})$ with $x$ being a non-zero function of $t$. (The $q$ term is actually smooth across $z = 0$.) With the preceding understood, define $q$ on the $z\neq 0$ and $t > \frac{1}{2}\delta$ part of $(0,\infty)\times\mathbb{R}^2$ by the rule

$$q = \tfrac{1}{4}\tfrac{1}{\mathfrak{a}_p}\left(\tfrac{1}{z^{m+p}} + \mu_{p-1}\tfrac{1}{z^{m+p-1}} + \cdots + \mu_1\tfrac{1}{z^{m+1}}\right)\varphi + \varpi^{\ddagger}_\delta q$$

(3.25)

where $\varpi^{\ddagger}_\delta$ denotes the following:

- $\varpi^{\ddagger}_\delta = \varpi_\delta$ which is $\chi(\frac{2t}{\delta}-1)$ *unless both* $m = 0$ *and* $\mu_1 \neq 0$ *in (3.24)*.
- $\varpi^{\ddagger}_\delta = \chi(\frac{2t}{\delta(1+|z|^2)^{1/4}}-1)$ *if both* $m = 0$ *and* $\mu_1 \neq 0$ *in (3.24)*.

(3.26)

To say something about the distinction in (3.26): The analysis in the Sections 4 and 5 need the square of the norm of $\nabla_{A_t}\mathfrak{a} - B_A - \frac{i}{2}[\varphi,\varphi^*]$ to have finite integral on the whole of $(0,\infty)\times\mathbb{R}^2$. This won't be the case if $\varpi_\delta$ is used for $\varpi^{\ddagger}_\delta$ in the case when $m = 0$ and $\mu_1 \neq 0$.

In either case in (3.26), the definition in (3.25) of $q$ agrees with (3.16) where $t < \frac{1}{2}\delta$. As was the case with (3.18), the singularity of $q$ along the $z = 0$ locus is of no consequence with regards to β and 𝔟: Because of (2.15), both β and 𝔟 are smooth across the z = 0 locus. To be explicit, β and 𝔟 are given where $t > \frac{1}{4}\delta$ by the rule whereby

$$\beta = \partial_t \varpi^{\ddagger}_\delta q + \varpi^{\ddagger}_\delta (\nabla_{\hat{A}t} - 2\alpha) q \quad and \quad \mathit{b} = -\tfrac{i}{4}((\partial_1 + i\partial_2)\varpi^{\ddagger}_\delta q + \varpi^{\ddagger}_\delta(\nabla_{\hat{A}1} + i\nabla_{\hat{A}2}) q \ .$$

(3.27)

An important point with regards to this formula: Supposing that $t > \tfrac{1}{2}\delta$, if $|z| > 2t$, then the complex function $\gamma_*$ that appears in (3.16) and (3.24) is zero in which case

$$q = -\tfrac{1}{4}(\tfrac{1}{a_p}\tfrac{1}{z^{m+p}} + \mu_{p-1}\tfrac{1}{z^{m+p-1}} + \cdots + \mu_1 \tfrac{1}{z^{m+1}})\varphi \ ,$$

(3.28)

As a consequence, β and 𝔟 where $t > \tfrac{1}{2}\delta$ and $|z| > 2t$ are given by

- $\beta = -\tfrac{1}{4}(\partial_t \varpi^{\ddagger}_\delta)(\tfrac{1}{a_p}\tfrac{1}{z^{m+p}} + \mu_{p-1}\tfrac{1}{z^{m+p-1}} + \cdots + \mu_1 \tfrac{1}{z^{m+1}})\varphi \ ,$
- $\mathit{b} = \tfrac{i}{4}((\partial_1 + i\partial_2)\varpi^{\ddagger}_\delta)(\tfrac{1}{a_p}\tfrac{1}{z^{m+p}} + \mu_{p-1}\tfrac{1}{z^{m+p-1}} + \cdots + \mu_1 \tfrac{1}{z^{m+1}})\varphi \ .$

(3.29)

In any event, both β and 𝔟 are zero where $t > \delta$ unless $m = 0$ and $\mu_1 \neq 0$ in (3.24). In the case where $m = 0$ and $\mu_1 \neq 0$, then both are zero where $t > \delta(1+|z|^2)^{1/4}$ and in any event, their support lies where t is between $\tfrac{1}{2}\delta(1+|z|^2)^{1/4}$ and $\delta(1+|z|^2)^{1/4}$ if $t > 4\delta$. This is because (3.29) describes β and 𝔟 in this region (which is the case because the inequalities $|z| < 2t$ and $t < \delta(1+|z|^2)^{1/4}$ are consistent only when $|z|$ is less than $3\delta$ (assuming $\delta < \tfrac{1}{100}$) which requires in turn that t be less than $4\delta$.)

**d) A summary of the properties of (A, 𝔞)**

The following lemma summarizes some of the salient features of the pair (A, 𝔞) just described. By way of notation, the lemma has $\kappa_{\ddagger}$ denoting Lemma 3.2's version of κ.

**Lemma 3.3**: *Fix a collection $(a_1, \ldots, a_p)$ of complex numbers with $a_p \neq 0$. Given these, there exists $\kappa > \kappa_{\ddagger}$ with the following significance: Fix a positive number δ less than $\tfrac{1}{\kappa}$ and use δ with $(a_1, \ldots, a_p)$ to define Section 3a-c's version of A, φ and $\mathfrak{a}_3$. Then the corresponding pair obeys the first three bullets in (1.3) and (1.10).*

*Proof of Lemma 3.3*: The pair (A, 𝔞) obeys (1.3) by construction. With regards to (1.10), the small t behavior is affirmed by Lemma 3.2. With regards to large t and small $\tfrac{t}{|z|}$, it follows from the form of w in (3.21) that it is sufficient to verify that $t|\beta|$ and $t|\mathit{b}|$ limit to zero as $t \to \infty$. With regrds to small $\tfrac{t}{|z|}$, it is sufficient to verify that $|\beta|$ and $|\mathit{b}|$ limit to zero as $\tfrac{t}{|z|} \to 0$. This is done in the next paragraphs.

Suppose first that either $m \neq 0$ or $\mu_1 = 0$: Unless $m = 0$ and $\mu_1 \neq 0$ in (3.24), the definition in (3.26)-(3.27) has $\beta$ and $\flat$ being zero for $t > \delta$. Thus the large t behavior is more than consistent with what is required by (1.10) if either $m \neq 0$ or $\mu_1 = 0$. In this same case, both $\beta$ and $\flat$ limit to zero as $\frac{t}{|z|} \to 0$ because the norms of both are bounded by $c_0 \frac{1}{|z|^2}$ where $|z|$ is greater than $2t$ and $t$ is between $\frac{1}{2}\delta$ and $\delta$ (where it has to be if $\beta$ and $\flat$ are non-zero).

Now suppose that $m = 0$ and $\mu_1 \neq 0$. Where $t > 4\delta$, the pair $(\beta, \flat)$ are depicted by (3.29). Therefore, they vanish where $t \geq \delta(1+|z|^2)^{1/4}$ and are supported where $t$ is between that and half of that. In this region, $|z| \sim t^2$ so $|\beta|$ is bounded by $c_0 |\mu_1| \frac{1}{t^4}$ and $|\flat|$ by bounded by $c_0 |\mu_1| \frac{1}{t^5}$. Indeed, their norms where they are described by (3.29) are bounded by a product of three factors: First, $|\mu_1| \frac{1}{|z|}$ which is $\mathcal{O}(\frac{1}{t^2})$; then a factor of $|\varphi|$ which is $\mathcal{O}(\frac{1}{t})$; and then the norms of first derivatives of $\varpi^{\ddagger}_\delta$ which are bounded by $c_0 \delta \frac{1}{|z|^{1/2}}$ in the case of the t-derivative and $c_0 \delta \frac{1}{|z|}$ for the $z_1$ and $z_2$ derivatives. Thus $|\nabla \varpi^{\ddagger}_\delta|$ is also $\mathcal{O}(\frac{1}{t})$.) As for the small $\frac{t}{|z|}$ behavior where $t > \frac{1}{2}\delta$: This requires that $|z|$ be large; and in this event, it follows from (3.29) that the norms of $\beta$ and $\flat$ are both bounded by $c_0 \delta |\mu_1| \frac{1}{t} \frac{1}{|z|^{3/2}}$ which is $\mathcal{O}(\frac{t}{|z|})^{3/2}$ where $t \geq \frac{1}{2}\delta$.

### e) On the size of $B_{A3}$ and $E_{A1}$ and $E_{A2}$

Theorem 1 refers to the conditions posed in (1.9) of which the third is a constraint on the behavior of $B_{A3}$ and $E_{A1}$ and $E_{A2}$ where $|z|$ has a positive lower bound. The proof that this constraint is obeyed will use the fact that the constraint in the third bullet of (1.9) is satisfied by data sets from Section 3. The next lemma makes a formal assertion to this effect. By way of notation: The lemma uses $\kappa_\diamond$ to denotes Lemma 3.3's version of $\kappa$.

**Lemma 3.4**: *Fix a collection $(a_1, \ldots, a_p)$ of complex numbers with $a_p \neq 0$. Given these, there exists $\kappa > \kappa_\diamond$ with the following significance: Fix a positive number $\delta$ less than $\frac{1}{\kappa}$ and use $\delta$ with $(a_1, \ldots, a_p)$ to define Section 3a-c's version of $A$, $\varphi$ and $\mathfrak{a}_3$. The pointwise norms of the corresponding $B_{A3}$, $E_{A1}$ and $E_{A2}$ are bounded by $\kappa \frac{1}{t^2}$. In addition, for any positive real number $R$,*

$$\int_{(0,\infty) \times \{z \in \mathbb{R}^2 : |z| > R\}} (|B_{A3}|^2 + |E_{A1}|^2 + |E_{A2}|^2) < \frac{\kappa}{R} .$$

***Proof of Lemma 3.3***:  Lemma 3.2 implies the pointwise norm bounds $t < \frac{1}{2}\delta$ and it implies that the part of the $(0, \infty) \times \{z \in \mathbb{R}^2 : |z| > R\}$ integral of $(|B_{A3}|^2 + |E_{A1}|^2 + |E_{A2}|^2)$ where $t < \frac{1}{2}\delta$ is bounded by $\frac{c_0}{R}$.

The pointwise norm bounds also hold where $t > \delta$ unless $m = 0$ and $\mu_1 \neq 0$ because the connection A there is Aut(P) equivalent to the model solution connection $A^{(m)}$. For the same reason, there is a $\frac{c_0}{R}$ bound for the contribution to the $(0, \infty) \times \{z \in \mathbb{R}^2 : |z| > R\}$ integral of $|B_{A3}|^2 + |E_{A1}|^2 + |E_{A2}|^2$ from the $t > \delta$ part of the domain unless $m = 0$ and $\mu_1 \neq 0$.

In the case when $m = 0$ and $\mu_1 \neq 0$, the pointwise norm bounds automatically hold where $t \geq \delta(1 + |z|^2)$ because A and $A^{(m)}$ are Aut(P) equivalent there. Likewise, there is a $\frac{c_0}{R}$ bound when $m = 0$ and $\mu_1 \neq 0$ for the contribution to the $(0, \infty) \times \{z \in \mathbb{R}^2 : |z| > R\}$ integral of $|B_{A3}|^2 + |E_{A1}|^2 + |E_{A2}|^2$ from the $t \geq \delta(1 + |z|^2)^{1/4}$ part of the integration domain.

The verification of the pointwise norm bounds and that there is a $\frac{c_0}{R}$ bound for the contribution to the integral of $|B_{A3}|^2 + |E_{A1}|^2 + |E_{A2}|^2$ from rest of $(0, \infty) \times \{z \in \mathbb{R}^2 : |z| > R\}$ can be done by writing $B_{A3}$, $E_{A1}$ and $E_{A2}$ in terms of the connection $\hat{A}$ and function $\alpha$ given by (3.22) and the pair $\beta$ and $\flat$ as depicted in (3.27). The formulae are given below; the notation uses the superscripts + and - to indicate the $\mathcal{L}^+$ and $\mathcal{L}^-$ projections of the indicated sections of $\text{ad}(P) \otimes_{\mathbb{R}} \mathbb{C}$.

- *The formula for* $B_{A3}$:
  a) $\langle \sigma B_{A3} \rangle = \langle \sigma B_{\hat{A}3} \rangle + 4|\flat|^2$.
  b) $B_{A3}^+ = -i(\nabla_{\hat{A}1} - i\nabla_{\hat{A}2})\flat$.
- *The formula for* $E_{A1} + iE_{A2}$:
  a) $\langle \sigma(E_{A1} + iE_{A2}) \rangle = -i(\frac{\partial}{\partial z_1} + i\frac{\partial}{\partial z_2})\alpha - 4\langle \beta^* \flat \rangle$.
  b) $(E_{A1} + iE_{A2})^+ = i(\nabla_{\hat{A}1} + i\nabla_{\hat{A}2})\beta + 2\nabla_{\hat{A}t}\flat$.
  c) $(E_{A1} + iE_{A2})^- = -i(\nabla_{\hat{A}1} + i\nabla_{\hat{A}2})\beta^*$.

(3.30)

A direct calculation using the formulas given for $(\hat{A}, \alpha, \beta, \flat)$ can now be done to verify all of the required bounds. Readers will be spared the details. What follows directly are the important points that enter.

Since $(\hat{A}, \alpha, \beta, \flat)$ are smooth with uniformly bounded derivatives to second order (any given order) where $t \geq \frac{1}{2}\delta$, the issue is that of the large $|z|$ behavior of $B_{A3}$, $E_{A1}$ and $E_{A2}$. (For the same reason, it is also sufficient to consider only the cases where R is greater than 1 for the $\frac{c_0}{R}$ integral bound.)  In this region, $B_{\hat{A}3}$ is -1 times that Laplacian along the $\mathbb{R}^2$ factor of $(0, \infty) \times \mathbb{R}^2$ of the function w that is depicted in (3.21). In this regard, w in the region of interest is has the form

$$w \sim -\tfrac{1}{2}\ln(\sqrt{2}\,t) + x(t)\frac{t^2}{|z|^2} + \cdots$$

(3.31)

where $x$ is a smooth function of t which is constant if $t < \tfrac{1}{2}\delta$ and if $t > \delta$. Meanwhile, the unwritten terms in (3.31) are $\mathcal{O}(\frac{t^4}{|z|^4})$, the norms of their derivatives are factors of $\mathcal{O}(\tfrac{1}{t})$ or $\mathcal{O}(\tfrac{1}{|z|})$ smaller than this (the former for t-derivatives, the latter for $\mathbb{R}^2$ derivatives), and the norms of the second derivatives are factors of $\mathcal{O}(\tfrac{1}{t^2})$ or $\mathcal{O}(\tfrac{1}{t}\tfrac{1}{|z|})$ or $\mathcal{O}(\tfrac{1}{|z|^2})$ smaller, and so on for higher order derivatives. As a consequence of (3.31), the norm of $|B_{\hat{A}3}|$ where $t \geq \tfrac{1}{2}\delta$ and $|z| > 1$ is bounded by $c_0 \tfrac{t^2}{|z|^4}$. Meanwhile, $\alpha$ is the t-derivative of w, so the norm of its derivatives along the $\mathbb{R}^2$ factor of $(0,\infty)\times\mathbb{R}^2$ is bounded by $c_0\tfrac{t}{|z|^3}$. These pointwise, large $|z|$ bounds for $B_{\hat{A}3}$ and $E_{\hat{A}1}$ and $E_{\hat{A}2}$ are more than sufficient to accommodate a $\tfrac{c_0}{R}$ bound for the lemma's integral.

The contributions of $\beta$ and $\mathfrak{b}$ to the lemma's integral in the region of interest can be bounded using their depictions in (3.29) where $t \geq \tfrac{1}{2}\delta$ and $|z| > 2t$. Suppose first that either $m > 0$ or $\mu_1 = 0$ in which case $\beta$ and $\mathfrak{b}$ are zero unless $t < \delta$. As can be seen from the depiction in (3.29), the norms of products of $\beta$ and $\mathfrak{b}$ in this region are bounded $c_0\tfrac{1}{|z|^4}$, the $\hat{A}$-covariant derivatives along the slices $\{t\}\times\mathbb{R}^2$ are bounded by $c_0 \tfrac{1}{|z|^3}$, and their t direction $\hat{A}$-covariant derivatives are bounded by $c_0\tfrac{1}{|z|^2}$. Thus, the squares of the norms of $B_{A3}$ and $E_{A1}$ and $E_{A2}$ are bounded in the region of interest by $c_0\tfrac{1}{|z|^4}$ whose large $|z|$ fall off is more than sufficient for the lemma's $c_0\tfrac{1}{R}$ bound.

As for the case where $m = 0$ and $\mu_1 \neq 0$: The key points here are that $|\beta| \sim \tfrac{1}{t}\tfrac{1}{|z|^{3/2}}$ and $|\mathfrak{b}| \sim \tfrac{1}{t}\tfrac{1}{|z|^2}$ where $t > \tfrac{1}{2}\delta$ and $|z| \geq 2t$, and that these are supported only where $t \sim |z|^{1/2}$ which is to say where $|z| \sim t^2$. Also, their respective t and $\mathbb{R}^2$ factor $\hat{A}$-covariant derivatives are smaller by factors of $\mathcal{O}(\tfrac{1}{t})$ and $\mathcal{O}(\tfrac{1}{|z|})$ respectively. As a consequence, the $\beta$ and $\mathfrak{b}$ contributions to $|B_{A3}|^2$ and $|E_{A1}|^2$ and $|E_{A2}|^2$ are more than consistent with a $c_0\tfrac{1}{t^2}$ pointwise bound.

The required $c_0\tfrac{1}{R}$ bound for the $|z| > \tfrac{1}{R}$ integral of $|B_{A3}|^2+|E_{A1}|^2+|E_{A2}|^2$ is obtained by integrating the squares of terms that are bounded by $\mathcal{O}(\tfrac{1}{t^2}\tfrac{1}{|z|^2})$ in the worst case. Doing the $|z|$ integration first with the $|z| > R$ constraint leads to a t-integral from $\tfrac{1}{2}\delta$ to infinity of a function that is bounded by $\mathcal{O}(\tfrac{1}{t^4}\tfrac{1}{R})$ when $R > 1$.

### f) On the size of $\nabla_t \mathfrak{a}_3 - B_{A3} - \tfrac{i}{2}[\varphi,\varphi^*]$

The issue in this paper is to find a data set $(A,\varphi,\mathfrak{a}_3)$ that obeys all four bullets in (1.3) plus (1.9) and (1.10). As noted in already, the first three bullets of (1.3) follow directly if the data $\varphi$ and A and $\mathfrak{a}_3$ can be written as in (2.18) and (2.20) in terms of a section of $P \times_{Ad(SU(2))} Sl(2;\mathbb{C})$. The issue then is to find a section so that $\varphi$, A and $\mathfrak{a}_3$ obey the fourth equation in (1.3). To say more about that, introduce by way of notation $\mathfrak{X}$ to denote the section of ad(P) over $(0,\infty) \times \mathbb{R}^2$ that is given by the rule

$$\mathfrak{X} \equiv \nabla_{At}\mathfrak{a}_3 - B_{A3} - \tfrac{i}{2}[\varphi, \varphi^*] \, .$$

(3.32)

The fourth equation in (1.3) is the assertion that $\mathfrak{X} \equiv 0$. This is the significance of $\mathfrak{X}$.

The following lemma gives bounds for the norm of the version of $\mathfrak{X}$ as defined by the data $\varphi$, A and $\mathfrak{a}_3$ from the preceding subsections. The lemma also has $\kappa_\diamond$ denoting the version of $\kappa$ from Lemma 3.3. These bounds play a central role in the rest of this paper.

**Lemma 3.5**: *Fix a collection $(a_1, \ldots, a_p)$ of complex numbers with $a_p \neq 0$. Given these, there exists $\kappa > \kappa_\diamond$ with the following significance: Fix a positive number $\delta$ less than $\tfrac{1}{\kappa}$ and use $\delta$ with $(a_1, \ldots, a_p)$ to define Section 3's version of $\varphi$, $\mathfrak{a}_3$ and A, and then use the latter to define $\mathfrak{X}$.*
- *The $(0,\infty) \times \mathbb{R}^2$ integral of $(1+t^2)|\mathfrak{X}|^2$ is finite.*
- *Supposing that $R > 0$, then the integral of $|\mathfrak{X}|^2$ over the $|z| > R$ part of $(0,\infty) \times \mathbb{R}^2$ is no greater than $\tfrac{\kappa}{R}$.*
- *Supposing that $R > 0$ and $t \in (0, \infty)$, then the integral of $|\mathfrak{X}|^2$ over the $|z| > R$ part of the slice $\{t\} \times \mathbb{R}^2$ is at most $\tfrac{\kappa}{R}$*

*Proof of Lemma 3.5*: As explained directly, the bounds in the lemma follow from the detailed behavior of the norm of $\mathfrak{X}$, and in particular, the observations that follow. With regards to notation: What is denoted by $c_\diamond$ is a number that is greater than 1 that is determined by $\delta$ and complex numbers $\{a_1, \ldots, a_p\}$.

- Where $t \leq \tfrac{1}{2}\delta$:
  a) $|\mathfrak{X}| \leq c_\diamond \, t^{2(k+1)} |z|^{2m}$ where $|z| \leq t^{(k+1)/p}$.
  b) $|\mathfrak{X}| \leq c_\diamond \frac{|z|^{2k-p}}{t^{k+1}}$ where $t^{(k+1)/p} \leq |z| \leq 4t$ (and thus $|\mathfrak{X}| \leq c_\diamond t^{k-p-1}$).
  c) $\mathfrak{X} \equiv 0$ where $|z| \geq 4t$.
- Where $t \in [\tfrac{1}{2}\delta, 4\delta]$ except in the case when $m = 0$ and $\mu_1 \neq 0$:
  a) $|\mathfrak{X}| \leq c_\diamond (1 + \tfrac{1}{\delta^2})$ where $t \in [\tfrac{1}{2}\delta, \delta]$ and $|z| \leq 4t$.
  b) $|\mathfrak{X}| \leq c_\diamond (1 + \tfrac{1}{\delta^2})\tfrac{1}{|z|^2}$ where $t \in [\tfrac{1}{2}\delta, \delta]$ and $|z| \geq t$.
- Where $t \geq 4\delta$:
  a) *In the case where $m > 0$ or $\mu_1 = 0$*: $|\mathfrak{X}| \equiv 0$ where $t \geq \delta$.
  b) *In the case where $m = 0$ and $\mu_1 \neq 0$*:
    1) $|\mathfrak{X}| = 0$ except where $t \in [\tfrac{1}{2}\delta(1+|z|^2)^{1/4}, \delta(1+|z|^2)^{1/4}]$.

2) $|\mathfrak{X}| \leq c_0 \frac{1}{\delta} \frac{1}{t^2} |\mu_1| \frac{1}{|z|(1+|z|^2)^{1/4}}$ where $t \in [\frac{1}{2}\delta(1+|z|^2)^{1/4}, \delta(1+|z|^2)^{1/4}]$.

(3.33)

These inequalities are justified momentarily. They are used directly to prove the lemma.

To prove the lemma: Except for the $m = 0, \mu_1 \neq 0$ case, the function $|\mathfrak{X}|$ is zero for $t > \delta$. Meanwhile, it is bounded for $t < \delta$ and has compact support for $t < \frac{1}{2}\delta$. For t between $\frac{1}{2}\delta$ and $\delta$, the integral of $|\mathfrak{X}|^2$ is finite on constant t slices because $|\mathfrak{X}|$ is bounded on these slices by a $c_0$ multiple of $\frac{1}{1+|z|^2}$. This bound and what was said just now about the $t < \frac{1}{2}\delta$ and $t > \delta$ behavior imply that its integral on the $|z| > R$ part of these slices is no greater than $c_0$ when $R < 1$, and that it is bounded by $c_0 \frac{1}{R^2}$ when $R > 1$. Because $\mathfrak{X} \equiv 0$ where $t > \delta$, these bounds imply directly that $(1+t^2)|\mathfrak{X}|^2$ has finite $(0,\infty) \times \mathbb{R}^2$ integral.

Turn next to the case when $m = 0$ and $\mu_1 \neq 0$. As in the previous case, $|\mathfrak{X}|$ is bounded with compact support where $t \leq \frac{1}{2}\delta$. Where t is greater than this, the integral of $|\mathfrak{X}|^2$ on the constant t slices is bounded by $c_0 \frac{1}{t^4}$ because $|\mathfrak{X}|^2$ on these slices is bounded by $\frac{1}{t^4}$ times a constant multiple of $\frac{1}{(1+|z|^3)}$ which is integrable (see Item b) of the third bullet in (3.33)). Thus $(1+t^2)$ times the integral of $|\mathfrak{X}|^2$ on the constant t slice are no greater than a constant times $\frac{1}{t^2}$ for $t \geq \frac{1}{2}\delta$ which has finite integral on $[\frac{1}{2}\delta, \infty)$. Moreover, the $c_0 \frac{1}{t^4} \frac{1}{(1+|z|^3)}$ bound for $|\mathfrak{X}|^2$ on these slices implies that its integral on the $|z| > R$ part any such slice is at most $c_0$ if $R < 1$ and at most $c_0 \frac{1}{R}$ if $R > 1$.

A parenthetical remark: The function $|\mathfrak{X}|^2$ would not have finite integral on $(0,\infty) \times \mathbb{R}^2$ for the case $m = 0$ and $\mu_1 \neq 0$ if $\varpi^{\ddagger}_\delta$ were equal to $\varpi_\delta$ because the integral of $|\mathfrak{X}|^2$ on disks of very large radius r about the origin in the constant t slices for $t \in [\frac{1}{2}\delta, \delta]$ would be $\mathcal{O}(|\mu_1|^2 \ln r)$.)

What follows from here is meant justifies the inequalities in (3.33). With regards to the first bullet: Consider first the domain where $t \leq \frac{1}{2}\delta$ and $|z| \leq t$. The connection A and $\mathfrak{a}_3$ in this region are the model solution $A^{(k)}$ and corresponding $\mathfrak{a}_3^{(k)}$ which are depicted in (1.5). This implies that $\nabla_{At}\mathfrak{a}_3 - B_{A3}$ on this domain is $\frac{1}{2}[\varphi^{(k)}, \varphi^{(k)*}]$. Thus,

$$\mathfrak{X} = \tfrac{1}{2}([\varphi^{(k)}, \varphi^{(k)*}] - [\varphi, \varphi^*]),$$

(3.34)

which implies that

$$|\mathfrak{X}| \leq c_0|\varphi|^2|\sigma_3 - \sigma| + c_0||\varphi|^2 - |\varphi^{(k)}|^2|$$

(3.35)

with $\sigma$ given by (3.6) and $\varphi$ given by (3.8). Therefore, from (3.8):

$$|\varphi|^2 = |\varphi^{(k)}|^2(1 + |\gamma|^2)^2,$$

(3.36)

in the region of interest; and thus, because $|\varphi^{(k)}| \le c_0 \frac{|z|^k}{t^{k+1}}$ where $|z| \le t$ and $|\gamma| \le c_0 \frac{t^{k+1}}{|z|^p}$,

- $|\mathfrak{X}| \le c_0 t^{2(k+1)} |z|^{2m}$      where $|z| \le t^{(k+1)/p}$,
- $|\mathfrak{X}| \le c_0 \frac{|z|^{2k-p}}{t^{k+1}} \le c_0 t^{k-p-1}$     where $t^{(k+1)/p} \le |z| \le t$,

(3.37)

(Keep in mind that $t \le \frac{1}{2}\delta$ also in both bullets.)

    Consider next the domain where $t \le \frac{1}{2}\delta$ and $|z| \ge t$. It follows from what is said by the second bullet of Lemma 3.2 and from (3.15) that $|\mathfrak{X}|$ on this domain obeys:

- $|\mathfrak{X}| \le c_0 t^{k-p-1}$  where $t \le |z| \le 4t$.
- $|\mathfrak{X}| = 0$      where $|z| \ge 4t$.

(3.38)

    With regards to the second and third bullets: Because the $|z| \le 4t$ part of the domain where $\frac{1}{2}\delta \le t \le \delta$ is compact, and because the data $(A, \varphi, \mathfrak{a}_3)$ is smooth, there is an upper bound for $|\mathfrak{X}|$, thus $|\mathfrak{X}| \le c_0(1 + \frac{1}{\delta^2})$. Meanwhile, the formula for $w$ in (3.31) and for $\beta$ and $\mathfrak{b}$ in (3.29) lead directly to the asserted bounds where $|z| > t$ and $\frac{1}{2}\delta \le t \le \delta$. Note in this regard that unless $m = 0$ and $\mu_1 \ne 0$, the pair $(A, \mathfrak{a}_3)$ where $t > \delta$ is $\mathrm{Aut}(P)$ equivalent to $(A^{(m)}, \mathfrak{a}^{(m)})$ which has $\mathfrak{X} \equiv 0$. And, were it not for $\beta$, this would also be true for $t \in [\frac{1}{2}\delta, \delta]$. The $\mathrm{Aut}(P)$ equivalence with $(A^{(m)}, \mathfrak{a}^{(m)})$ also holds in the $m = 0$ and $\mu_1 \ne 0$ case where $t > \delta$ and $|z| > t$ unless $t$ is between $\frac{1}{2}\delta(1 + |z|^2)^{1/4}$ and $\delta(1 + |z|^2)^{1/4}$; and this is solely because $\beta$ and $\mathfrak{b}$ are non-zero there. In any event, the formula in (3.29) for $\beta$ and $\mathfrak{b}$ can be used to obtain a bound for $|\mathfrak{X}|$ in this region. In this regard, the worst case terms comes from the $\nabla_{\hat{A}t}$ derivative of $\beta$ whose norm is $\mathcal{O}(\frac{1}{t^2} \frac{1}{|z|^{3/2}})$ and $\hat{A}$-covariant derivatives of $\mathfrak{b}$ along the $\mathbb{R}^2$ factor whose norms are also $\mathcal{O}(\frac{1}{t^2} \frac{1}{|z|^{3/2}})$.

### 4. A 1-parameter family of deformations of $(A, \mathfrak{a})$

    The preceding section described $\mathbb{R}$-invariant pairs $(A, \mathfrak{a})$ that are defined on the whole of $(0, \infty) \times \mathbb{R}^2 \times \mathbb{R}$ that obey the first three bullets of (1.3) and all of (1.9) and (1.10). The plan for this section is to deform an initial choice for $(A, \mathfrak{a})$ so that if it comes from Section 3, then the result of the deformation will be a new pair that obeys all of the bullets in (1.3) and all bullets of (1.9) and (1.10).

#### a) The deformation equations

    By way of background: Suppose for the moment that $(A_\diamond, \varphi_\diamond, \mathfrak{a}_{\diamond 3})$ is a data set consisting of a connection on P, a section of $\mathrm{ad}(P) \otimes_\mathbb{R} \mathbb{C}$ and $\mathrm{ad}(P)$. Any given smooth

section of $P\times_{Ad(SU(2))} Sl(2;\mathbb{C})$ can be used to define a new data set as follows: Letting $\mathfrak{g}$ denote the section, then the new data set (denoted by $(A,\varphi,\mathfrak{a}_3)$) is defined by the rule:

- $(A - A_\diamond)_t - i(\mathfrak{a}_3 - \mathfrak{a}_{\diamond 3}) = -(\nabla_{A_{\diamond t}} \mathfrak{g})\mathfrak{g}^{-1} + i[\mathfrak{a}_{\diamond 3}, \mathfrak{g}]\mathfrak{g}^{-1}$.
- $(A - A_\diamond)_1 + i(A - A_\diamond)_2 = -(\nabla_{A_{\diamond 1}} + i\nabla_{A_{\diamond 2}})\mathfrak{g}\,\mathfrak{g}^{-1}$
- $\varphi = \mathfrak{g}\varphi_\diamond \mathfrak{g}^{-1}$.

(4.1)

The essential point of the preceding construction (which was noted at the very end of Section 2) is that $(A,\varphi,\mathfrak{a}_3)$ obeys the first three bullets of (1.3) if (and only if) the data set $(A_\diamond, \varphi_\diamond, \mathfrak{a}_{\diamond 3})$ obeys those same three bullets.

A 1-parameter family of data sets is described below which is constructed from a corresponding 1-parameter family of sections of $P\times_{Ad(SU(2))} Sl(2;\mathbb{C})$ using the rule in (4.1) starting with a data set that obeys the first three bullets of (1.3). (In the relevant instances, the starting data set is from Section 3.) The point of using (4.1) to define the deformation of the starting data set is to guarantee that each member of the family also obeys the first three bullets of (1.3). The family of sections of $P\times_{Ad(SU(2))} Sl(2;\mathbb{C})$ is parametrized by the half-line $[0,\infty)$ and denoted by $\{\mathfrak{g}|_s\}_{s\in[0,\infty)}$. The plan is to choose this family so that the $s\to\infty$ limit data set obeys the fourth bullet of (1.3).

The desired family $\{\mathfrak{g}|_s\}_{s\in[0,\infty)}$ is defined using a 1-parameter family of smooth sections ad(P) denoted by $\{u|_s\}_{s\in[0,\infty)}$ by integrating point-wise the differential equation

$$\tfrac{\partial}{\partial s} \mathfrak{g} = iu\mathfrak{g}$$

(4.2)

with $\mathfrak{g}|_{s=0}$ being the identity section. Meanwhile, the section $u$ at any given value of s is chosen to obey two constraints. The first constraint is the requirement that

$$(\nabla_{At}^2 + \nabla_{A1}^2 + \nabla_{A2}^2)u - [\mathfrak{a}_3,[u,\mathfrak{a}_3]] - [\mathfrak{a}_1,[u,\mathfrak{a}_1]] - [\mathfrak{a}_2,[u,\mathfrak{a}_2]] = -\mathfrak{X}$$

(4.3)

with $\mathfrak{X} \equiv \nabla_{At}\mathfrak{a}_3 - B_{A3} - \tfrac{i}{2}[\varphi,\varphi^*]$; and with it understood that $(A,\mathfrak{a})$ at the given value of s is defined from $\mathfrak{g}$ at that value of s by (4.1). (The ad(P) components $\mathfrak{a}_1$ and $\mathfrak{a}_2$ of $\mathfrak{a}$ are obtained from $\varphi$ by setting $\mathfrak{a}_1 \equiv \tfrac{1}{2}(\varphi+\varphi^*)$ and $\mathfrak{a}_2 = -\tfrac{i}{2}(\varphi-\varphi^*)$.) The second constraint is that $u$ at any given value of s must have the following asymptotic behavior on $(0,\infty)\times\mathbb{R}^2$:

- $u \to 0$ *as* $t \to 0$ *and as* $t \to \infty$ *uniformily with respect to* $|z|$.
- $u \to 0$ *as* $\tfrac{t}{|z|} \to 0$.

(4.4)

With regards to (4.3): As explained momentarily, if $u$ obeys (4.3), then the $(A, \varphi, \mathfrak{a}_3)|_s$ version of $\mathfrak{X}$ is given in terms of the original version by the rule

$$\mathfrak{X}|_s = e^{-s} \mathfrak{X}|_{s=0} .$$

(4.5)

Thus, the point-wise norm of $\mathfrak{X}|_s$ decreases to zero as $s \to \infty$ at an exponential rate. To see about (4.4), first use (4.1) to see that the s-derivatives of $A$, $\varphi$ and $\mathfrak{a}_3$ are as follows:

- $\frac{\partial}{\partial s} \varphi = i[u, \varphi]$
- $\frac{\partial}{\partial s} \mathfrak{a}_3 = \nabla_{At} u$,
- $\frac{\partial}{\partial s} A_t = -[\mathfrak{a}_3, u]$,
- $\frac{\partial}{\partial s} A_1 = \nabla_{A2} u$,
- $\frac{\partial}{\partial s} A_2 = -\nabla_{A1} u$.

(4.6)

Now use (4.6) and (4.1) to see that

$$\frac{\partial}{\partial s} \mathfrak{X} = (\nabla_{At}^2 + \nabla_{A1}^2 + \nabla_{A2}^2) u - [\mathfrak{a}_3, [u, \mathfrak{a}_3]] - [\mathfrak{a}_1, [u, \mathfrak{a}_1]] - [\mathfrak{a}_2, [u, \mathfrak{a}_2]] ,$$

(4.7)

which implies that $\frac{\partial}{\partial s} \mathfrak{X} = -\mathfrak{X}$ if $u$ obeys (4.3); and integrating the latter gives (4.5).

With regards to (4.4): This rule is enforced so that the each member of the family $\{(A, \varphi, \mathfrak{a}_3)|_s\}_{s \in [0, \infty)}$ and the limit member obeys the conditions in (1.10) if the starting member obeys (1.10) and certain extra constraints that are described in what follows.

**b) A construction of $\mathfrak{g}$**

The existence of the desired path $s \to \mathfrak{g}|_s$ obeying (4.2) with $u$ given by (4.3) and (4.4) is demonstrated via an open/closed argument. To set this up, let $(A_0, \varphi_0, \mathfrak{a}_{03})$ denote the starting data set which is assumed to obey the first three bullets in (1.3). With the starting data set given, let $\mathcal{I} \subset [0, \infty)$ denote the set with the following property: If $s_\diamond \in \mathcal{I}$, then there is a differentiable path $s \to \mathfrak{g}|_s$ parametrized by $s \in [0, s_\diamond]$ with $\mathfrak{g}|_{s=0}$ being the identity section and such that (4.2) holds with $u$ obeying (4.3) and (4.4) at each $s \in [0, s_\diamond]$. The set $\mathcal{I}$ is non-empty because $0$ is in $\mathcal{I}$. Therefore, if $\mathcal{I}$ is both open and closed, then $\mathcal{I}$ must be the whole half-line $[0, \infty)$. As explained momentarily, the assertion that $\mathcal{I}$ is an open set (given an extra condition on the starting data set) is an instance of the upcoming Proposition 4.1 Section 5 contains the proof that $\mathcal{I}$ is closed (see Proposition 5.1) and, given that $\mathcal{I} = [0, \infty)$, that there is an $s \to \infty$ limit section of $P \times_{\text{Ad}(SU(2))} Sl(2; \mathbb{C})$. Propositions 5.1, 5.2 and 6.1 say more about the $(0, \infty) \times \mathbb{R}^2$ asymptotics of each member of the family $\{\mathfrak{g}|_s\}_{s \in [0, \infty)}$ and the asymptotics of the $s \to \infty$ limit section.

**Proposition 4.1**: *Let $(A_\diamond, \varphi_\diamond, \mathfrak{a}_{\diamond 3})$ denote a data set that obeys the first three bullets in (1.3) and whose version of $\mathfrak{X}$ (denoted by $\mathfrak{X}_\diamond$, thus $\nabla_{A_{\diamond t}}\mathfrak{a}_{\diamond 3} - B_{A_\diamond} - \frac{i}{2}[\varphi_\diamond, \varphi_\diamond^*]$) is such that $(1+t^2)|\mathfrak{X}_\diamond|^2$ has finite integral on $(0,\infty) \times \mathbb{R}^3$. Granted these assumptions, then there exists $\varepsilon > 0$ and a smooth map $r \to u|_r$ from $(-\varepsilon, \varepsilon)$ to the space of smooth sections of $\mathrm{ad}(P)$ with $u|_0 \equiv 0$ and such that the following is true: Define a map $r \to \mathfrak{h}|_r$ from $(-\varepsilon, \varepsilon)$ to the space of sections of $P \times_{\mathrm{Ad}(SU(2))} \mathrm{Sl}(2; \mathbb{C})$ by solving the ordinary differential equation $\frac{\partial}{\partial r}\mathfrak{h} = i u \mathfrak{h}$ with $\mathfrak{h}|_{r=0}$ being the identity section. If $r \in (-\varepsilon, \varepsilon)$ and if $(A|_r, \varphi|_r, \mathfrak{a}_3|_r)$ is defined from $\mathfrak{h}|_r$ via*

- $(A - A_\diamond)_t - i(\mathfrak{a}_3 - \mathfrak{a}_{\diamond 3}) = -(\nabla_{A_{\diamond t}} \mathfrak{h})\mathfrak{h}^{-1} + i[\mathfrak{a}_{\diamond 3}, \mathfrak{h}]\mathfrak{h}^{-1}$.
- $(A - A_\diamond)_1 + i(A - A_\diamond)_2 = -(\nabla_{A_{\diamond 1}} + i\nabla_{A_{\diamond 2}})\mathfrak{h}\,\mathfrak{h}^{-1}$.
- $\varphi - \varphi_\diamond = -[\varphi_\diamond, \mathfrak{h}]\mathfrak{h}^{-1}$.

*then $u|_r$ obeys the $(A, \varphi, \mathfrak{a}_3) = (A|_r, \varphi|_r, \mathfrak{a}_3|_r)$ version of (4.3) and (4.4) and each $r \in (-\varepsilon, \varepsilon)$ member of the family $\{(A|_r, \varphi|_r, \mathfrak{a}_3|_r)\}_{r \in (-\varepsilon, \varepsilon)}$ obeys the $s = r$ version of (4.6) and has corresponding $\mathfrak{X}|_r$ equal to $e^{-r}\mathfrak{X}_\diamond$.*

Section 4c presents the analysis that is needed for the proof of Proposition 4.1. The proof itself is in Section 4d modulo two lemmas that are proved in Sections 4e and 4f.

To apply Proposition 4.1 to the question of whether $\mathcal{I}$ is open: Starting with a data set $(A_0, \varphi_0, \mathfrak{a}_{30})$ that obeys (1.3) with the $(0, \infty) \times \mathbb{R}^2$ integral of $(1+t^2)|\mathfrak{X}_\diamond|^2$ being finite, there exists a non-negative number $s_\diamond$ such that the family $s \to \mathfrak{g}|_s$ is defined on the interval $[0, s_\diamond]$ so as to have the desired properties, which is that $\{\mathfrak{g}|_s\}_{s \in [0, s_\diamond]}$ obeys (4.2) with $u|_s$ obeying (4.3)-(4.4) when $(A|_s, \varphi|_s, \mathfrak{a}_3|_s)$ is defined from $\mathfrak{g}|_s$ via (4.1). Note that $s_\diamond$ can equal 0. Let $\mathfrak{g}_\diamond$ denote the $s = s_\diamond$ end member of the family $\{\mathfrak{g}|_s\}_{s \in [0, s_\diamond]}$ and let $(A_\diamond, \varphi_\diamond, \mathfrak{a}_{\diamond 3})$ denote the end member of the corresponding family $\{(A|_s, \varphi|_s, \mathfrak{a}_3|_s)\}_{s \in [0, s_\diamond]}$. Use this $(A_\diamond, \varphi_\diamond, \mathfrak{a}_{\diamond 3})$ as input to Proposition 4.1 to define a positive $\varepsilon$ and the map $s \to \mathfrak{h}|_s$ as described in the proposition. Then, define $\mathfrak{g}|_s$ for $s \in [s_\diamond, s_\diamond + \varepsilon)$ by setting $\mathfrak{g}|_s = \mathfrak{h}|_{r=s-s_\diamond}\mathfrak{g}_\diamond$. The extended family $\{\mathfrak{g}|_s\}_{s \in [0, s_\diamond + \varepsilon]}$ is continuous because the respective $s \to s_\diamond$ limits from above and from below are $\mathfrak{g}_\diamond$.

To prove that the extended family $\{\mathfrak{g}|_s\}_{s \in [0, s_\diamond + \varepsilon]}$ is differentiable across $s_\diamond$, note first that the s-derivative at $s_\diamond$ of the respective $[0, s_\diamond]$ part of the family and the $[s_\diamond, s_\diamond + \varepsilon]$ part exist, the former by assumption and the guaranteed by Proposition 4.1. By virtue of (4.2), the former can be written as $iu_<\mathfrak{g}|_{s_\diamond}$ and the latter as $iu_>\mathfrak{g}|_{s_\diamond}$ with $u_<$ and $u_>$ being sections of $\mathrm{ad}(P)$ that both obey (4.3) and (4.4). Because both $u_<$ and $u_>$ obey (4.3), their difference obeys the modification of (4.3) that replaces $\mathfrak{X}$ with 0 on the right hand side. As a consequence of that equation, the function $|u_< - u_>|^2$ obeys $\Delta|u_< - u_>|^2 \geq 0$ on $(0,\infty) \times \mathbb{R}^2$ which implies via the maximum principle that $|u_< - u_>|^2 = 0$ since $u_< - u_>$ also obeys (4.4). Therefore, the respective derivatives of the $[0, s_\diamond]$ and $[s_\diamond, s_\diamond + \varepsilon)$ parts of the extended family agree at $s_\diamond$ and so the extended family is differentiable across $s = s_\diamond$.

### c) Analytic background for the proof of Proposition 4.1

The four parts of this subsection supply the analytic tools that are used subsequently to solve (4.3)–(4.4). In this regard, the focus in here is on the equation

$$-(\nabla_{At}^2 + \nabla_{A1}^2 + \nabla_{A2}^2)v + [\mathfrak{a}_3, [v, \mathfrak{a}_3]] + [\mathfrak{a}_1, [v, \mathfrak{a}_1]] + [\mathfrak{a}_2, [v, \mathfrak{a}_2]] = \mathfrak{P}$$
(4.8)

when $\mathfrak{P}$ is a given section of ad(P) and when $(A, \varphi, \mathfrak{a}_3)$ is any given data set comprising a connection on P, a section of $ad(P) \otimes_\mathbb{R} \mathbb{C}$ and section of ad(P). As always, $\mathfrak{a}_1 = \frac{1}{2}(\varphi + \varphi^*)$ and $\mathfrak{a}_2 = \frac{i}{2}(\varphi - \varphi^*)$. The analysis for (4.8) is summarized by the upcoming Lemma 4.2.

*Part 1*: A preliminary digression is required to introduce the Green's function for the standard Laplacian on $(0,\infty) \times \mathbb{R}^2$ with Dirichlet boundary condition at $t = 0$. (The standard Laplacian acts on functions as the sum of the second derivatives in the coordinate direction. It is denoted by $\Delta$ in what follows.) The version of the Green's function with pole at a designated point $q \in (0,\infty) \times \mathbb{R}^2$ is denoted by $G_q$. Writing $q$ as $(t_q, z_q)$, then $G_q$ at any point $(t, z)$ other than $q$ is given by the formula

$$G_q(t,z) = \frac{1}{4\pi} \frac{1}{((t-t_q)^2 + |z-z_q|^2)^{1/2}} - \frac{1}{4\pi} \frac{1}{((t+t_q)^2 + |z-z_q|^2)^{1/2}}.$$
(4.9)

An important point for later: This Green's function $G_q$ is zero where $t = 0$ and it is strictly positive on the complement of $q$ in $(0, \infty) \times \mathbb{R}^2$. What follows are other important points to keep in mind about $G_q$:

- $G_q > 0$ *on* $(0, \infty) \times \mathbb{R}^2$ *and* $G_q(p) \geq \frac{1}{c_0} \frac{1}{|p-q|}$ *if* $|p-q| \leq \frac{1}{4} t_q$.
- $G_q(p) \leq c_0 \frac{1}{|p-q|}$ *at any given* $p \in (0,\infty) \times \mathbb{R}^2$.
- $G_q(p) \leq c_0 \frac{t_p t_q}{|p-q|^3}$ *at any given* $p \in (0,\infty) \times \mathbb{R}^2$ *with* $|p-q| \geq \frac{1}{4} t_q$ *or* $|p-q| \geq \frac{1}{4} t_p$.
- $|\nabla G_q(p)| \leq c_0 \frac{1}{|p-q|^2}$ *at any given* $p \in (0,\infty) \times \mathbb{R}^2$.

(4.10)

And, here is a key consequence of the second and third bullets of (4.10):

$$\int_{(0,\infty) \times \mathbb{R}^2} \frac{1}{1+t^2} G_q^2 \leq c_0 \min(t_q, \frac{1}{t_q}(1 + |\ln t_q|)).$$
(4.11)

To prove this, consider the integral over the domain where $|p-q| < \frac{1}{2} t_q$ and then the integral over the complimentary domain. Use the bound in the second bullet of the lemma for the former and the bound in the third bullet for the latter. (It also helps to integrate first over the constant t slices and then integrate over t.

*Part 2* A Banach space for the analyis of (4.8) is defined in this part of subsection using the Green's function. The first step towards the definition is to introduce an auxilliary Hilbert space, which is the completion of the space of smooth, compactly supported sections of ad(P) using the norm whose square is the function

$$v \to \int_{(0,\infty)\times\mathbb{R}^2} (|\nabla_A v|^2 + |[\mathfrak{a},v]|^2) \,. \tag{4.12}$$

This Hilbert space is denoted by $\mathbb{H}$ and its norm is denoted by $\|\cdot\|_{\mathbb{H}}$. The desired Banach space is a vector subspace in $\mathbb{H}$ which is denoted by $\mathbb{B}$; and a given element $v \in \mathbb{H}$ is in $\mathbb{B}$ when the function

$$q \to \int_{(0,\infty)\times\mathbb{R}^2} (1+G_q)(|\nabla_A v|^2 + |[\mathfrak{a},v]|^2) + \left( \int_{(0,\infty)\times\mathbb{R}^2} G_q |\Delta_{(A,\mathfrak{a})} v| \right)^2 \tag{4.13}$$

is bounded on $(0,\infty)\times\mathbb{R}^2$. Here and subsequently, $\Delta_{(A,\mathfrak{a})} u$ denotes the expression on the left hand side of (4.8). The square of the norm on $\mathbb{B}$ is the supremum of the function on $(0,\infty)\times\mathbb{R}^2$ that is depicted on the right hand side of (4.13). This Banach space norm is denoted in what follows by $\|\cdot\|_{\mathbb{B}}$.

An important point about the Hilbert space $\mathbb{H}$: Because $|\nabla_A(\cdot)| \geq |\nabla|(\cdot)||$, the norm of any element in $\mathbb{H}$ is in the Hilbert completion of the space of compactly supported functions on $(0,\infty)\times\mathbb{R}^2$ using the norm whose square is the function

$$f \to \int_{(0,\infty)\times\mathbb{R}^2} |\nabla f|^2 \,. \tag{4.14}$$

For the present purposes, the important point with regards to the latter Hilbert space is that its elements obey Hardy's inequality [HLP], which is this: If $f$ is in the Hilbert space, then the function $\frac{1}{t^2} f^2$ has finite integral on $(0,\infty)\times\mathbb{R}^2$ and that integal obeys the bound

$$\int_{(0,\infty)\times\mathbb{R}^2} \tfrac{1}{t^2} f^2 \leq 4 \int_{(0,\infty)\times\mathbb{R}^2} |\nabla f|^2 \,. \tag{4.15}$$

The proof amounts to an integration by parts (see the appendix for the derivation). To summarize: If u is in $\mathbb{B}$, or just in $\mathbb{H}$, then $f = |u|$ obeys (4.15).

An important point about the Banach space $\mathbb{B}$: If $v \in \mathbb{B}$, then the norm of $v$ is a priori bounded on the whole of $(0,\infty) \times \mathbb{R}^2$ because

$$|v|(q) \leq \int_{(0,\infty) \times \mathbb{R}^2} G_q |\Delta_{(A,\mathfrak{a})} v|$$

(4.16)

whose right hand side is (tautologically) no greater than $\|v\|_{\mathbb{B}}$.

To prove the preceding inequalitity, invoke the identity

$$-\tfrac{1}{2} \Delta |v|^2 + |\nabla_A v|^2 + |[\mathfrak{a}, v]|^2 = -\langle v \Delta_{(A,\mathfrak{a})} v \rangle$$

(4.17)

which implies in turn this inequality:

$$-\Delta |v| \leq |\Delta_{(A,\mathfrak{a})} v| .$$

(4.18)

Fix a point $q \in (0, \infty) \times \mathbb{R}^2$ and then a very small but positive number (to be denoted by $\delta$) which is much less than the lesser of $t_q$ and $\tfrac{1}{t_q}$. Having done that, let $\pi_\delta$ denote the bump function $\chi(2(1-\tfrac{t}{\delta}))\chi(\delta(t^2+|z|^2)^{1/2}-1)$. To be sure: This function is equal to 1 where t is greater than $\delta$ and $(t^2+|z|^2)^{1/2}$ is less than $\tfrac{1}{\delta}$ and equal to zero where t is less than $\tfrac{1}{2}\delta$ or $(t^2+|z|^2)^{1/2}$ is greater than $2\tfrac{1}{\delta}$. Thus, it has compact support in $(0,\infty) \times \mathbb{R}^2$.) Multiply both sides of (4.18) by $\pi_\delta G_q$ and integrate the result over $(0,\infty) \times \mathbb{R}^2$; then integrate by parts twice on the left hand side of the resulting inequality to bound the norm of $v$ at q by the integral of $G_q |\Delta_{(A,\mathfrak{a})} v|$ plus terms that have derivatives on $\pi_\delta$. Then, take $\delta$ ever smaller with limit zero and use Hardy's inequality and the inequalities for $G_q$ in (4.10) to see that the terms with derivatives on $\pi_\delta$ have limit zero as $\delta \to 0$.

An extension of the preceding definitions: Suppose now that V is a vector space with a representation of SU(2). There is a version of the Banach space $\mathbb{B}$ and its $\mathbb{B}$-norm for sections of the associated vector bundle $P \times_{SU(2)} V$. The definitions and the surrounding discussions are the same as those given above for the case when the associated vector bundle is ad(P) with it understood that $u$ is a section of the bundle $P \times_{SU(2)} V$. The generic associated vector bundle version of the $\mathbb{B}$-norm will not be distinguished notationally from the ad(P) version. (The relevant case has V being the vector space of 2×2 complex matrices with its adjoint representation of SU(2).)

*Part 3*: This part of the subsection states and proves a lemma that provides the basic analytic tool for the proof of Proposition 4.1. To set things up, suppose that $\mathfrak{P}_0$ and $\mathfrak{P}_1$ are sections of ad(P) and that $\mathfrak{b} = (\mathfrak{b}_t, \mathfrak{b}_1, \mathfrak{b}_2)$ and $\mathfrak{c} = (\mathfrak{c}_1, \mathfrak{c}_2, \mathfrak{c}_3)$ are 3-tuples of sections of ad(P); and that these have the following respective properties:

- *The function $(1+t^2)|\mathfrak{P}_0|^2$ has finite integral on $(0,\infty) \times \mathbb{R}^2$.*
- *The $(0,\infty) \times \mathbb{R}^2$ integral of any function from the set $\{(1+G_q)|\mathfrak{P}_1|\}_{q \in (0,\infty) \times \mathbb{R}^2}$ is finite, and this set of integrals is bounded.*
- *The $(0,\infty) \times \mathbb{R}^2$ integrals of $|\mathfrak{b}|^2$ and $|\mathfrak{c}|^2$ are finite, as are the integrals of the functions from the set $\{G_q|\nabla_{Ai}\mathfrak{b}_i + [\mathfrak{a}_{3a}, \mathfrak{c}_a]|\}_{q \in (0,\infty) \times \mathbb{R}^2}$; and the latter set of $(0,\infty) \times \mathbb{R}^2$ integrals has an upper bound.*

(4.19)

(There is an implicit sum over the indexing set $\{t, 1, 2\}$ for the index i in the third bullet, and likewise an implicit sum over the indexing set $\{1, 2, 3\}$ for the index a.)

The upcoming lemma uses $\mathcal{Z}$ to denote an upper bound for the set of integrals in (4.19)'s second bullet and the set of integrals with $G_q$ from the third bullet; and also for the square root of the integral of $(1+t^2)|\mathfrak{P}_0|^2$ from (4.19)'s first bullet and also for the square root of the integrals of $|\mathfrak{b}|^2$ and $|\mathfrak{c}|^2$ from (4.19)'s third bullet.

**Lemma 4.2**: *There exists $\kappa > 1$ with the following significance: Fix a data set $(A, \varphi, \mathfrak{a}_3)$ with A being a connection on P, and with $\varphi$ and $\mathfrak{a}_3$ denoting respective sections of $\mathrm{ad}(P) \otimes_{\mathbb{R}} \mathbb{C}$ and $\mathrm{ad}(P)$. Supposing that $\mathfrak{P}_0$ and $\mathfrak{P}_1$ and $\mathfrak{b}$ and $\mathfrak{c}$ are described by (4.19), then there exists a unique element $v \in \mathbb{B}$ obeying the equation*

$$-(\nabla_{At}^2 + \nabla_{A1}^2 + \nabla_{A2}^2)v + [\mathfrak{a}_3, [v, \mathfrak{a}_3]] + [\mathfrak{a}_1, [v, \mathfrak{a}_1]] + [\mathfrak{a}_2, [v, \mathfrak{a}_2]] = \mathfrak{P}_0 + \mathfrak{P}_1 + \nabla_{Ai}\mathfrak{b}_i + [\mathfrak{a}_a, \mathfrak{c}_a] .$$

*Moreover, the $\mathbb{B}$ norm of $u$ obeys $\|v\|_{\mathbb{B}} \leq \kappa \mathcal{Z}$.*

*Proof of Lemma 4.2*: Since the equation is linear, it is sufficient to consider only the separate cases where one of $\mathfrak{P}_0$ and $\mathfrak{P}_1$ and $(\mathfrak{b}, \mathfrak{c})$ is not identically zero. The proof has six steps. Since the case where $\mathfrak{P}_1 \equiv 0$ and $(\mathfrak{b}, \mathfrak{c}) \equiv 0$ is simplest, it is treated first, this being the content of the Steps 1-3. This case is needed in any event for the other cases. The remaining steps treat those cases.

Step 1: The desired $v$ will be obtained by minimizing a certain function on the Hilbert space $\mathbb{H}$. This step in the proof defines the function and derives some of its salient properties.

The function in question is denoted by $\mathcal{E}$; it is defined first on the space of compactly supported sections of ad(P) by the rule

$$x \to \mathcal{E}(x) \equiv \tfrac{1}{2} \int_{(0,\infty)\times\mathbb{R}^2} (|\nabla_A x|^2 + |[\mathfrak{a},x]|^2) - \int_{(0,\infty)\times\mathbb{R}^2} \langle x\mathfrak{P}_0 \rangle \ .$$

(4.20)

The left most integral on the right hand side is half of the norm on $\mathbb{H}$; and the right most integral defines a bounded, linear function on $\mathbb{H}$ whose norm is bounded by

$$( \int_{(0,\infty)\times\mathbb{R}^2} (1+t^2)|\mathfrak{P}_0|^2 \ )^{1/2} \|x\|_{\mathbb{H}} \ .$$

(4.21)

This bound follows from Hardy's inequality in (4.15).

<u>Step 2</u>: This step of the proof explains why $\mathcal{E}$ has a unique minimizer in $\mathbb{H}$. This minimizer is the desired solution.

By virtue of the bound in (4.21) for the norm of the right most integral in (4.20), the function $\mathcal{E}$ on $\mathbb{H}$ is bounded from below; and so it has an infimum on $\mathbb{H}$ which is a priori negative. Because the left most integral on the right hand side of (4.20) is a positive multiple of the norm on $\mathbb{H}$, any minimizing sequence for $\mathcal{E}$ in $\mathbb{H}$ has bounded norm; which implies that it converges weakly in $\mathbb{H}$ and also that the value of $\mathcal{E}$ on this weak limit is the infimum of $\mathcal{E}$ on $\mathbb{H}$. This fact implies in turn that the weak convergence is actually strong convergence. (The space $\mathbb{H}$ is separable because the space of smooth, compactly supported sections of ad(P) is separable.)

Let $v$ denote the weak limit. Because $v$ is a critical point of $\mathcal{E}$, it a priori obeys the desired equation in a distributional sense, which means this: If w is also in $\mathbb{H}$, then

$$\int_{(0,\infty)\times\mathbb{R}^2} (\langle\nabla_A w, \nabla_A v\rangle + \langle[\mathfrak{a},w],[\mathfrak{a},v]\rangle) - \int_{(0,\infty)\times\mathbb{R}^2} \langle w\mathfrak{P}_0 \rangle = 0 \ .$$

(4.22)

This equation implies via standard elliptic regularity arguments (see e.g. Theorem 5.6.2 in [M]]) that the lemma's equation is obeyed as an identity between smooth sections of ad(P). (Remember that $\mathfrak{P}_0$ is assumed to be a smooth section of ad(P).)

To see that $v$ is unique, suppose for the moment that $v´$ is in $\mathbb{H}$ and that (4.22) holds for any w from $\mathbb{H}$ with $v$ replaced by $v´$. Then $v - v´$ obeys the $\mathfrak{P}_0 \equiv 0$ version of (4.22) for any choice of w from $\mathbb{H}$, and in particular, for $w = v - v´$. The latter version of (4.22) says that $v - v´ = 0$.

Step 3: This step in the proof derives the asserted bound on $\|v\|_{\mathbb{B}}$. To do this, note first that the bound

$$\int_{(0,\infty)\times\mathbb{R}^2} (|\nabla_A v|^2 + |[\mathfrak{a},v]|^2) \leq c_0 \int_{(0,\infty)\times\mathbb{R}^2} (1+t^2)|\mathfrak{P}_0|^2$$

(4.23)

follows from the $w = v$ version of (4.22) and the $x = v$ version of (4.21)'s bound for the norm of the right most integral on the right hand side of (4.20).

To obtain the other parts of the norm: Fix $q \in (0,\infty)\times\mathbb{R}^2$. Because the integral of $G_q|\mathfrak{P}_0|$ is bounded by $c_0$ times the product of the square roots of the integrals of $\frac{1}{1+t^2}G_q^2$ and $(1+t^2)|\mathfrak{P}_0|^2$, it is bounded by (see (4.11)) the square root of

$$c_0 \min(t_q, \tfrac{1}{t_q}(1+|\ln t_q|)) \int_{(0,\infty)\times\mathbb{R}^2} (1+t^2)|\mathfrak{P}_0|^2 \;.$$

(4.24)

This is a q-independent upper bound for the $(0,\infty)\times\mathbb{R}^2$ integral of $G_q|\Delta_{(A,\mathfrak{a})}v|$ which is that the square of this integral is at most $c_0$ times what appears on the right hand side of (4.23).

To continue: The inequality in (4.17) in this case says that

$$-\tfrac{1}{2}\Delta|v|^2 + |\nabla_A v|^2 + |[\mathfrak{a},v]|^2 = \langle v\mathfrak{P}_0\rangle \;.$$

(4.25)

Multiply both sides of this by $G_q$ and then integrate result over $(0,\infty)\times\mathbb{R}^2$. As was the case when this was done with (4.18), integration by parts can be applied twice to the left hand side, and doing that leads to the inequality

$$\tfrac{1}{2}|v|^2(q) + \int_{(0,\infty)\times\mathbb{R}^2} G_q(|\nabla_A v|^2 + |[\mathfrak{a},v]|^2) \leq \int_{(0,\infty)\times\mathbb{R}^2} G_q|v||\mathfrak{P}_0| \;.$$

(4.26)

Granted the bound for $|v|$ given by (4.16) and (4.24), this inequality leads to the bound

$$\int_{(0,\infty)\times\mathbb{R}^2} G_q(|\nabla_A v|^2 + |[\mathfrak{a},v]|^2) \leq c_0 \min(\sqrt{t_q}, \tfrac{1}{\sqrt{t_q}}(1+|\ln t_q|)^{1/2}) \int_{(0,\infty)\times\mathbb{R}^2} (1+t^2)|\mathfrak{P}_0|^2$$

(4.27)

which gives the rest of the asserted $\|\cdot\|_{\mathbb{B}}$-norm bound for $v$.

Step 4: This step and the Step 5 consider the case when $\mathfrak{P}_0 \equiv 0$ and $(\mathfrak{b},\mathfrak{c}) \equiv 0$ and $\mathfrak{P}_1 \neq 0$. As with the previous case, it is sufficient to assume that $\mathfrak{P}_1$ is smooth.

To start the story on $\mathfrak{P}_1$, fix for the moment a number $\delta \in (0, 1)$ and use it to define the bump function $\pi_\delta$ (which is $\chi(2(1-\frac{t}{\delta})\chi(\delta(t^2+|z|^2)^{1/2}-1))$. Since $\pi_\delta$ has compact support, the function $(1+t^2) \pi_\delta^2 |\mathfrak{P}_1|^2$ has finite integral on $(0,\infty)\times\mathbb{R}^2$. As a consequence of this (and what was proved in Steps 1-3), there exists a unique $v_\delta$ in $\mathbb{H}$ obeying the version of the lemma's differential equation with $\mathfrak{P}_1$ replaced by zero and with $\mathfrak{P}_0$ replaced by $\pi_\delta \mathfrak{P}_1$.

To see about the $\|\cdot\|_\mathbb{B}$ norm of $v_\delta$: It follows from the assumption about the norm of $\mathfrak{P}_1$ in the second bullet of (4.19) that any $q \in (0,\infty)\times\mathbb{R}^2$ version of $G_q|\Delta_{(A,\mathfrak{a})}u|$ has finite integral on $(0,\infty)\times\mathbb{R}^2$ with its integral bounded by $c_0$ times what the lemma calls $\mathcal{Z}$. As a consequence of this and (4.16), the norm of $|v_\delta|$ at any given point in $(0,\infty)\times\mathbb{R}^2$ is also bounded by $c_0\mathcal{Z}$.

With the latter bound understood, note that (4.22) holds with $v$ replaced by $v_\delta$ and with $\mathfrak{P}_0$ replaced by $\pi_\delta \mathfrak{P}_1$. Taking $w = v_\delta$ in this inequality and using the $c_0\mathcal{Z}$ bound for the sup-norm of $|v_\delta|$ leads to the $\delta$-independent bound $\|v_\delta\|_\mathbb{H} \leq c_0\mathcal{Z}$. By the same token, (4.26) holds with $v = v_\delta$ and $\mathfrak{P}_0$ replaced by $\pi_\delta \mathfrak{P}_1$. This version of (4.26) with the $c_0\mathcal{Z}$ sup-norm bound on $|v|$ leads to a $c_0\mathcal{Z}^2$ bound for the $(0,\infty)\times\mathbb{R}^2$ integral of $G_q(|\nabla_A v_\delta|^2 + |[\mathfrak{a}, v_\delta]|^2)$. Taking the bounds in the preceding paragraph together leads to an $\delta$-independent bound by $c_0\mathcal{Z}$ bound for $\|v_\delta\|_\mathbb{B}$.

<u>Step 5</u>: The desired section $v$ for the $\mathfrak{P}_0 \equiv 0$ and $\mathfrak{P}_1 \neq 0$ version of the lemma is obtained by taking the $\delta \to 0$ limit of the collection $\{v_\delta\}_{\delta \in (0,1)}$. In this regard, weak convergence of this sequence in $\mathbb{H}$ follows because these all lie in the radius $c_0\mathcal{Z}$ ball about the origin. In fact, the convergence in $\mathbb{H}$ is strong convergence, a fact which proved as follows: Fix $\delta > 0$ and then any $\delta' \in (0, \delta)$. Having done that, subtract the $v = v_\delta$ and $\mathfrak{P}_0 \equiv \pi_\delta \mathfrak{P}_1$ version of (4.22) from the corresponding version with $\delta$ replaced by $\delta'$ in all occurrences. Because both $|v_\delta|$ and $|v_{\delta'}|$ are bounded by $c_0\mathcal{Z}$, the $w = v_{\delta'} - v_\delta$ version of the difference of the $\delta$ and $\delta'$ versions of (4.20) leads to the following:

$$\|v_{\delta'} - v_\delta\|_\mathbb{H}^2 \leq c_0 \mathcal{Z} \int_{(0,\infty)\times\mathbb{R}^2} (1-\pi_\delta)|\mathfrak{P}_1|$$

(4.28)

The fact that the $\delta \to 0$ limit of the right hand side of (4.28) is zero (the dominated convergence theorem) implies the strong convergence claim.

Let $v$ denote the $\delta \to 0$ limit of the $\{v_\delta\}_{\delta > 0}$. This obeys the lemma's differential equation so it is smooth (assuming $\mathfrak{P}_1$ is smooth). It also has a $c_0\mathcal{Z}$ bound to the $(0,\infty)\times\mathbb{R}^2$ integral of $G_q|\Delta_{(A,\mathfrak{a})}v|$ and $|v|(q)$ for any given $q \in (0,\infty)\times\mathbb{R}^2$; and also a $c_0\mathcal{Z}^2$ bound to the $(0,\infty)\times\mathbb{R}^2$ integral of $G_q(|\nabla_A v|^2 + |[\mathfrak{a},v]|^2)$ for any given $q \in (0,\infty)\times\mathbb{R}^2$. To elaborate with regards to the latter bounds, fix some small $\delta$ and then $\delta' \in (0,\delta)$. Then (4.26) holds with $v$ replaced by $v_{\delta'} - v_\delta$ and $\mathfrak{P}_0$ replaced by $(\pi_{\delta'} - \pi_\delta)\mathfrak{P}_1$; and that leads in turn to any given $q \in (0,\infty)\times\mathbb{R}^2$ version of the following

$$\tfrac{1}{2}|v_{\delta'} - v_\delta|^2(q) + \int_{(0,\infty)\times\mathbb{R}^2} G_q(|\nabla_A(v_{\delta'} - v_\varepsilon)|^2 + |[\mathfrak{a}, v_{\delta'} - v_\delta]|^2) \le c_0 Z \int_{(0,\infty)\times\mathbb{R}^2} G_q(1-\pi_\delta)|\mathfrak{P}_1|$$

(4.29)

Since the right hand side of (4.29) has limit zero as $\delta \to 0$ for any fixed $q$, and since these limits are uniform with respect to $q$ if it varies in any given compact set in $(0,\infty)\times\mathbb{R}^2)$, it follows as a direct consequence that

$$\lim_{\delta \to 0} |v_{\delta'} - v_\delta|(q) = 0 \quad and \quad \lim_{\delta \to 0} \int_{(0,\infty)\times\mathbb{R}^2} G_q(|\nabla_A(v_{\delta'} - v_\varepsilon)|^2 + |[\mathfrak{a}, v_{\delta'} - v_\delta]|^2) = 0,$$

(4.30)

with the limits being uniform with respect to $q$ as $q$ varies in compact subsets of the domain $(0,\infty)\times\mathbb{R}^2$. These limits implies that $\|v\|_\mathbb{B}$ is bounded and thus that $v$ is in $\mathbb{B}$.

<u>Step 6</u>: This last step considers the case where $\mathfrak{P}_0 \equiv 0$ and $\mathfrak{P}_1 \equiv 0$ with $\mathfrak{b}$ and/or $\mathfrak{c}$ not identically zero. In this case $u$ is found by minimizing the function

$$x \to \mathcal{E}_*(x) \equiv \tfrac{1}{2} \int_{(0,\infty)\times\mathbb{R}^2} (|\nabla_A x|^2 + |[\mathfrak{a}, x]|^2) + \int_{(0,\infty)\times\mathbb{R}^2} (\langle \nabla_A x, \mathfrak{b}\rangle + \langle [\mathfrak{a}, x], \mathfrak{c}\rangle)$$

(4.31)

Since this function is bounded from below by

$$\tfrac{1}{4}\|x\|_\mathbb{H}^2 - \int_{(0,\infty)\times\mathbb{R}^2} (|\mathfrak{b}|^2 + |\mathfrak{c}|^2),$$

(4.32)

Step 1's argument with only cosmetic changes proves that there is a unique minimizer in $\mathbb{H}$ for $\mathcal{E}_*$ and that said minimizer obeys the desired equation in a weak sense. Because the various $q \in (0,\infty)\times\mathbb{R}^2$ versions of the $(0,\infty)\times\mathbb{R}^2$ integral of $G_q|\nabla_{Ai}\mathfrak{b}_i + [\mathfrak{a}_a, \mathfrak{c}_a]|$ have a q-independent upper bound, cosmetic modifications to the arguments in Steps 4 and 5 prove that the minimizer of $\mathcal{E}_*$ in $\mathbb{H}$ comes from the Banach space $\mathbb{B}$ and that it has the asserted $\mathbb{B}$-norm bound.

**d) Proof of Proposition 4.1**

As in Proposition 4.1, let $(A_\diamond, \varphi_\diamond, \mathfrak{a}_{\diamond 3})$ denote a data set obeying the first three bullets of (1.3) with the $(0,\infty) \times \mathbb{R}^2$ integral of $(1+t^2)|\mathfrak{X}_\diamond|^2$ being finite. (Remember that $\mathfrak{X}_\diamond$ is $\nabla_{A_\diamond t}\mathfrak{a}_{\diamond 3} - B_{A_\diamond} - \frac{i}{2}[\varphi_\diamond, \varphi_\diamond^*]$.) Use $(A_\diamond, \varphi_\diamond, \mathfrak{a}_{\diamond 3})$ to define the Banach space $\mathbb{B}$. By way of a look ahead, the proof of Proposition 4.1 uses a contraction mapping theorem on the Banach space of differentiable maps $r \to x|_r$ from an interval around 0 in $\mathbb{R}$ to the Banach space $\mathbb{B}$. The details of this are given below in Parts 1-4 of what follows.

*Part 1*: Suppose for the moment that $\varepsilon$ is positive and that $u$ is a smooth map from the interval $(-\varepsilon, \varepsilon)$ into the space of sections of $\text{ad}(P)$. The solution to the equation

$$\tfrac{\partial}{\partial r} \mathfrak{h} = i u \mathfrak{h}$$

(4.33)

which is the identity section at $r = 0$ can be written as the convergent series

$$\mathfrak{h} = \mathbb{I} + \int_{0 \le r_1 \le r} u|_{r_1} + \int_{0 \le r_2 \le r_1 \le r} u|_{r_1} u|_{r_2} + \cdots$$

(4.34)

The data $(A, \varphi, \mathfrak{a}_3)$ at any given $r \in (-\varepsilon, \varepsilon)$ can then be written in terms of $\mathfrak{h}$ and the starting data set $(A_\diamond, \varphi_\diamond, \mathfrak{a}_{\diamond 3})$ using the formula in the three bullets of Proposition 4.1. The task is to solve (4.3) and (4.4) for the map $r \to u|_r$ (for small, non-zero r) with $(A, \varphi, \mathfrak{a}_3)$ as functions of the map $u$ as just described, and with $\mathfrak{X} = e^{-r}\mathfrak{X}_\diamond$.

The important point for the purposes at hand is that (4.3) at any given value of r has the schematic form

$$-\Delta_\diamond u - [A_i, \nabla_{A_\diamond i} u] - [\hat{\mathfrak{a}}_a, [\mathfrak{a}_{\diamond a}, u]] - [A_i, [A_i, u]] - [\hat{\mathfrak{a}}_a [\hat{\mathfrak{a}}_a, u]] + \nabla_{A_\diamond i} [A_i, u] + [\mathfrak{a}_{\diamond a}, [\hat{\mathfrak{a}}_a, u]] = e^{-r} \mathfrak{X}_\diamond$$

(4.35)

with $\{A_i: i \in \{t, 1, 2\}\}$ being the ad(P) components of $A - A_\diamond$ and $\{\hat{\mathfrak{a}}_a: a \in \{1, 2, 3\}\}$ being those of $\mathfrak{a} - \mathfrak{a}_\diamond$ which, when written in terms of $\mathfrak{h}$ using the bullets in Proposition 4.1 have the properties that are listed below in (4.36). (There is an implicit sum in (4.35) for both instances of the repeated index i as it ranges over the set $\{t, 1, 2\}$, and for both instances of the repeated index a as it ranges over the set $\{1, 2, 3\}$.)

- *The* ad(P) *sections from the sets* $\{A_i\}_{i \in \{t,1,2\}}$ *and* $\{\hat{\mathfrak{a}}_a\}_{a \in \{1,2,3\}}$ *at any given point are canonical $\mathbb{R}$-linear functions of the components of the vector* $(\nabla_{A_\diamond} \mathfrak{h}, [\mathfrak{a}_\diamond, \mathfrak{h}])$ *and the components of $\mathfrak{h}$ at that point with norms obeying*

$$|A| + |\hat{\mathfrak{a}}| \le c_0 (1 + |\mathfrak{h}|)(|\nabla_{A_\diamond} \mathfrak{h}| + |[\mathfrak{a}_\diamond, \mathfrak{h}]|).$$

- *The norm of $\nabla_{A_\Diamond i} A_i + [\mathfrak{a}_{\Diamond a}, \hat{\mathfrak{a}}_a]$ obeys the bound*

$$|\nabla_{A_\Diamond i} A_{i\,i} + [\mathfrak{a}_{\Diamond a}, \hat{\mathfrak{a}}_a]| \le (1+|\mathfrak{h}|)(|\Delta_\Diamond \mathfrak{h}| + |[\mathfrak{X}_\Diamond, \mathfrak{h}]| + |\nabla_{A_\Diamond} \mathfrak{h}|^2 + |[\mathfrak{a}_\Diamond, \mathfrak{h}]|^2) \; .$$

(4.36)

The preceding observations are derived by writing the instance of A and $\mathfrak{a}$ on the left hand side of (4.3) as $A_\Diamond + A$ and $\mathfrak{a}_3$ as $\mathfrak{a}_\Diamond + \hat{\mathfrak{a}}$ and then writing A and $\hat{\mathfrak{a}}$ in terms of $\mathfrak{h}$ as indicated in the three bullets of the proposition. Note in this regard that $\mathfrak{h}$ and $\mathfrak{h}^{-1}$ have the same norm because $\mathfrak{h}$ has determinant 1.

With regards to the third bullet's bound: The appearance here of $\Delta_\Diamond \mathfrak{h}$ (as opposed to a generic linear combination of second derivatives of $\mathfrak{h}$) is a happy consequence of the formulas for A and $\hat{\mathfrak{a}}$ in the three bullets of Proposition 4.1. To derive the appearance of $\Delta_\Diamond \mathfrak{h}$, act on both sides of the top bullet in Proposition 4.1 (which defines $A_t - i\hat{\mathfrak{a}}_3$) by the operator as $-(\nabla_{A_{\Diamond t}} + i[\mathfrak{a}_{\Diamond 3}, \cdot\,])$, and act on both sides of the second bullet (which defines $A_1 + iA_2$) by $-(\nabla_{A_{\Diamond 1}} - i\nabla_{A_{\Diamond 2}})$, and act on both sides of the third bullet (which defines $\hat{\mathfrak{a}}_1 - i\hat{\mathfrak{a}}_2$) by as $[\varphi_\Diamond^*, \cdot\,]$. The anti-Hermitian part of the sum of the left hand sides of these identities is $-\nabla_{A_{\Diamond i}} A_i - [\mathfrak{a}_{\Diamond \alpha}, \hat{\mathfrak{a}}_a]$, Meanwhile, the part of the right hand side that has second derivatives of $\mathfrak{h}$ or double commutators of components of $\mathfrak{a}_{\Diamond a}$ on $\mathfrak{h}$ can be written in terms of $\Delta_\Diamond \mathfrak{h}$ and $[\mathfrak{X}_\Diamond, \mathfrak{h}]$.

The important point with regards to (4.36) is that (4.34) can be used to write the $A_\Diamond$–covariant derivatives of $\mathfrak{h}$ and the commutators of $\mathfrak{h}$ with $\mathfrak{a}_\Diamond$ in terms of integrals of corresponding $A_\Diamond$-covariant derivatives of $u$ and $\mathfrak{a}_\Diamond$ commutators with $u$ at values of r between the given value and $r = 0$. As a consequence, the terms in (4.35) with A and $\hat{\mathfrak{a}}$ are formally quadratic or higher order as functions of $u$ and its covariant derivative and commutators. This fact suggests that perturbation theoretic techniques can be brought to bear to find the desired solution map $r \to u|_r$ for values of r near $r = 0$. And, this is the approach taken.

*Part 2*: To set up the perturbation theory, fix for the moment a number $\varepsilon \in (0, 1]$ (an upper bound is given subsequently) and let *C* denote the Banach space of differentiable maps from the interval $(-\varepsilon, \varepsilon)$ into the Banach space $\mathbb{B}$ which are identically zero at the $r=0$ point in $(-\varepsilon, \varepsilon)$. (Remember that $\mathbb{B}$ here is defined using $(A_\Diamond, \varphi_\Diamond, \mathfrak{a}_{\Diamond 3})$.) The norm of a map $r \to x|_r$ in *C* is denoted by $\|x\|_C$ and it is defined by the rule

$$\|x\|_C \equiv \sup_{r \in (-\varepsilon, \varepsilon)} \| \tfrac{\partial}{\partial r} x \|_\mathbb{B}$$

(4.37)

Because maps in *C* are zero at $r = 0$, the map at any given value of r obeys the bound

$$\|x|_r\|_{\mathbb{B}} \leq |r| \|x\|_C .$$
(4.38)

An element $x$ in $C$ is said to be a *smooth* element if $x$ and $\frac{\partial}{\partial r} x$ at any given value of r are smooth sections of ad(P). It is sufficient in the constructions that follow to consider smooth elements in $C$ (the constructions that follow work in general--but it is less of a notational hassle to assume that a given element $x$ from $C$ is smooth because then the $(0,\infty)\times\mathbb{R}$ derivatives of $x$ and $\frac{\partial}{\partial r} x$ are defined at each point.)

Let $x$ denote a smooth element in $C$. Because the pointwise norm of $\frac{\partial}{\partial r} x$ on $(0,\infty)\times\mathbb{R}^2$ at any given value of $r \in (-\varepsilon, \varepsilon)$ is bounded by $\|x\|_C$ (see (4.16)), it follows that $x$ defines a map from $(-\varepsilon, \varepsilon)$ to the space of sections of $P\times_{\mathrm{Ad}(SU(2))}\mathrm{Sl}(2;\mathbb{C})$ by setting $u$ to be $\frac{\partial}{\partial r} x$ in (4.34). What with (4.35)–(4.36), that formula leads to the pointwise bounds in in the upcoming (4.39) for $x$'s version of $\mathfrak{h}$ if r is such that $|r|\|x\|_C \leq 1$. These bounds hold for any point $q \in (0,\infty)\times\mathbb{R}^2$.

- $|[\mathbb{I} - \mathfrak{h}|_r](q)| \leq c_0 r \|x\|_C$.
- $\displaystyle\int_{(0,\infty)\times\mathbb{R}^2} (1+G_q)(|\nabla_{A_\Diamond}\mathfrak{h}|_r|^2 + |[\mathfrak{a}_\Diamond, \mathfrak{h}|_r]|^2) + \Big(\int_{(0,\infty)\times\mathbb{R}^2} G_q |\Delta_{(A_\Diamond,\mathfrak{a}_\Diamond)}\mathfrak{h}|_r|\Big)^2 \leq c_0 r^2 \|x\|_C^2$

(4.39)

(When $\mathfrak{h}$ is viewed as a section of the associated vector bundle to P with fiber the space of 2×2 complex matrices (associated via the adjoint representation of SU(2)), then (4.39) says in effect that the $\mathbb{B}$-norm of $\mathfrak{h}$ is bounded by $c_0 r$ times the $C$-norm of $x$.)

*Part 3*: Let $x$ denote again a smooth element in $C$ and define $\mathfrak{h}$ from $x$ using the $u = \frac{\partial}{\partial r} x$ version of (4.34) as instructed in Part 2. For any given $r \in (-\varepsilon, \varepsilon)$, define $(A|_r, \mathfrak{a}|_r)$ from $\mathfrak{h}|_r$ using the rule from the three bullets of Proposition 4.1. The plan now is to use $x$ to define a new smooth element in $C$ to be denoted by $w$ (sometimes $w[x]$) by first solving the equation below at each r for the corresponding $v|_r$ and then setting $w$ at an given value of r to be the integral of the function $v|_{(\cdot)}$ from 0 to that given value of r.

$$-(\nabla_{At}^2 + \nabla_{A1}^2 + \nabla_{A2}^2)v + [\mathfrak{a}_3, [v, \mathfrak{a}_3]] + [\mathfrak{a}_1, [v, \mathfrak{a}_1]] + [\mathfrak{a}_2, [v, \mathfrak{a}_2]] = e^{-r}\mathfrak{X}_\Diamond .$$
(4.40)

Lemma 4.2 can be used directly (with $\mathfrak{P}_0 = e^{-r}\mathfrak{X}_\Diamond$) to see that there is a unique solution to this equation in the $(A|_r, \mathfrak{a}|_r)$ version of the Banach space $\mathbb{B}$. Even so, for purposes of comparing solutions for different choices of r and of $x$, it proves useful to write (4.40) using $(A_\Diamond, \mathfrak{a}_\Diamond)$ and $\mathfrak{h}$ as done below in (4.41) and then invoke Lemma 4.2 for the $(A_\Diamond, \mathfrak{a}_\Diamond)$ version of $\mathbb{B}$ (which is the operative version in this subsection):

$$-\Delta_{\Diamond} \nu - [A_i, \nabla_{A_{\Diamond i}} u] - [\hat{\mathfrak{a}}_a, [\mathfrak{a}_{\Diamond a}, u]] - [A_i, [A_i, u] + [\hat{\mathfrak{a}}_a [u, \hat{\mathfrak{a}}_a]] - \nabla_{A_{\Diamond i}} [A_i, u] + [\mathfrak{a}_{\Diamond a}, [u, \hat{\mathfrak{a}}_a]] = e^{-r} \mathfrak{X}_{\Diamond}.$$

(4.41)

To be sure: The ad(P) sections $\{A_I\}_{i \in \{t,1,2\}}$ and $\{\hat{\mathfrak{a}}_a\}_{a \in \{1,2,3\}}$ that appear here are defined via $u$'s version of $\mathfrak{h}$ (from (4.34)) by the same rules that define them in (4.35) and (4.36). (As in (4.35), both instances of the repeated index i are sums over the set $\{t, 1, 2\}$, and both instances of the repeated index a are sums over the set $\{1, 2, 3\}$.)

The following lemma says in part that (4.41) can be solved where $|r| \|x\|_C < 1$ with $\nu$ defining a continuous map from $(-\varepsilon, \varepsilon)$ to the space of smooth elements in $\mathbb{B}$.

**Lemma 4.3**: *There exists $\kappa > 1$ which is independent of $(A_{\Diamond}, \mathfrak{a}_{\Diamond})$ and $\varepsilon$, and with the following significance: Suppose that $x$ is from $C$ and $r \in (-\varepsilon, \varepsilon)$ obeys $|r| \|x\|_C \leq 1$. There exists a unique solution to (4.41) in $\mathbb{B}$ at that value of $r$; and this solution $\nu|_r$ obeys:*

$$\|\nu|_r\|_{\mathbb{B}} \leq \kappa |r| \|x\|_C^2 + \kappa \left( \int_{(0,\infty)\times\mathbb{R}^2} (1+t^2) |\mathfrak{X}_{\Diamond}|^2 \right)^{1/2}.$$

*Moreover, the assignment $r \to \nu|_r$ defines a continuous map to $\mathbb{B}$ from the $|r| \|x\|_C < 1$ part of the interval $(-\varepsilon, \varepsilon)$. And, if $x$ and $x'$ are any two elements in $C$, and if $r \in (-\varepsilon, \varepsilon)$ is such that $|r| (\|x\|_C + \|x'\|_C) \leq 1$, then the $\mathbb{B}$-norm of the difference between the respective $x$ and $x'$ versions of $\nu$ obeys*

$$\|\nu[x]|_r - \nu[x']|_r\|_{\mathbb{B}} \leq \kappa |r| (\|x\|_C + \|x'\|_C) \|x - x'\|_C.$$

*Proof of Lemma 4.3*: The existence of a unique solution $\nu|_r$ is an instance of Lemma 4.2 which is applied using $\mathfrak{P}_0, \mathfrak{P}_1$ and $\mathfrak{b}$ and $\mathfrak{c}$ as follows:

- $\mathfrak{P}_0 = e^{-r} \mathfrak{X}_{\Diamond}$.
- $\mathfrak{P}_1 = [A_i, \nabla_{A_{\Diamond i}} u] + [\hat{\mathfrak{a}}_a, [\mathfrak{a}_{\Diamond a}, u]] + [A_i, [A_i, u]] + [\hat{\mathfrak{a}}_a [\hat{\mathfrak{a}}_a, u]]$.
- $\mathfrak{b} = [A_i, u]$ and $\mathfrak{c} = [\hat{\mathfrak{a}}_a, u]$.

(4.42)

It follows from (4.36) and (4.39) that these meet the requirements for use by Lemma 4.2.

The continuity with regards to varying r can be seen in two ways: First, continuity follows from the fact that any given r version of $\nu|_r$ is unique. To elaborate, assume the alternative to continuity: There is a sequence of values of r (denote it by $\{r_i\}_{i \in \mathbb{N}}$) that converges to some given value (call it $r_0$) with the corresponding sequence of $\nu|_r$'s lacking a subsequence that converges in $\mathbb{B}$ to the $r = r_0$ version of $\nu|_r$. To obtain nonsense from this, note first that the sequence of $\mathbb{B}$-norms of the sequence $\{\nu|_{r_i}\}_{i \in \mathbb{N}}$ are uniformly bounded which implies that their $\mathbb{H}$-norms are uniformly bounded. That implies in turn that there is a subsequence with a weak limit in $\mathbb{H}$. Both the weak limit and the $r = r_0$ version of $v|_r$ will obey the $r = r_0$ version of (4.40). This implies that the weak limit *is* the $r = r_0$ version of $v|_r$. This weak convergence can be parlayed into strong convergence in $\mathbb{H}$ and then convergence in $\mathbb{B}$ using arguments that are much like the ones used in Step 5 of the proof of Lemma 4.1. Strong convergence of a subsequence is the desired nonsensical conclusion because it contradicts the initial assumption that characterized the sequence $\{r_i\}_{i \in \mathbb{N}}$.

The second proof uses the polynomial nature of the dependence of the ad(P) sections $\{A_i\}_{i \in \{t,1,2\}}$ and $\{\hat{a}_a\}_{a \in \{1,2,3\}}$ on $\mathfrak{h}$, $\nabla_{A_\Diamond} \mathfrak{h}$ and $[\mathfrak{a}_\Diamond, \mathfrak{h}]$, and on the algebraic dependence of $\mathfrak{h}$ on $\frac{\partial}{\partial r} x$ (the $u = \frac{\partial}{\partial r} x$ version of (4.34)) to see that any given $r$ and $r'$ version of $v|_r - v|_{r'}$ obeys an equation that has the form of the equation depicted in Lemma 4.2 with $Z$ obeying the bound $Z \leq c_0(r - r')(1 + \|x\|_C)$. (The equation is the one that is obtained by subtracting the $r'$ version of (4.41) from the $r$ version.)

To see about the norm of $v[x] - v[x']$, fix $r$ with the indicated constraints and subtract that $v[x]|_r$ version of (4.41) from the $v[x']|_r$ version. The difference of the two equations can be written as an equation for $u \equiv v[x] - v[x']$ that has the form depicted in Lemma 4.2 with $A = A_\Diamond$ and $\mathfrak{a} = \mathfrak{a}_\Diamond$, and with the corresponding value of $Z$ no greater than $c_0 r \|x - x'\|_C (\|x\|_C + \|x'\|_C)$. This bound for $Z$ is again due to the polynomial nature of the dependence of the ad(P) sections $\{A_i\}_{i \in \{t,1,2\}}$ and $\{\hat{a}_a\}_{a \in \{1,2,3\}}$ that appear in (4.2) on $\mathfrak{h}$, $\nabla_{A_\Diamond} \mathfrak{h}$ and $[\mathfrak{a}_\Diamond, \mathfrak{h}]$, and then the algebraic dependence of $\mathfrak{h}$ on $\frac{\partial}{\partial r} x$.

*Part 4*: Lemma 4.3 can be summarized as follows: Let L denote the square root of the $(0, \infty) \times \mathbb{R}^2$ integral of the function $(1 + t^2)|\mathfrak{X}_\Diamond|^2$. There exists $c > 1$ (but less than $c_0$) such that the following is true if $\varepsilon < \frac{1}{cL}$. Let $C_L$ denote the radius L ball in $C$ about the origin. Given any element $x$ from this ball and any $r \in (-\varepsilon, \varepsilon)$, there is a solution in $\mathbb{B}$ to the corresponding version of (4.41). These solutions vary continuously with $r$ and thus define a continuous map from $(-\varepsilon, \varepsilon)$ to $\mathbb{B}$. Let $v[x]$ denote this map and let $w[x]$ denote the corresponding map from $(-\varepsilon, \varepsilon)$ to $C$ whose $r$ derivative is $v[x]$. Then, the assignment $x \to w[x]$ maps $C_L$ to itself as a contraction mapping.

With the preceding understood, it then follows that the map $x \to w[x]$ has a unique fixed point in $C_L$. This fixed point is denoted by $x_*$ and its $r$-derivative is denoted by $u_*$. This $u_*$ satisfies the requirements of Proposition 4.1 (it obeys (4.3) with $(A, \mathfrak{a})$ and $\mathfrak{h}$ as described by the proposition) except possibly for two issues that are discussed momentarily.

With regards to constructing the fixed point $x_*$: As explained directly, this fixed point can be obtained as the limit in $\mathbb{B}$ of a sequence of smooth elements in $C_L$ (which is why it is sufficient up to this point to consider only smooth elements): Start with $x_1 \equiv 0$, then set $x_2 = w[x_1]$, then set $x_3 = w[x_2]$, and so on. Since $x_0$ is smooth, so is $x_1$, and then so is $x_2$, and so on. This sequence is a Cauchy sequence in $C_L$ because of the instance of Lemma 4.2 that says in effect (if $c > c_0$) that

$$\|x_k - x_{k-1}\|_C \leq \tfrac{1}{2} \|x_{k-1} - x_k\|_C$$

(4.43)

for $k > 1$. Iterating this inequality leads to the bound $\|x_k - x_{k-1}\|_C \leq c_0 2^{-k} L$.

With regards to the two as yet unresolved issues: The first is a concern about whether $u_*$ at any given value of r is a smooth section of ad(P). One is tempted in this regard to mumble something to the effect that $u_*$ being a smooth section of ad(P) at any given value of r is proved using 'standard' elliptic bootstrapping arguments because it satisfies (4.3) which is a uniformly elliptic equation. Or (to better intimidate the reader) say something to the effect that $u_*$ being smooth is an instance of theorems in Chapters 5 and 6 of C. B. Morreys' book [M] by virtue of (4.3) being uniformly elliptic. Or even better: The fixed point $u_*$ is smooth, see Morrey [M].

To be honest, it is not clear to the author that 'standard' elliptic bootstrapping or a theorem in Morrey's book can be brought to bear, at least initially. The problematic issue is that of starting the bootstrapping: If $|\nabla_{A_\diamond} u_*|$ is known a priori to be a locally bounded function, then indeed the bootstrapping is a standard operation (and there are theorems in Morrey's book to quote—Theorem 5.6.2 for example). The issue is getting from what is known initially about $\nabla_{A_\diamond} u_*$ given only that $u_*$ is in $\mathbb{B}$ (and hence (4.3) holds in a sort of distributional sense) to knowing that the norm of $|\nabla_{A_\diamond} u_*|$ is locally bounded. The next subsection addresses the differentiablility question. The lemma that follows states formally the result.

**Lemma 4.4**: *There exists $\kappa \geq 1$ (which is independent of $(A_\diamond, \varphi_\diamond, \mathfrak{a}_{\diamond 3})$) with the following significance: If the number c is greater than $\kappa$ and if $\varepsilon < \tfrac{1}{cL}$, then the unique fixed point of the map $x \to w[x]$ in $C_L$ is a smooth section of* ad(P) *at each value of* r.

This lemma is proved in Section 4e.

The second issue to address concerns the pointwise asymptotics of $u$ (see (4.4) where $t \to 0$, $t \to \infty$ and $\tfrac{t}{|z|} \to 0$. The following lemma formally asserts that the pointwise asymptotic constraints in (4.4) are satisfied.

**Lemma 4.5**: *There exists $\kappa \geq 1$ (which is independent of $(A_\diamond, \varphi_\diamond, \mathfrak{a}_{\diamond 3})$) with the following significance: If the number $c$ is greater than $\kappa$ and if $\varepsilon < \frac{1}{cL}$, then both the unique fixed point of the map $x \to w[x]$ in $C_L$ and its r-derivative limit uniformly to zero as $t \to 0$, $t \to \infty$ and $\frac{t}{|z|} \to 0$.*

This lemma is proved in Section 4f.
  These two lemmas complete the proof Proposition 4.1.

### e) The proof of Lemma 4.4

The four parts of this subsection prove the following: If the number $c$ from Part 4 of the previous subsection is greater than $c_0$, and if $\varepsilon < \frac{1}{cL}$, then the unique fixed point to the map $x \to w[x]$ in $C_L$ has the following property:

> *The respective $A_\diamond$–covariant derivatives of $x_*$ and $u_* = \frac{\partial}{\partial r} x_*$ at any given $r \in (-\varepsilon, \varepsilon)$ have bounded norm on compact subsets of $(0, \infty) \times \mathbb{R}^2$.*

(4.44)

As noted in the paragraph prior to Lemma 4.4, this conclusion is sufficient to initiate a standard elliptic bootstrapping arguments that prove Lemma 4.4.
  The proof of the assertion in (4.44) has six parts. By way of a look ahead, the proof exploits the fact that $x$ is the contraction mapping fixed point. Also, by way of a look ahead: The complexity of the arguments that follow are due for the most part to the fact that $(0, \infty) \times \mathbb{R}^2$ is non-compact.
  With regards to proving (4.44): It is sufficient to prove that $u_*$ has locally bounded $A_\diamond$–covariant derivative because that implies directly–by integrating with respect to r–that this is also the case for $x_*$.

  *Part 1*: Since $(A_\diamond, \varphi_\diamond, \mathfrak{a}_{\diamond 3})$ is assumed to be smooth, there exists, for any given positive number $\delta$, a smooth function on $(0, \infty) \times \mathbb{R}^2$ to be denoted by $\mathfrak{d}$ with the property that if $q \in (0, \infty) \times \mathbb{R}^2$, then

- $\mathfrak{d}(q) \leq \frac{1}{100} t_q$ *if* $t_q \leq 1$ *and* $\mathfrak{d}(q) \leq \frac{1}{100}$ *if* $t_q \geq 1$.
- $\mathfrak{d}(q)^2 (|B_{A_\diamond 3}| + |E_{A_\diamond 3}| + |\nabla_{A_\diamond} \mathfrak{a}_\diamond| + |\mathfrak{a}_\diamond|^2) \leq \delta^2$ *on the ball of radius $\mathfrak{d}(q)$ centered at* q.

(4.45)

The fact that $\mathfrak{d}$ depends on $\delta$ is left implicit in the notation because $\delta$ can eventually be fixed to have size $\mathcal{O}(c_0^{-1})$.

If δ is sufficiently small (less than $c_0^{-1}$), then the radius $\mathfrak{d}(q)$ ball centered at any given point $q \in (0,\infty) \times \mathbb{R}^2$ is well inside the 'Uhlenbeck' radius ball. For the present purposes, this implies that Karen Uhlenbeck's gauge fixing theorem [U] can be brought to bear in the radius $4\mathfrak{d}(q)$ ball: There exists an isomorphism from the product principle SU(2) bundle over this ball to P that pulls $A_\diamond$ back as a connection that can be written as $\theta_0 + \mathfrak{a}_\diamond$ with $\theta_0$ denoting the product connection and with $\mathfrak{a}_\diamond$ being a Lie algebra valued 1-form that obeys the conditions

- $\frac{\partial}{\partial t} A_{\diamond t} + \frac{\partial}{\partial z_1} A_{\diamond 1} + \frac{\partial}{\partial z_2} A_{\diamond 2} = 0$ *on the whole of the radius* $4\mathfrak{d}(q)$ *ball centered at* q.
- $|A_\diamond| \leq c_0 \mathfrak{d}(q) \sup_B |F_{A_\diamond}|$ *on the concentric, radius* $2\mathfrak{d}(q)$ *ball centered at* q.

(4.46)

The subsequent notation has B denoting the radius $2\mathfrak{d}(q)$ ball centered at q.

*Part 2*: Use the bundle isomorphism on B from Part 1 to write $\Delta_\diamond$ on B as

$$\Delta_\diamond = \Delta + 2A_\diamond \cdot \nabla + Q(\cdot)$$

(4.47)

with with $\Delta$ denoting the standard Laplacian (as before) and with $A_\diamond \cdot \nabla$ and Q as follows:

$$A_\diamond \cdot \nabla \equiv [A_{\diamond t}, \tfrac{\partial}{\partial t}(\cdot)] + [A_{\diamond j}, \tfrac{\partial}{\partial z_j}(\cdot)] \quad and \quad Q(\cdot) \equiv [A_{\diamond t},[A_{\diamond t}, \cdot\,]] + [A_{\diamond j},[A_{\diamond j}, \cdot\,]] + [\mathfrak{a}_{\diamond a},[\mathfrak{a}_{\diamond a},(\cdot)]] \,.$$

(4.48)

(The notation has the two instance of the repeated index j being summed over the set {1,2} and the repeated index a being summed over the set {1, 2, 3}. The point of writing $\Delta_\diamond$ as in (4.47) is this: If δ is small (less than $c_0^{-1}$), then the Dirichelet Green's function of $\Delta_\diamond$ on B is well approximated by the Dirichelet Green's function for the standard Laplacian $\Delta$. In particular, if $\mathfrak{P}$ is a given section of ad(P) on B, then adequate bounds near q for solutions to the equation $-\Delta_\diamond(\cdot) = \mathfrak{P}$ can be had using $\Delta$'s Green's function.

To elaborate, fix a point p in the radius $\mathfrak{d}(q)$ ball with center q and let $G_{B,p}$ denote the Dirichelet Greens function on B for $\Delta$ with pole at p. For the present purposes, the important point about $G_{B,p}$ is that it is positive on B and, supposing that $p' \in B$, then

$$G_{B,p}(p') \leq c_0 \tfrac{1}{|p-p'|} \quad and \quad |\nabla G_{B,p}|(q) \leq \tfrac{1}{|p-p'|^2} \,.$$

(4.49)

Let $\varpi$ now denote the bump function $\chi(2 \frac{|(\cdot)-q|}{\mathfrak{d}(q)} - 1)$. This function is equal to 1 where the distance to q is less than $\tfrac{1}{2}\mathfrak{d}(q)$ and zero where it is greater than $\mathfrak{d}(q)$. Suppose that ν a section of ad(P) over B that obeys the equation $\Delta_\diamond v = -\mathfrak{P}$. Multiply both sides of this equation by $\varpi$ and use the fact that $\varpi \Delta v = \Delta(\varpi v) - 2\nabla\varpi \cdot \nabla v - \Delta\varpi\, v$ and that $\varpi v$ has compact support in B to write ν at any point in the radius $\tfrac{1}{2}\mathfrak{d}(q)$ ball centered at q as

$$v|_p = \int_B G_{B,p}(\varpi \mathfrak{P} + 2\varpi A_\diamond \cdot \nabla v + Q(\varpi v) + 2\langle \nabla \varpi, \nabla v \rangle + \Delta \varpi \, v) \ .$$

(4.50)

The section $v$ is smooth since $\mathfrak{P}$ is smooth (which is assumed). As a consequence, the preceding identity can be differentiated to obtain a corresponding identity for $\nabla v$ where the distance to q is less than $\frac{1}{2} \eth(q)$. These identities for $v$ and $\nabla v$ can be used to bound $\nabla_{A_\diamond} v$ because $\nabla v$ differs in B from $\nabla_{A_\diamond} v$ by at most $c_0 \delta \frac{1}{\eth(q)} |v|$. (See the second bullet of (4.46) and the second bullet of (4.45).)

*Part 3*: Use the identity in (4.50) and its derivative with the bounds in the lower bullet of (4.45) and those in (4.49) to see that

$$|\nabla v|(p) \le c_0 \int_{|q-(\cdot)|<\eth(q)} \frac{1}{|p-(\cdot)|^2} \left(|\mathfrak{P}| + \frac{(1+\varepsilon)}{\eth(q)} |\nabla v|\right) + c_0 \frac{1}{\eth(q)} \sup_B |v|$$

(4.51)

if the distance from p to q is less than $\frac{1}{2} \eth(q)$.

Now let p´ denote a point with distance at most $\frac{1}{4} \eth(q)$ from q. Multiply both sides of (4.51) by $\frac{1}{|p´-p|^2}$ and integrate with respect to p over the radius $\frac{1}{4} \eth(q)$ ball centered at p´. Doing that leads to this identity:

$$\int_{|p´-(\cdot)|<\frac{1}{4}\eth(q)} \frac{1}{|p´-(\cdot)|^2} |\nabla v| < c_0 \int_B \frac{1}{|p´-(\cdot)|} |\mathfrak{P}| + c_0 \left(\int_B \frac{1}{|p´-(\cdot)|} |\nabla v|^2\right)^{1/2} + c_0 \sup_B |v| \ .$$

(4.52)

(This is because the integral with respect to p of $\frac{1}{|p´-p|^2} \frac{1}{|p-(\cdot)|^2}$ over the whole of $\mathbb{R}^3$ is no greater than $c_0 \frac{1}{|p´-(\cdot)|}$.)

Looking ahead: A crucial point to note with regards to (4.52) is this: If $v$ is from the Banach space $\mathbb{B}$, then the both of the terms with $v$ that appear on the right hand side of (4.52) are bounded by $c_0 \|v\|_\mathbb{B}$. Remember in this regard that bounds for $|\nabla v|$ and $|v|$ in B lead to bounds for $\nabla_{A_\diamond} v$ because the latter differs from $\nabla v$ in B by a term with norm at most $c_0 \delta \frac{1}{\eth(q)} |v|$. (This is by virtue of the second bullets in (4.46) and (4.45).)

Supposing that $v$ is from $\mathbb{B}$, then taking p to be q (4.51) and using (4.52) with p´ being the point q leads to this bound for $|\nabla v|$ at q:

$$|\nabla v|(q) \le c_0 \int_{|q-(\cdot)|<\mathfrak{d}(q)} \frac{1}{|q-(\cdot)|^2}|\mathfrak{P}| + c_0 \frac{1}{\mathfrak{d}(q)} \|v\|_{\mathbb{B}} \,.$$

(4.53)

Granted what was said to the effect that $|\nabla_{A_\diamond} v - \nabla v| \le c_0 \mathfrak{d}(q)|v|$, the right hand side of this also bounds the norm at q of $|\nabla_{A_\diamond} v|$ with a possibly larger version of $c_0$.

*Part 4*: Fix for the moment a number $c > 100$ so that if $\varepsilon < \frac{1}{cL}$, then the map $x \to w[x]$ is a contraction mapping on the corresponding version of $C_L$. Supposing that $x$ is a smooth element in $C_L$, then the equation in (4.41) for $v[x]$ at a given value of r has the form $-\Delta_\diamond v|_r = \mathfrak{P}$. As a consequence, the analysis in Part 3 and (4.53) in particular can be used to bound $|\nabla v|$ and $|\nabla_{A_\diamond} v|$ at any given point q. In this instance, $\mathfrak{P}$ is a sum to be denoted by $\mathfrak{Q}_1 + \mathfrak{Q}_2 + \mathfrak{Q}_3$ with their norms obeying the bounds below in (4.54) (see (4.36) in this regard). By way of notation, $u$ is $\frac{\partial}{\partial r} x$ and I(r) denotes the interval in $(-\varepsilon, \varepsilon)$ between r and 0.

- $|\mathfrak{Q}_1| \le c_0 |\mathfrak{X}_\diamond|(1+|u|)$.
- $|\mathfrak{Q}_2| \le c_0 |r| \sup_{r' \in I(r)} |[\mathfrak{a}_\diamond, u|_{r'}]|^2$.
- $|\mathfrak{Q}_3| \le c_0 |r| \sup_{r' \in I(r)} (|\nabla_{A_\diamond} u|_{r'}|^2 + |\Delta_\diamond u|_{r'}||u|_{r'}|)$.

(4.54)

With (4.54) understood, then (4.53) leads to this bound:

$$|\nabla_{A_\diamond} v|_r|(q) \le c_0(1+\|x\|_C) \int_{|q-(\cdot)|<\mathfrak{d}(q)} \frac{1}{|q-(\cdot)|^2}|\mathfrak{X}_\diamond| + c_0 \int_{|q-(\cdot)|<\mathfrak{d}(q)} \frac{1}{|q-(\cdot)|^2}(|\mathfrak{Q}_2| + |\mathfrak{Q}_3|) + c_0 \frac{1}{\mathfrak{d}(q)} \|x\|_C \,.$$

(4.55)

With regards to the $|\mathfrak{X}_\diamond|$ integral on the right hand side of (4.55): With no added assumptions on $\mathfrak{X}_\diamond$, one has only the $|\mathfrak{X}_\diamond| \le c_0 \frac{1}{\mathfrak{d}(q)^2} \delta^2$ bound from (4.45) which leads to a $c_0 \frac{1}{\mathfrak{d}(q)}$ bound for the left most integral on the right hand side of (4.55).

With regards to the $|\mathfrak{Q}_2|$ integral on the right hand side of (4.55): Since the integrand is bounded by $c_0 \frac{1}{\mathfrak{d}(q)^2} \delta^2 r \|x\|_C^2$ (see (4.45)), this integral can't be larger than $c_0 \frac{1}{\mathfrak{d}(q)} r \|x\|_C^2$ which is at most $c_0 \frac{1}{\mathfrak{d}(q)} \frac{1}{c} \|x\|_C$ because $\|x\|_C \le L$ and $|r| < \varepsilon$ and $\varepsilon$ is at most than $\frac{1}{cL}$.

*Part 5*: The problematic terms in (4.55) is the $|\mathfrak{Q}_3|$ integral on the right hand side. To make something of the right most integral in (4.63), suppose that there exists a constant $M > 1$ such that

$$\sup_{r \in (-\varepsilon, \varepsilon)} |\nabla_{A_\diamond} u| < M \tfrac{1}{\partial(q)} \quad \text{and} \quad \sup_{r \in (-\varepsilon, \varepsilon)} |\Delta_\diamond u| \leq c_0 M \tfrac{1}{\partial(q)^2} \ .$$

(4.56)

at any given $q \in (0,1) \times \mathbb{R}^2$. With this bound assumed, return now to (4.63). Fix $\mu \in (0,1)$ for the moment and decompose the $|\mathfrak{Q}_3|$ integral on the right hand side of (4.55) as a sum of the respective integrals where the distance to q is less than $\mu \partial(q)$ and where it is greater than $\mu \partial(q)$. These respective integrals are bounded by

$$c_0 \mu \, |r| (M^2 + M \|x\|_C) \tfrac{1}{\partial(q)} \quad \text{and} \quad c_0 \tfrac{1}{\mu} |r| \|x\|_C^2 \tfrac{1}{\partial(q)} \ .$$

(4.57)

Minimizing their sum over $\mu$ (take $\mu \sim \tfrac{1}{M} \|x\|_C$) leads to this:

$$\int_{|q - (\cdot)| < \partial(q)} \tfrac{1}{|q - (\cdot)|^2} |\mathfrak{Q}_3| \leq c_0 \, |r| (M + \|x\|_C) \|x\|_C \tfrac{1}{\partial(q)} \ .$$

(4.58)

Therefore, since $\|x\|_C < L$ and $|r| < \varepsilon$ which is less than $\tfrac{1}{cL}$, if $c > 100 c_0$, then (4.58) with what has been said already about the other terms in (4.55) lead to the $A_\diamond$-derivative bound

$$|\nabla_{A_\diamond} v[x]|_r |(q) \leq \tfrac{1}{100} M \tfrac{1}{\partial(q)} + c_0 (1 + L) \tfrac{1}{\partial(q)}$$

(4.59)

There is one more critical observation here, which is that $|\Delta_\diamond v[x]|$ can also be bounded in terms of M and $\partial(q)$. Indeed, such a bound comes directly from the equation $\Delta_\diamond v[x] = \mathfrak{Q}_1 + \mathfrak{Q}_2 + \mathfrak{Q}_3$ using the bounds for these $\mathfrak{Q}$'s in (4.54) with (4.56) and (4.45):

$$|\Delta_\diamond v[x]|(q) \leq c_0 \tfrac{1}{\partial(q)^2} (1 + L) + c_0 |r| M(M + \|x\|_C) \ .$$

(4.60)

This bound and the bound in (4.59) are the crucial inputs to the proof of (4.44).

*Part 6*: To finish the proof of (4.44): Remember that $u_*$ is the r-derivative of the fixed point in $C_L$ of the map $x \to w(x))$ and that the fixed point $x_*$ is the limit in $\mathbb{B}$ of the sequence $\{x^{(0)} \equiv 0, x^{(1)}, \ldots\}$ where $x^{(1)}$ is $w(x_0)$ and $x^{(2)}$ is $w(x^1)$ and so on. Thus, the version of v for $x^{(k)}$ (this is $v[x^k]$) obeys (4.59) and (4.60) with M given by the version of (4.56) with $u = v[x^{k-1}]$. Therefore, each successive $x^{(k)}$ has a version of M (to be denoted by $M_{(k)}$) which, according to (4.59) and (4.60), is at most

$$M_{(k)} \leq \tfrac{1}{100} (1 + |r| M_{(k-1)}) M_{(k-1)} + c_0 (1 + L).$$

(4.61)

Now suppose that $|r|M_j$ is less than 1 for all $j < k$. Then interating this bound leads to a bound of the form $M_k \leq c_0(1+L)$. Hence, $rM_k$ is also less than 1 if $\varepsilon L < c_0^{-1}$ which will happen if $c > c_0$ since $\varepsilon < \frac{1}{cL}$. Assume henceforth that this lower bound for $c$ holds.

Supposing this lower bound for $c$, and noting that $M_1 \equiv 0$ since $x_0 \equiv 0$, it follows that $|r|M_1$ is less than 1, and thus so is $|r|M_2$, and thus $|r|M_3$ and so on. Granted that, then

$$|\nabla_{A_\diamond} v[x_k]|_r|(q) \leq c_0(1+L)\frac{1}{\partial(q)}$$
(4.62)

holds for all k and for all $r \in (-\varepsilon, \varepsilon)$ and all $q \in (0,\infty)\times\mathbb{R}^2$. This implies in turn that $|\nabla_{A_\diamond} u_*|_r|(q)$ obeys this same upper bound for all $r \in (-\varepsilon, \varepsilon)$ and $q \in (0,\infty)\times\mathbb{R}^2$.

**f) Proof of Lemma 4.5**

The proof is given only for the case of $u_*$ because the case for $x_*$ follows by integrating with respect to r. There are three parts to the proof.

*Part 1*: Since $u_*$ is smooth, it obeys the differential equation in (4.40) which says (in the terminology of Section 4c) that $-\Delta_{(A,\mathfrak{a})}u_* = e^{-r}\mathfrak{X}_\diamond$ with $(A,\mathfrak{a})$ as described in Proposition 4.1. The point to make here is that $u_*$ is an element in $(A,\mathfrak{a})$'s version of the Banach space $\mathbb{B}$. This is because the $(A,\mathfrak{a})$ version of $\mathbb{B}$ and the $(A_\diamond, \mathfrak{a}_\diamond)$ version contain the same elements, their norms being commensurate. The fact that their norms are commensurate follows from (4.39) and the fact that the sup-norm of any element in either version of $\mathbb{B}$ is bounded by $c_0$ times that version's $\mathbb{B}$-norm (see in particular (4.16)).

*Part 2*: Because $u_*$ is in the $(A,\mathfrak{a})$ version of $\mathbb{B}$, the bound in (4.16) can be invoked which in this case says that

$$|u_*|(q) \leq c_0 e^{-r} \int_{(0,\infty)\times\mathbb{R}^2} G_q |\mathfrak{X}_\diamond|.$$
(4.63)

Granted this inequality, then (4.11) can be invoked to see that

$$|u_*|(q) \leq c_0 \min(\sqrt{t_q}, \frac{1}{\sqrt{t_q}}(1+|\ln t_q|^{1/2}) e^{-r}(\int_{(0,\infty)\times\mathbb{R}^2}(1+t^2)|\mathfrak{X}_\diamond|^2)^{1/2}.$$
(4.64)

This bound implies in turn that $|u_*|$ limits uniformly to zero as $t \to 0$ and as $t \to \infty$.

*Part 3*: This part considers the $\frac{t}{|z|} \to 0$ asymptotics of $|u_*|$. To this end, suppose first that $q \in (0,\infty) \times \mathbb{R}^2$ and that $\frac{t_q}{|z_q|} < 1$. Either $t_q < (\frac{t_q}{|z_q|})^{1/2}$ or not. If so, then (4.64) applies and therefore

$$|u_*|(q) \le c_0 (\frac{t_q}{|z_q|})^{1/4} e^{-r} ( \int_{(0,\infty) \times \mathbb{R}^2} (1+t^2)|\mathcal{X}_\diamond|^2 )^{1/2} .$$

(4.65)

Now suppose that $t_q \ge (\frac{t_q}{|z_q|})^{1/2}$. This implies in particular that $|z_q| \ge (\frac{|z_q|}{t_q})^{1/2}$ so $|z_q|$ is large if $\frac{t_q}{|z_q|}$ is small. Keeping this in mind, break the integral on the right hand side in (4.63) into the region where $|z| < \frac{1}{2} |z_q|$ and the complementary region. If $|z| < \frac{1}{2} |z_q|$, then $|z - z_q| > \frac{1}{2} |z_q|$ in which case the third bullet in (4.10) can be used to bound the integral of $G_q |\mathcal{X}_\diamond|$ over this region by $c_0$ times the expression on the right hand side of (4.65). Meanwhile, the integral of $G_q |\mathcal{X}_\diamond|$ over the $|z| \ge \frac{1}{2} |z_q|$ region is at most

$$c_0 e^{-r} ( \int_{(0,\infty) \times \{z \in \mathbb{R}^2 : |z| \ge \frac{1}{2} (\frac{|z_q|}{|t_q|})^{1/2}\}} (1+t^2)|\mathcal{X}_\diamond|^2 )^{1/2}$$

(4.66)

which limits to zero uniformly as $\frac{t_q}{|z_q|} \to 0$.

## 5. The $s \to T$ limit of $\{\mathfrak{g}|_s\}_{s \in [0,T)}$

Fix a data set $(A_\diamond, \varphi_\diamond, \mathfrak{a}_{3\diamond})$ that obeys the first three bullets of (1.3) and such that the corresponding $\mathcal{X}_\diamond$ obeys Proposition 4.1's requirement that $(1+t^2)|\mathcal{X}_\diamond|^2$ has finite integral on $(0,\infty) \times \mathbb{R}^2$. For example, this could be a data set from Section 3. Let T denote either a positive number or $\infty$, and suppose that there exists a smooth map $s \to \mathfrak{g}|_s$ from the half-open interval $[0, T)$ to the space of sections of $P \times_{Ad(SU(2))} Sl(2; \mathbb{C})$ that is described by (4.2)–(4.4) with $\mathfrak{g}|_{s=0}$ being the identity. The upcoming Proposition 5.1 asserts (in part) that the 1-parameter family of sections $\{\mathfrak{g}|_s\}_{s \in [0,T)}$ of $P \times_{Ad(SU(2))} Sl(2; \mathbb{C})$ has a smooth limit as $s \to T$.

By way of terminology in the proposition: The phrase "$\{\mathfrak{g}|_s\}_{s \in [0,T)}$ extends smoothly to $[0, T]$" in the case when T is $\infty$ means that the $s \to \infty$ limit of $\{\mathfrak{g}|_s\}_{s \in [0,\infty)}$ exists as a section of $P \times_{Ad(SU(2))} Sl(2; \mathbb{C})$ and that the associated reparameterized family $\{\mathfrak{g}|_{s=-\ln(1-r)}\}_{r \in [0,1]}$ defines a smooth map from $[0,1]$ to the space of sections of $P \times_{Ad(SU(2))} Sl(2; \mathbb{C})$. This reparametrization is used implicitly when T is $\infty$.

The upcoming proposition introduces a function $q \to \mu_q$ on $(0,\infty) \times \mathbb{R}^2$ that is defined as follows: The value of $\mu_{(\cdot)}$ at any given point $q = (t_q, z_q) \in (0,\infty) \times \mathbb{R}^2$ is the largest number that is less than the smaller of $\frac{1}{t_q}$ and $\frac{1}{100}$ subject to the constraint

$$\mu_q{}^2 t_q{}^2 (|B_{A_\Diamond 3}| + |E_{A_\Diamond 3}| + |\nabla_{A_\Diamond} \mathfrak{a}_\Diamond| + |\mathfrak{a}_\Diamond|^2) \leq 1$$

(5.1)

at all points in the radius $\frac{1}{4} t_q$ ball centered at q. By way of comparison, any given positive $\delta$ version of the function $\mathfrak{d}$ from (4.45) can be taken to be $\delta \mu_q t_q$.

**Proposition 5.1**: *Suppose that $(A_\Diamond, \varphi_\Diamond, \mathfrak{a}_{\Diamond 3})$ with its corresponding $\mathfrak{X}_\Diamond$ and the 1-parameter family $\{\mathfrak{g}|_s\}_{s \in [0,T)}$ are as described in the opening paragraph of this section. Then the family $\{\mathfrak{g}|_s\}_{s \in [0,T)}$ extends to $[0, T]$ to define a smooth map from the closed interval $[0, T]$ into the space of smooth sections of $P \times_{\mathrm{Ad}(SU(2))} \mathrm{Sl}(2; \mathbb{C})$. Moreover, there exist $\kappa > 1$ such that*

$$|\mathfrak{g}|_s| < \kappa \quad \text{and} \quad |(\nabla_{A_\Diamond} \mathfrak{g})|_s| < \frac{\kappa}{\mu_{(\cdot)}^\kappa t} \;.$$

*at all points in $(0, \infty) \times \mathbb{R}^2$ and for all $s \in [0, T]$. And, given $\varepsilon > 0$, there exists $\kappa_\varepsilon > 1$ such that if $(t, z)$ is a point in $(0, \infty) \times \mathbb{R}^2$ with $t < \frac{1}{\kappa_\varepsilon}$ or $t > \kappa_\varepsilon$ or $\frac{t}{|z|} < \frac{1}{\kappa_\varepsilon}$, then at that point,*

$$|\mathfrak{g}|_s - \mathbb{I}| < \varepsilon \quad \text{and} \quad |(\nabla_{A_\Diamond} \mathfrak{g})|_s| < \varepsilon \frac{1}{\mu_{(\cdot)}^\kappa t} \;.$$

*for all $s \in [0, T]$.*

Sections 5b–5f contain the proof of this proposition. Section 5a contains a preliminary result that is used in the proof

Proposition 5.1 has the corollary that follows about the family $\{(A|_s, \varphi|_s, \mathfrak{a}_3|_s)\}_{s \in [0,T)}$ which is defined by (4.1) from the family $\{\mathfrak{g}|_s\}_{s \in [0,T)}$.

**Proposition 5.2**: *Suppose that $(A_\Diamond, \varphi_\Diamond, \mathfrak{a}_{\Diamond 3})$ with its corresponding $\mathfrak{X}_\Diamond$ and the 1-parameter family $\{\mathfrak{g}|_s\}_{s \in [0,T)}$ are as described in the opening paragraph of this section. Define the corresponding 1-parameter data set $\{(\mathcal{A}|_s \equiv A|_s - A_\Diamond, \varphi|_s, \mathfrak{a}_3|_s)\}_{s \in (0,T)}$ via (4.1) from the family $\{\mathfrak{g}|_s\}_{s \in [0,T)}$. Then the family $\{(\mathcal{A}|_s, \varphi|_s, \mathfrak{a}_3|_s)\}_{s \in [0,T)}$ extends to $[0, T]$ to define a smooth, 1-parameter family of $\mathrm{ad}(P)$ valued 1-form, section of $\mathrm{ad}(P) \otimes_{\mathbb{R}} \mathbb{C}$ and section of $\mathrm{ad}(P)$. Moreover if the function $\mu_{(\cdot)}$ from Proposition 5.1 is bounded on $(0, \infty) \times \mathbb{R}^2$, then the following is true: There exists $\kappa > 1$ such that the norms of $|\mathcal{A}|_s|$ and $|\mathfrak{a}|_s - \mathfrak{a}_\Diamond|$ are bounded by $\kappa \frac{1}{t}$ on the whole of $(0, \infty) \times \mathbb{R}^2$ and for any $s \in [0, T]$. Moreover, given $\varepsilon > 0$, there exists $\kappa_\varepsilon > 1$ such that the norms of any $s \in [0, T]$ version of $\mathcal{A}|_s$ and $|\mathfrak{a}_3|_s - \mathfrak{a}_{\Diamond 3}|$ are at most $\varepsilon \frac{1}{t}$ on the parts of $(0, \infty) \times \mathbb{R}^2$ where either $t < \frac{1}{\kappa_\varepsilon}$ or $t > \kappa_\varepsilon$ or $\frac{t}{|z|} < \frac{1}{\kappa_\varepsilon}$. Meanwhile, the norm of any $s \in [0, T]$ version of $\varphi|_s - \varphi_\Diamond$ is at most $\varepsilon |\varphi_\Diamond|$ on these same parts of $(0, \infty) \times \mathbb{R}^2$. In particular, if $\mu_{(\cdot)}$ is bounded on $(0, \infty) \times \mathbb{R}^2$ and if $(A_\Diamond, \varphi_\Diamond, \mathfrak{a}_{\Diamond 3})$ is described by (1.10), then so is each member of the family $\{(A|_s, \varphi|_s, \mathfrak{a}_{3s}|\}_{s \in [0,T]}$.*

By way of an example: Any small δ version of a data sets from Section 3 is described by (1.10) and its associated function $\mu_{(\cdot)}$ is bounded. See Lemmas 3.3 and 3.4 in this regard.

*Proof of Proposition 5.2*: By virtue of the definition of the family in (4.1), the convergence of $\{\mathfrak{g}|_s\}_{s \in [0,T)}$ as described by Proposition 5.1 leads directly to the manner of convergence of $\{(A|_s, \varphi|_s, \mathfrak{a}_3|_s)\}_{s \in [0,T)}$ asserted by Proposition 5.2.

A notational convention: In the context of discussing $\{\mathfrak{g}|_s\}_{s \in [0,,T)}$ and/or the corresponding family of data sets $\{(A|_s, \varphi|_s, \mathfrak{a}_3|_s)\}_{s \in [0,T)}$ from (4.1), this section uses on occasion $c_\diamond$ to denote a number that is greater than 1 and independent of T and any chosen $s \in [0,T]$ and any chosen point in $(0, \infty) \times \mathbb{R}^2$. However, $c_\diamond$ can depend on global (which means not point in $(0,\infty) \times \mathbb{R}^2$ dependent) properties of the initial $(A_\diamond, \varphi_\diamond, \mathfrak{a}_{\diamond 3})$. For example, $c_\diamond$ can depend on the $(0,\infty) \times \mathbb{R}^2$ integral of $(1+t^2)|\mathfrak{X}_\diamond|^2$. As with the ubiquitous $c_0$, the number $c_\diamond$ can be assumed to increase between successive appearances.

**a) The section *u* and the Banach space $\mathbb{B}$**

In the context of Proposition 5.1, any given $s \in [0, T)$ version of $(A, \varphi, \mathfrak{a}_3)$ along the flow (they are defined via (4.1) from $\{\mathfrak{g}|_s\}_{s \in [0, T)}$) has a corresponding version of the Banach space $\mathbb{B}$. And, the corresponding $u$ that is defined by (4.2)–(4.4) at that value of s must be in this space with its norm obeying

$$\|u\|_{\mathbb{B}} \leq c_0 e^{-s} \left( \int_{(0,\infty) \times \mathbb{R}^2} (1+t^2)|\mathfrak{X}_\diamond|^2 \right)^{1/2}.$$

(5.2)

This assertion that $u$ is in the $(A, \varphi, \mathfrak{a}_3)$ version of $\mathbb{B}$ and the inequality in (5.2) is the central observation in this subsection.

To prove this observation: First invoke Lemma 4.2 with $\mathfrak{P}_0$ being $-e^{-s}\mathfrak{X}_\diamond$ and $\mathfrak{P}_1$ and $\mathfrak{b}$ and $\mathfrak{c}$ being zero to obtain a solution to (4.3) that is in the $(A, \varphi, \mathfrak{a}_3)$ version of $\mathbb{B}$. Denote this solution by $u_*$. A replay of the arguments from Lemma 4.5's proof proves that $|u_*|$ limits uniformly to zero as $t \to 0$, as $t \to \infty$ and as $\frac{t}{|z|} \to 0$ (see (4.64)–(4.66).)

Since both $|u_*|$ and $|u|$ limit to zero uniformly as $t \to 0$, as $t \to \infty$ and as $\frac{t}{|z|} \to 0$, this is also the case for $|u - u_*|$. Keeping this in mind: Note that $u - u_*$ obeys the version of (4.3) with $\mathfrak{X} \equiv 0$; and so it follows that $|u - u_*|$ obeys $\Delta|u - u_*| \geq 0$ which implies via the maximum principle that $u = u_*$ because of the aforementioned asymptotic limits of $|u - u_*|$.

### b) $C^0$ convergence and an outline of the proof

This subsection proves that $\{\mathfrak{g}|_s\}_{s\in[0,T]}$ converges as $s \to T$ in the $C^0$ topology on the space of sections of $P\times_{Ad(SU(2))} Sl(2;\mathbb{C})$ and that the various members of the family $\{\mathfrak{g}|_s\}_{s\in[0,T]}$ have the asserted asymptotic limits to the identity section.

With regards to convergence: Because (4.16) holds for elements in $\mathbb{B}$, the sup-norm of an element is no greater than its $\mathbb{B}$-norm. Since $u$ at any given $s \in [0, T)$ obeys (5.2) and because the $(0,\infty)$-integral of $e^{-s}$ is absolutely convergent, the integral of the norm of $u$ on the interval $[0, T]$ at any fixed point in $(0,\infty)\times\mathbb{R}^2$ is finite with an upper bound that is independent of the point and the number $T$. This implies that $\mathfrak{g}|_s$ for $s\in[0,T)$ when viewed as an $SU(2)$-equivariant map from $P$ to $Sl(2;\mathbb{C})$ maps $P$ into an $s$-independent, compact subset of $Sl(2;\mathbb{C})$. Moreover, it implies that the $s \to T$ limit of the family $\{\mathfrak{g}|_s\}_{s\in[0,T)}$ (which I will call $\mathfrak{g}_T$) is a $C^0$ section of $P\times_{Ad(SU(2))} Sl(2;\mathbb{C})$ because the space of $C^0$ sections is closed with respect to the sup-norm.

With respect to the asymptotics: Because $u = u_*$ in (4.64)–(4.66), it follows that this section $\mathfrak{g}_T$ limits uniformly to the identity section $\mathbb{I}$ as $t \to 0$ and as $t \to \infty$ and as $\frac{t}{|z|} \to 0$ with an upper bound for the rate of decay determined a priori by integrals of the function $(1+t^2)|\mathfrak{X}_\Diamond|^2$ over the indicated parts in (4.64)–(4.66) of $(0,\infty)\times\mathbb{R}^2$.

As a parenthetical remark for now: The $C^0$ convergence as $s \to T$ of $\{\mathfrak{g}|_s\}_{s\in[0,T)}$ leads directly to the $C^0$ convergence as $s \to T$ of $\{\varphi|_s\}_{s\in[0,T)}$ because $\varphi|_s$ at any given value of $s \in [0,T)$ is $(\mathfrak{g}|_s)\,\varphi_\Diamond(\mathfrak{g}|_s)^{-1}$. Moreover, what was said in the preceding paragraph about the asymptotics of any given $s \in [0, T]$ version of $\mathfrak{g}|_s$ implies the following: Given $\varepsilon > 0$, there exists $c_\varepsilon > 1$ such that for any $s \in [0, T]$, the bound $|\varphi|_s - \varphi_\Diamond| < \varepsilon|\varphi_\Diamond|$ holds where $t < \frac{1}{c_\varepsilon}$ or $t > c_\varepsilon$ or $\frac{t}{|z|} < \frac{1}{c_\varepsilon}$. The analogous conclusions apply for the family $\{\lambda|_s\}_{s\in[0,T)}$ with regards to convergence as $s \to T$ and with regards to their asymptotics relative to $\lambda_\Diamond$. (Remember that $\lambda$ is $\sigma - 4q$ and that each $s \in [0,T]$ version is $(\mathfrak{g}|_s)^{-1}\lambda_\Diamond(\mathfrak{g}|_s)$.)

OUTLINE OF THE PROOF OF PROPOSITION 5.1:

With $C^0$ convergence understood, the remainder of the proof of Proposition 5.1 amounts to a step by step bootstrapping argument that successively strengthens the convergence norm for the family $\{\mathfrak{g}|_s\}_{s\in[0,T]}$. The input for this bootstrapping is the equation in (4.3) for $u$, the equations in (4.6) that describe how $A$ and $\mathfrak{a}_3$ change with $s$, and the identity $\varphi = \mathfrak{g}^{-1}\varphi_\Diamond \mathfrak{g}$. The equations in the first three bullets of (1.3) are also used as is the defining equation for $\mathfrak{X}$ (which is $\nabla_{A_t}\mathfrak{a}_3 - B_{A_3} - \frac{i}{2}[\varphi, \varphi^*]$). The arguments that follow present only the bootstrapping up to $C^1$ convergence because the subsequent steps going from $C^1$ to $C^2$ and so on are along the same lines as this initial $C^0$ to $C^1$ step.

The first step to going from $C^0$ convergence to $C^1$ convergence is to prove $L^2$ convergence for $\{A|_s\}_{s\in[0,T)}$ and $\{\mathfrak{a}_3|_s\}_{s\in[0,T)}$ as $s \to T$ in balls in $(0,\infty)\times\mathbb{R}^2$ with compact

closure and then use that convergence to obtain s and T-independent bounds for the integrals of $|\nabla_{A_\diamond} u|^2$ over these same balls. This is done in Section 5c. The subsequent Section 5d uses some of what is done in Section 5c to derive s-independent $L^2$ bounds for the curvature of A and the A-covariant derivative of $\mathfrak{a}$ on balls with compact closure. These are used in Section 5e to prove $L^2_1$ convergence for $\{A|_s\}_{s\in[0,T)}$ and $\{\mathfrak{a}_3|_s\}_{s\in[0,T)}$ as $s \to T$ on these balls. The latter convergence is used to obtain T-independent bounds for the integrals over these balls of $|\nabla_{A_\diamond}\nabla_{A_\diamond} u|^2$. Section 5f takes Section 5e's bounds as input and derives $C^1$ convergence on compact sets in $(0,\infty)\times\mathbb{R}^2$ and the asserted pointwise asymptotics for the $A_\diamond$-covariant derivative of $\mathfrak{g}_T$.

### c) $L^2$ convergence of $\{A|_s\}_{s\in[0,T)}$ and $\{\mathfrak{a}_3|_s\}_{s\in[0,T)}$ and $L^2_1$ bounds for $u$

This subsection has five parts which derive various bounds for the integrals on balls in $(0,\infty)\times\mathbb{R}^2$ of the square of the norms of $\mathfrak{a}_3|_s$ and $A|_s \equiv A|_s - A_\diamond$ and $(\nabla_{A_\diamond} u)|_s$.

*Part 1*: Fix a point $q \in (0,\infty)\times\mathbb{R}^3$ and let B denote the ball of radius $\frac{1}{4} t_q$ centered at q. By virtue of (4.6)'s second bullet, if s and s´ are from [0, T) with s > s´, then the difference between the corresponding $\mathfrak{a}_3|_s$ and $\mathfrak{a}_3|_{s'}$ obeys the bound

$$\int_B |\mathfrak{a}_3|_s - \mathfrak{a}_3|_{s'}|^2 \le c_0 (\int_{s'}^s (\int_B |(\nabla_A u)|_{(\cdot)}|^2)^{1/2})^2 .$$

(5.3)

Use the definition of the $\|\cdot\|_\mathbb{B}$ norm to bound the integral over B on the right hand side of (5.3) by the product of $c_0 \min(t_q, 1)$ and $\|u|_{(\cdot)}\|_\mathbb{B}^2$. Then invoke (5.1) and (4.23) to see that

$$\int_B |\mathfrak{a}_3|_s - \mathfrak{a}_3|_{s'}|^2 \le c_0 \min(t_q, 1)(e^{-s} - e^{-s'})^2 \int_{(0,\infty)\times\mathbb{R}^2} (1+t^2)|\mathfrak{X}_\diamond|^2 .$$

(5.4)

This bound implies directly that $\{\mathfrak{a}_3|_s\}_{s\in[0,T)}$ converges as $s \to T$ on B in the Hilbert space that is obtained by completing the space of smooth sections of ad(P) on B using the norm whose square is the B-integral of the square of the pointwise norm.

Much the same argument using the third, fourth and fifth bullets of (4.6) lead to an analogous integral bound for $A|_s - A|_{s'}$ and the analogous manner of convergence for the $s \to T$ limit of $\{A|_s\}_{s\in[0,T)}$ which is that

$$\int_B |A|_s - A|_{s'}|^2 \le c_0 (\int_{s'}^s (\int_B (|(\nabla_A u)|_{(\cdot)}|^2 + |[\mathfrak{a}_3, u]|_{(\cdot)}|^2))^{1/2})^2 .$$

(5.5)

This in turn leads to the bound

$$\int_B |A|_s - A|_{s'}|^2 \leq c_0 \min(t_q, 1)(e^{-s} - e^{-s'})^2 \int_{(0,\infty)\times\mathbb{R}^2} (1+t^2)|\mathfrak{X}_\diamond|^2 \ .$$

(5.6)

With regards to convergence: It follows directly from (5.6) that the 1-parameter family $\{A|_s\}_{s\in[0,T)}$ converges as $s \to T$ on B in the Hilbert space completion of the space of smooth, ad(P)-valued 1-forms on B using the norm whose square is the B-integral of the square of the pointwise norm.

*Part 2*: More can be said with regards to an upper bound for the right hand side integrals in (5.3) and (5.5) when $t_q$ is small relative to 1. To elaborate: If $t_q < 1$, then the integral over B of any $s \in [0, T)$ version of the function $|\nabla_A u|^2 + |[\mathfrak{a}, u]|^2$ has the bound

$$\int_B (|\nabla_A u|^2 + |[\mathfrak{a}, u]|^2) \leq c_0 e^{-2s} t_q^2 \int_{(0,\infty)\times\mathbb{R}^2} (1+t^2)|\mathfrak{X}_\diamond|^2 \ .$$

(5.7)

This comes from the $v \equiv u$ and $\mathfrak{P}_0 \equiv e^{-s}\mathfrak{X}_\diamond$ version of (4.26) which in this instance says

$$\tfrac{1}{2}|u|^2(q) + \int_{(0,\infty)\times\mathbb{R}^2} G_q(|\nabla_A u|^2 + |[\mathfrak{a}, u]|^2) \leq \int_{(0,\infty)\times\mathbb{R}^2} G_q |u| |\mathfrak{P}_0| \ .$$

(5.8)

To derive (5.7) from this, restrict the left hand side integral to B and then use the bound $G_q(p) \geq \tfrac{1}{c_0} \tfrac{1}{|p-q|}$ which holds in B (see (4.10)). Meanwhile, bound $|u|$ from above in the right hand side integral using the $\mathbb{B}$-norm bound in (5.2) and by bounding the integral of $(1+t^2)G_q^2$ from above by (4.11). The bounds in (5.3), (5.5), (5.7) lead to the $t_q < 1$ bound:

$$\int_B |\mathfrak{a}_3|_s - \mathfrak{a}_3|_{s'}|^2 + \int_B |A|_s - A|_{s'}|^2 \leq c_0 t_q^2 (e^{-s} - e^{-s'})^2 \int_{(0,\infty)\times\mathbb{R}^2} (1+t^2)|\mathfrak{X}_\diamond|^2$$

(5.9)

whose right hand side is a factor of $t_q$ smaller than those of (5.3) and (5.6).

*Part 3*: Here is an analog of (5.9) when $\tfrac{t_q}{|z_q|}$ is small: If $\tfrac{t_q}{|z_q|} \leq 1$, then

$$\int_B (|\nabla_A u|^2 + |[\mathfrak{a}, u]|^2) \leq c_0 e^{-2s} t_q \left( \left(\tfrac{t_q}{|z_q|}\right)^{1/2} \int_{(0,\infty)\times\mathbb{R}^2} (1+t^2)|\mathfrak{X}_\diamond|^2 + \int_{(0,\infty)\times\{|z|>\frac{1}{8}(\frac{|z_q|}{t_q})^{1/2}\}} (1+t^2)|\mathfrak{X}_\diamond|^2 \right).$$

(5.10)

(This assertion is proved momentarily.) The inequality in (5.16) leads to a corresponding analog of (5.15) that says this: Given $\varepsilon > 0$, there exists $r_\varepsilon > 1$ which is independent of s (it depends only on the initial $(A_0, \varphi_0, \mathfrak{a}_{03})$) such that

$$\int_B ||\mathfrak{a}_3|_s - \mathfrak{a}_3|_{s'}|^2 + \int_B ||A|_s - A|_{s'}|^2 \leq c_0 \varepsilon^2 t_q (e^{-s} - e^{-s'})^2 \int_{(0,\infty) \times \mathbb{R}^2} (1+t^2)|\mathfrak{X}_0|^2 \quad \text{where } |z_q| > r_\varepsilon t_q .$$

(5.11)

The proof of (5.10) introduces a new bump function whose definition follows directly: Fix $q \in (0, \infty) \times \mathbb{R}^2$ and let B again denote the ball of radius $\frac{1}{4} t_q$ centered at q. Now let $\iota_q$ denote the function with compact support on B that is given by the rule

$$\iota_q(\cdot) = \chi(\tfrac{16|(\cdot) - q|}{t_q} - 1).$$

(5.12)

This function is equal to 1 where the distance to q is less than $\frac{1}{16} t_q$ and equal to zero where that distance is greater than $\frac{1}{8} t_q$.

With $\iota_q$ in hand, here is how to prove (5.10): There are two cases to consider: The first is $t_q \leq (\frac{t_q}{|z_q|})^{1/2}$, in which case (5.9) can be invoked and (5.9) implies (5.10). In the second case, $t_q \geq (\frac{t_q}{|z_q|})^{1/2}$ in which case $|z| \geq (\frac{|z_q|}{t_q})^{1/2}$. In this case, take the inner product of both sides of (4.3) with $\iota_q^2 u$ and then integrate both sides of the result over B. Having done that, then integrate by parts on the left hand side to obtain this:

$$\int_B \iota_q^2 (|\nabla_A u|^2 + |[\mathfrak{a}, u]|^2) \leq c_0 t_q \sup_B |u|^2 + c_0 e^{-2s} \int_B (1+t^2)|\mathfrak{X}_0|^2$$

(5.13)

This inequality with the $u_* = u$ version of the bound in (4.66) for $|u_*|$ leads directly to (5.10) with a different $c_0$ when $t_q \geq (\frac{t_q}{|z_q|})^{1/2}$.

*Part 4*: It is past time to introduce short-hand notation for subsequent use. To this end, fix $q \in (0,\infty) \times \mathbb{R}^2$ and let B again denote the ball of radius $\frac{1}{4} t_q$ centered at q. Now introduce $\mathcal{X}_q$ to denote the number

$$\mathcal{X}_q \equiv \sup_{s \in [0,T)} e^s (\sup_B |u|^2 + \tfrac{1}{t_q} \int_B (|\nabla_A u|^2 + |[\mathfrak{a}, u]|^2))^{1/2} .$$

(5.14)

By way of an explanation for the factor of $\frac{1}{t_q}$ multiplying the integral on the right hand side of (5.14): This is present so that the two terms on the right hand side of (5.14) scale the same under constant coordinate rescalings of the form $(t, z) \to (\lambda t, \lambda z)$.)

It follows from the definition of the $\|\cdot\|_{\mathbb{B}}$ norm and (4.16) and (5.2) that

$$\mathcal{X}_q^2 \leq c_0 \int_{(0,\infty)\times\mathbb{R}^2} (1+t^2)|\mathfrak{X}_\diamond|^2 \ .$$

(5.15)

However, as illustrated by (5.9) and (4.64) when $t_q \ll 1$ or by (5.10) and (4.65)–(4.66) when $\frac{t_q}{|z_q|} \ll 1$, the number $\mathcal{X}_q$ can be much less than the integral on the right hand side of (5.15). Also, in the event that $t_q > 1$, the number $\mathcal{X}_q$ is smaller than $c_0 \frac{1}{t_q}|\ln t_q|$ times the right hand side of (5.15) which is a tautological consequence of the definition of the $\mathbb{B}$-norm (for the integral part of $\mathcal{X}_q$) and the bound in (4.64) for the $\sup_{\mathbb{B}}|u|$ part.

By way of an sample application of the notation: The conclusions of Section 5b can be summarized as follows: There are T independent upper bounds for $u$ and for both $\mathfrak{g}$ and $\mathfrak{g}^{-1}$ on $(0,\infty)\times\mathbb{R}^2$; and moreover,

- $|u|_s \leq c_0 e^{-s} \mathcal{X}_q$
- $|\mathfrak{g}|_s - \mathfrak{g}|_{s'}| \leq c_\diamond |e^{-s} - e^{-s'}| \mathcal{X}_q$

(5.16)

at any point q in $(0,\infty)\times\mathbb{R}^2$ and for any numbers s and s´ from $[0, T)$.

*Part 5*: With $\mathcal{X}_q$ so defined, then the inequalities in (5.3) and (5.5) lead to this:

$$\int_B ||\mathfrak{a}_3|_s - \mathfrak{a}_3|_{s'}|^2 + \int_B ||A|_s - A|_{s'}|^2 \leq c_0 t_q (e^{-s} - e^{-s'})^2 \mathcal{X}_q^2 \ .$$

(5.17)

As explained momentarily, this last inequality with (5.16) leads to the bound

$$\int_B |\nabla_{A_\diamond} u|_s|^2 \leq c_\diamond t_q e^{-2s} \mathcal{X}_q^2 \quad \text{if } s \in [0, T) \ .$$

(5.18)

Indeed, this follows because $\nabla_{A_\diamond} u = \nabla_A u - [A, u]$ and because (5.17) leads directly to a $c_\diamond t_q$ bound for the integral over B of $|A|^2$. (Remember also that $|u| \leq \|u\|_{\mathbb{B}}$.)

By way of a parenthetical remark: The bound in (5.18) and (5.16) lead to the following bounds for $\mathfrak{g}$:

- $\int_B |\nabla_{A_\diamond} \mathfrak{g}|_s|^2 \leq c_0 t_q \mathcal{X}_q^2 \quad \text{if } s \in [0, T)$

- $\int_B |\nabla_{A_\diamond} \mathfrak{g}|_s - \nabla_{A_\diamond} \mathfrak{g}|_{s'}|^2 \leq c_0 t_q (e^{-s} - e^{-s'})^2 \mathcal{X}_q^2$ *for any* s, s' $\in [0, T)$.

(5.19)

This is because $\frac{\partial}{\partial s} \nabla_{A_\diamond} \mathfrak{g} = \nabla_{A_\diamond} u \, \mathfrak{g} + u \nabla_{A_\diamond} \mathfrak{g}$. (The top bullet in (5.19) follows directly from this with the bound in (5.18) for $u|_s$ and the lower bullet follows from this given the top bullet in (5.19) with that same bound for $u|_s$.)

With regards to ramifications: Let $\mathbb{M}(2;\mathbb{C})$ denote the vector space of $2 \times 2$ complex matrices. Given that $\{\mathfrak{g}|_s\}_{s \in [0,T)}$ converges as $s \to T$ in the $C^0$ topology on B, the preceding bounded, implies that the family $\{\mathfrak{g}_s\}_{s \in [0,T)}$ also converges as $s \to T$ in the Hilbert space that is obtained by completing the space of sections of $P \times_{Ad(SU(2))} \mathbb{M}(2;\mathbb{C})$ using the norm whose square is the integral over B of the sum of the squares of the section and it's $A_\diamond$-covariant derivative. (This norm is the $L^2_1$ norm.)

### d) Local $L^2$ bounds for the curvature of A and the A-covariant derivative of $\mathfrak{a}$

This subsection invokes the bounds in Section 5c to bound the integral of the square of the norm of the curvature of $A|_s$ and the $A|_s$-covariant derivatives of $\mathfrak{a}|_s$ on balls in $(0, \infty) \times \mathbb{R}^2$.

To start, fix $s \in [0,T)$ and let $(A, \varphi, \mathfrak{a}_3)$ denote $(A|_s, \varphi|_s, \mathfrak{a}_3|_s)$. Now write (4.3) using the definition of $\mathcal{X}$ as $\nabla_{At}\mathfrak{a}_t - B_{A3} - \frac{i}{2}[\varphi, \varphi^*] = e^{-s}\mathcal{X}_\diamond$. Having done that, act on both sides by $\nabla_{At}$ and use the Bianchi identities with the first and third bullets in (1.3) to write the result as:

$$(\nabla_{At}^2 + \nabla_{A1}^2 + \nabla_{A2}^2)\mathfrak{a}_3 - [\mathfrak{a}_1,[\mathfrak{a}_3,\mathfrak{a}_1]] - [\mathfrak{a}_2,[\mathfrak{a}_3,\mathfrak{a}_2]] = e^{-s}\nabla_{At}\mathcal{X}_\diamond .$$

(5.20)

Fix $q \in (0, \infty) \times \mathbb{R}^2$ and take the inner product of both sides of (5.20) with $\iota_q^2 \mathfrak{a}_3$ and integrate the result over B. Integration by parts (one on the left hand side and one on the right hand side) leads to this inequality:

$$\int_B \iota_q^2 (|\nabla_A \mathfrak{a}_3|^2 + |[\mathfrak{a},\mathfrak{a}_3]|^2) \leq c_0 \frac{1}{t_q^2} \int_B |\mathfrak{a}_3|^2 + c_0 e^{-2s} \int_B |\mathcal{X}_\diamond|^2 .$$

(5.21)

As a consequence of this and (5.17), there is an s-independent upper bound:

$$\int_B \iota_q^2 (|\nabla_A \mathfrak{a}_3|^2 + |[\mathfrak{a},\mathfrak{a}_3]|^2) \leq c_0 \frac{1}{t_q^2} \int_B |\mathfrak{a}_{\diamond 3}|^2 + c_0 \frac{1}{(1+t_q^2)} e^{-2s} \int_B (1+t^2) |\mathcal{X}_\diamond|^2 .$$

(5.22)

Because $B_{A3}$ can be written as $\nabla_{At}\mathfrak{a}_3 + \frac{i}{2}[\varphi, \varphi^*] - e^{-s}\mathcal{X}_\diamond$, and because the top bullet of (1.3) identifies $E_A$ with A-covariant derivatives of $\mathfrak{a}_3$, the preceding leads to this:

$$\int_B \iota_q^2(|B_{A3}|^2 + |E_A|^2) \leq c_0 \tfrac{1}{t_q^2}\left(\int_B |\mathfrak{a}_{\diamond 3}|^2 + t_q^2\int_B|\varphi_\diamond|^4\right) + c_0 \tfrac{1}{(1+t_q^2)} e^{-2s}\int_B (1+t^2)|\mathfrak{X}_\diamond|^2 \ .$$

(5.23)

There is also a bound of this sort for the integral over B of $\iota_q^2|\nabla_A\varphi|^2$:

$$\int_B \iota_q^2|\nabla_A\varphi|^2 \leq c_0 \tfrac{1}{t_q^2}\left(\int_B |\mathfrak{a}_\diamond|^2 + t_q^2 \sup_B|\varphi_\diamond|^2 \int_B|\mathfrak{a}_\diamond|^2\right) + c_0 \tfrac{1}{(1+t_q^2)} e^{-2s}\int_B (1+t^2)|\mathfrak{X}_\diamond|^2 \ .$$

(5.24)

To derive the latter inequality, first add the $\nabla_{At}$ derivative of the third equation in (1.3) to the $\nabla_{A1} - i\nabla_{A2}$ derivative of the second and use the definition of $\mathfrak{X}$ to derive this:

$$(\nabla_{At}^2 + \nabla_{A1}^2 + \nabla_{A2}^2)\varphi - [\mathfrak{a}_3, [\varphi, \mathfrak{a}_3]] + \tfrac{1}{2}[[\varphi, \varphi^*], \varphi] - i[e^{-s}\mathfrak{X}_\diamond, \varphi] = 0 \ .$$

(5.25)

To obtain (5.24), take the $\langle\ \rangle$ inner product of the latter identity with $\iota_q^2 \varphi^*$, then integrate the result over B, and then integrate by parts.

By way of notation: The inequalities in (5.11) and (5.22)–(5.24) lead to slightly weaker inequalities that can be expressed in terms of what is denoted by $\mu_q$ in (5.1):

$$\int_B \iota_q^2(|\nabla_A\mathfrak{a}|^2 + |[\mathfrak{a},\mathfrak{a}]|^2 + \tfrac{1}{t_q^2}|\mathfrak{a}|^2 + |B_{A3}|^2 + |E_A|^2) \leq c_0 \mu_q^{-4} \tfrac{1}{t_q} \ .$$

(5.26)

A parenthetical observation is worth making now (it doesn't play a subsequent role). The observation is that the a priori bounds from (5.26) for the integral of the square of the norm of the curvature of A and for the A-covariant derivative of $\mathfrak{a}$ on balls in $(0,\infty)\times\mathbb{R}^2$ can be used in conjunction with Karen Uhlenbeck's compactness theorem in [U] to obtain smooth $s \to T$ limits of subsequences of $\{(A|_s, \varphi|_s, \mathfrak{a}_3|_s)$ after acting by suitably chosen, s-dependent automorphisms of P. The following proposition makes a formal statement to this effect:

**Proposition 5.3**: *Fix a data set $(A_\diamond, \varphi_\diamond, \mathfrak{a}_{3\diamond})$ that obeys the first three bullets of (1.3) and with $\mathfrak{X}_\diamond$ such that the $(0,\infty)\times\mathbb{R}^2$ integral of $(1+t^2)|\mathfrak{X}_\diamond|^2$ is finite. Let T denote either a positive number or $\infty$, and suppose that there is a corresponding version of the family $s \to \mathfrak{g}|_s$ that obeys (4.2)-(4.4) for values of s in the interval $[0,T)$. Then the $[0,T)$ parametrized family of data sets $\{(A|_s, \varphi|_s, \mathfrak{a}_3|_s)\}_{s\in[0,T)}$ that is defined via $\{\mathfrak{g}|_s\}_{s\in[0,T)}$ by (4.1) has the following property: Let $\Omega$ denote any given open set in $(0,\infty)\times\mathbb{R}^2$ with compact closure. There exists a sequence $\{s_k\}_{k\in\mathbb{N}} \subset [0,T)$ converging to T and a corresponding sequence $\{g_k\}_{k\in\mathbb{N}}$ of automorphisms of $P|_\Omega$ such that the sequence* $\{g_k{}^*(A|_{s_k}, \varphi|_{s_k}, \mathfrak{a}_3|_{s_k})\}_{k \in \mathbb{N}}$ *converges in the $C^\infty$ topology on $\Omega$ to a data set $(A, \varphi, \mathfrak{a}_3)$ that obeys the top three bullets in (1.3) with the norm of the corresponding $\mathfrak{X}$ being $e^{-T}|\mathfrak{X}_\diamond|$.*

***Proof of Proposition 5.3***: This is a very brief summary: The assertion that there is convergence on compact subsets of a subsequence in the Sobolev $L^2_1$ weak topology after acting termwise by suitable automorphisms of P follows directly from Karen Uhlenbeck's theorem in [U] given the bounds on the curvature integrals in (5.26) and the covariant derivative integrals in (5.26). Convergence in the stronger topologies follows from the fact that the equations in (5.20) and (5.25) are elliptic equations whose non-linearities are controlled given $L^2_1$ norm bounds on $\mathfrak{a}$ and on the connection A on small radius balls in $\Omega$ (after pull-back by a suitable automorphism of P). The details of the proof are omitted because the proposition isn't used here. (See for example [T1] and [T2] and [GU] for analogous convergence assertions for sequences of solutions to (1.1).)

**e) $L^2_1$ convergence of $\{A|_s\}_{s \in [0, T)}$ and $\{\mathfrak{a}_3|_s\}_{s \in [0, T)}$ and $L^2_2$ bounds for $u$**

The square of the $L^2_2$ norm in question is the integral over B of the sum of the squares of the pointwise norms of the section, its $A_\diamond$-covariant derivative, and the latter's $A_\diamond$-covariant derivative. (Thus, derivatives to order 2 are involved.)

The discussion that follows has four parts.

*Part 1*: This part presents a Sobolev inequality that plays a starring role. To set the stage for this, let q denote for the moment a point in $\mathbb{R}^3$, let $\rho$ denote a positive number, and let $\mathcal{B}$ denote the ball of radius $\rho$ centered at q. Suppose in what follows that $f$ is a function on $\mathcal{B}$ with $|\nabla f|^2$ and $f^2$ having finite integral over $\mathcal{B}$. The function space inequality concern $f$. A preliminary inequality is the Sobolev inequality:

$$(\int_\mathcal{B} f^6)^{1/3} \leq c_0 \int_\mathcal{B} (|\nabla f|^2 + \tfrac{1}{\rho^2}|f|^2) \ .$$

(5.27)

(See e.g. [A] for a proof.) This Sobolev inequality plays a secondary role since it is only used to to bound the integral of the fourth power of $f$:

$$(\int_\mathcal{B} f^4)^{1/2} \leq c_0 (\int_\mathcal{B} (|\nabla f|^2 + \tfrac{1}{\rho^2}|f|^2))^{3/4} (\int_\mathcal{B} |f|^2)^{1/4} \ .$$

(5.28)

The preceding inequality is the crucial one for the purposes at hand.

*Part 2*: Return now to the equation in (4.3) for $u$ to derive bounds for the integrals of $|\nabla_A \nabla_A u|^2$ and $|[\mathfrak{a}_3, [u, \mathfrak{a}_3]]|^2$ over balls in $(0,\infty) \times \mathbb{R}^2$. To this end, fix again a point q in $(0,\infty) \times \mathbb{R}^2$ and reintroduce the ball B centered at q with radius $\frac{1}{4} t_q$. A new bump function is needed with compact support where the function $\iota_q$ is equal to 1. This function is denoted by $o_q$ and it is given by the rule

$$o_q(\cdot) = \chi(\tfrac{32|(\cdot) - q|}{t_q} - 1)$$

(5.29)

Fix $s \in [0, T)$ and multiply both sides of the $u|_s$ version of (4.3) by $o_q$, square the respective norms and then integrate the resulting identity over B. Integration by parts and commuting derivatives leads to this bound:

$$\int_B o_q^2(|\nabla_A \nabla_A u|^2 + |[\mathfrak{a}_3, [u, \mathfrak{a}_3]]|^2) \leq c_0 \tfrac{1}{t_q^2} \int_B (|\nabla_A u|^2 + |[\mathfrak{a}_3, u]|^2)$$
$$+ c_0 \int_B o_q^2(|B_{A3}| + |E_A| + |\nabla_A \mathfrak{a}| + |\mathfrak{a}|^2)|\nabla_A u|^2 + c_0 \int_B o_q^2(|B_{A3}|^2 + |E_A|^2 + |\nabla_A \mathfrak{a}|^2)|u|^2 \ .$$

(5.30)

The Sobolev inequality in (5.28) for the cases $f = |\mathfrak{a}|$ and $f = o_q|\nabla_A u|$ with $\mathcal{B}$ being the ball of radius $\frac{1}{16} t_q$ centered at q leads from (5.30) to the following bound:

$$\int_B o_q^2(|\nabla_A \nabla_A u|^2 + |[\mathfrak{a}_3, [u, \mathfrak{a}_3]]|^2) \leq c_\diamond \tfrac{1}{t_q}(1 + \mu_q^{-c_0})e^{-2s} \mathcal{X}_q^2 \ .$$

(5.31)

To say more about the derivation: The left most integral on the right hand side of (5.30) is bounded directly by a T-independent multiple of the right hand side of (5.31) by virtue of the definition of $\mathcal{X}_q$. Meanwhile, the right most integral on the right hand side of (5.30) is bounded by

$$c_0 \int_B o_q^2(|B_{A3}|^2 + |E_A|^2 + |\nabla_A \mathfrak{a}|^2) \sup_B |u|^2$$

(5.32)

which is bounded by a T-independent multiple of the right hand side of (5.31) using (5.26) with the definition of $\mathcal{X}_q$. As for the middle integral: It is no larger than the product of the square roots of two integrals, these square roots being

$$(\int_B \iota_q^4(|B_{A3}|^2 + |E_A|^2 + |\nabla_A \mathfrak{a}|^2 + |\mathfrak{a}|^4))^{1/2} \quad and \quad (\int_B o_q^4 |\nabla_A u|^4)^{1/2}.$$

(5.33)

(Remember that $o_q$ has its support where $\iota_q$ is equal to 1.) The square root of the integral of $\iota_q^4(|B_{A3}|^2 + |E_A|^2 + |\nabla_A \mathfrak{a}|^2)$ is bounded via (5.26), and the square root of the integral of $\iota_q^4|\mathfrak{a}|^4$ is bounded using (5.28) by

$$c_0 \left( \int_B (\iota_q^2 |\nabla_A \mathfrak{a}|^2 + \tfrac{1}{t_q^2}|\mathfrak{a}|^2) \right)^{3/4} \left( \int_B |\mathfrak{a}|^2 \right)^{1/4},$$

(5.34)

which is bounded in turn using (5.26). Meanwhile, an analogous appeal to (5.28) bounds the right most integral in (5.33) by

$$c_0 \left( \int_B (o_q^2 |\nabla_A \nabla_A u|^2 + \tfrac{1}{t_q^2}|\nabla_A u|^2) \right)^{3/4} \left( \int_B |\nabla_A u|^2 \right)^{1/4}.$$

(5.35)

The important point now is that the integral of $o_q^2|\nabla_A \nabla_A u|^2$ appears here with the exponent ¾ whereas it appears on the left hand side of (5.30) with the larger exponent 1. As a consequence, the middle integral in (5.30) can be bounded by a sum of terms with one being ½ of the integral that appears on the left hand side of (5.30) and the other being a T-independent multiple of what appears on the right hand side of (5.31).

With regards to the preceding two applications of (5.28): Keep in mind that $|\nabla_A e| \geq |\nabla|e||$ for any ad(P)-valued section or 1-form $e$, and thus for $e = \mathfrak{a}$ and for $e = \nabla_A u$ in particular. Also, keep in mind that $\iota_q|\nabla_A e| \leq |\nabla_A(\iota_q e)| + |\nabla \iota_q||e|$ and that $|\nabla \iota_q| \leq c_0 \tfrac{1}{t_q}$.

*Part 3*: The inequality in (5.31) with the Sobolev inequality in (5.28) leads to an apriori bound on the integral of $o_q^4|\nabla_A u|^4$ which leads via (4.1) directly to this bound

$$\left( \int_B o_q^4 |\mathfrak{a}_3|_s - \mathfrak{a}_3|_{s'}|^4 + \int_B o_q^4 |A|_s - A|_{s'}|^4 \right)^{1/2} \leq c_\diamond \tfrac{1}{t_q}(1 + \mu_q^{-c_0})(e^{-s} - e^{-s'})^4 X_q^4.$$

(5.36)

The latter implies in particular that the square root of the integral over B of $o_q^4|A|^4$ has an s-independent upper bound by $c_0$ times what appears on the right hand side of (5.36). And, that bound implies that the integral over B of $o_q^4|A|^2|\nabla_A u|^2$ also has an s-independent upper bound. Therefore, so does $o_q^4|\nabla_{A_\diamond}(\nabla_A u)|^2$ because $o_q^4|\nabla_{A_\diamond}(\nabla_A u)|$ no greater than the sum of $o_q^4|\nabla_A \nabla_A u|$ and $2o_q^4|A||\nabla_A u|$. Then, what with (4.1), the latter leads to this:

$$\left( \int_B o_q^4 |\nabla_{A_\diamond}(\mathfrak{a}_3|_s - \mathfrak{a}_3|_{s'})|^2 + \int_B o_q^4 |\nabla_{A_\diamond}(A|_s - A|_{s'})|^2 \right)^{1/2} \leq c_\diamond \tfrac{1}{t_q}(1 + \mu_q^{-c_0})(e^{-s} - e^{-s'})^2 X_q^2.$$

(5.37)

Note in particular that (5.37) with (5.17) implies that the 1-parameter families $\{\mathfrak{a}_3|_s\}_{s \in [0,T)}$ and $\{A|_s\}_{s \in [0,T)}$ both converge as $s \to T$ in the $L^2_1$ Sobolev topology on open sets in $(0, \infty) \times \mathbb{R}^2$ with compact closure.

*Part 4*: To see about the integral of $|\nabla_{A_\diamond} \nabla_{A_\diamond} u|^2$, start by introducing yet a third bump function, this one is denoted by $\hat{o}_q$ and it is defined by the rule

$$\hat{o}_q(\cdot) \equiv \chi(\tfrac{64|(\cdot)-q|}{t_q} - 1)$$

(5.38)

This function is equal to 1 where the distance to q is less than $\tfrac{1}{64} t_q$ and it has compact support where $o_q$ is equal to 1.

With $\hat{o}_q$ in hand, square the inequality

$$|\nabla_{A_\diamond} \nabla_{A_\diamond} u| \leq |\nabla_A \nabla_A u| + 2|A||\nabla_A u| + (4|A|^2 + 2|\nabla_{A_\diamond} A|)|u|,$$

(5.39)

then multiply both sides of the result by $o_q^2$ and then integrate over the ball B. Doing that bounds the integral over B of $\hat{o}_q^2 |\nabla_{A_\diamond} \nabla_{A_\diamond} u|^2$ by $c_0$ times the sum of the integrals over B of

$$o_q^2 |\nabla_A \nabla_A u|^2 \quad and \quad (o_q|A|)^2 (\hat{o}_q|\nabla_A u|)^2 \quad and \quad (o_q^4 |A|^4 + o_q^2 |\nabla_{A_\diamond} A|^2)|u|^2$$

(5.40)

Note in this regard that replacing $\hat{o}_q$ by $o_q$ is allowed since $\hat{o}_q$ has support where $o_q$ is equal to 1. Now, the integral over B of the left most term in (5.40) is bounded by

$$c_\diamond \tfrac{1}{t_q}(1 + \mu_q^{-c_0}) e^{-2s} \mathcal{X}_q^2.$$

(5.41)

courtesy of (5.31); and the integral of the right most term is bounded by the same courtesy of (5.36) and (5.37) and the definition of $\mathcal{X}_q$. As for the integral of the middle term in (5.40), it is bounded by the product of the square roots of the integrals of $o_q^4 |A|^4$ and $\hat{o}_q^4 |\nabla_A u|^4$; and the former is bounded courtesy of (5.36) and the latter courtesy of $f = \hat{o}_q |\nabla_A u|$ version of the Sobolev inequality in (5.28) and then (5.31). Put all of these bounds together to see that

$$\int_B \hat{o}_q^2 |\nabla_{A_\diamond} \nabla_{A_\diamond} u|^2 \leq c_\diamond \tfrac{1}{t_q}(1 + \mu_q^{-c_0}) e^{-2s} \mathcal{X}_q^2.$$

(5.42)

By way of a parenthetical remark: This bound plus the Sobolev inequality in (5.28) plus (5.16) and (5.18) imply $L^2_2$–norm convergence as $s \to T$ on the support of $\hat{o}_q$ for the family $\{g|_s\}_{s \in (0,T]}$. A precise statement statement is that

$$\int_B \tilde{o}_q^2 |\nabla_{A_\diamond} \nabla_{A_\diamond}(g|_s - g|_{s'})|^2 \le c_\diamond \frac{1}{t_q}(1 + \mu_q^{-c_0})(e^{-s} - e^{-s'})^2 \mathcal{X}_q^2 ,$$

(5.43)

with $\tilde{o}_q$ being a bump function with support where $\hat{o}_q$ is equal to 1. (It is defined by the analog of (5.38) with the factor 64 replaced by 132.)

**f) The $C^1$ convergence of $\{g|_s\}_{s \in [0,T)}$ as $s \to T$ on open sets with compact closure**

The proof for $C^1$ convergence exploits the definition of $u$ as $-i \frac{\partial}{\partial s} g g^{-1}$ to write

$$\tfrac{\partial}{\partial s}(\nabla_{A_\diamond} g g^{-1}) = i \nabla_{A_\diamond} u .$$

(5.44)

Since $g$ and $g^{-1}$ are uniformly bounded (see Section 5b), and since $\{g|_s\}_{s \in [0,T)}$ converges as $s \to T$ as depicted in (5.16), this identity implies first that $|\nabla_{A_\diamond} g|_s|$ for $s \in [0, T)$ has a T-independent upper bound on a given compact set if $\{|\nabla_{A_\diamond} u|_s|\}_{\in [0,T]}$ is bounded on that set; and supposing this, then (5.44) and (5.16) imply that if $s$ and $s'$ are in $[0, T)$, then

$$|\nabla_{A_\diamond} g|_s - \nabla_{A_\diamond} g|_{s'}| \le c_\diamond \int_s^{s'} |\nabla_{A_\diamond} u| + c_\diamond \sup_{[0,T]} |\nabla_{A_\diamond} g| |e^{-s} - e^{-s'}| \mathcal{X}_q .$$

(5.45)

This subsection supplies a bound for $|\nabla_{A_\diamond} u|$ at any given $q \in (0, \infty) \times \mathbb{R}^2$:

$$|\nabla_{A_\diamond} u|(q) \le c_\diamond \frac{1}{t_q}(1 + \mu_q^{-c_0}) e^{-s} \mathcal{X}_q .$$

(5.46)

The latter with (5.45) with (5.16) lead directly to the following convergence assertion:

$$|\nabla_{A_\diamond} g|_s - \nabla_{A_\diamond} g|_{s'}|(q) \le c_\diamond \frac{1}{t_q}(1 + \mu_q^{-c_0}) |e^{-s} - e^{-s'}| \mathcal{X}_q ;$$

(5.47)

and this implies in turn the $s \to T$ convergence of $\{g|_s\}_{s \in [0,T]}$ in the $C^1$ topology on compact sets in $(0, \infty) \times \mathbb{R}^2$.

The four parts of what follows derive the bound in (5.46). Looking ahead, the arguments in these parts of the proof are similar in most respects to those in Section 4e for the proof of Lemma 4.4. In particular, notation for Section 4e is used in what follows.

*Part 1*: A version of Hardy's inequality plays the central role in subsequent arguments. To set the stage: Fix $q \in (0,\infty) \times \mathbb{R}^2$ and $\rho \in (0, \frac{1}{2} t_q)$; and then let $\mathcal{B}$ denote the ball in $(0,\infty) \times \mathbb{R}^2$ with center q and radius $\rho$. The required Hardy's inequality says that if $f$ is any given smooth function on $\mathcal{B}$ (or just in the $L^2_1$ Sobolev space), then

$$\int_{\mathcal{B}} \frac{1}{|(\cdot)-q|^2} f^2 \leq c_0 \int_{\mathcal{B}} (|\nabla f|^2 + \frac{1}{\rho^2}|f|^2) \ .$$

(5.48)

The proof of Hardy's inequality amounts to little more than two applications of integration by parts in spherical coordinates, see the appendix for a derivation.

*Part 2*: Having fixed $\varepsilon \in (0, 10^{-6}]$, define a function on $(0,\infty) \times \mathbb{R}^2$ to be denoted by $\eth$ by the rule $\eth = \varepsilon \mu_{(\cdot)} t_{(\cdot)}$. With the upper bound on $\varepsilon$ understood, then the radius $4\eth$ ball centered at q will be well inside the radius $\frac{1}{512} t_q$ ball centered at q. If $\varepsilon$ is smaller if necessary (but still greater than $c_0^{-1}$), then (4.45) is obeyed and the radius $\eth(q)$ ball about any given point q is well inside the Uhlenbeck radius of that ball for the connection $A_\diamond$. There is, as a consequence, an isomorphism from the product principal SU(2) bundle over this ball to P that pulls $A_\diamond$ back as a connection that can be written as $\theta_0 + A_\diamond$ with $\theta_0$ denoting the product connection and with $A_\diamond$ being a Lie algebra valued 1-form that obeys the conditions in (4.46). The connection A then appears as

$$A = \theta_0 + A_\diamond + \mathcal{A} \ ,$$

(5.49)

where $\mathcal{A}$ here is the pull-back via the isomorphism of the ad(P)-valued 1-form $A - A_\diamond$.

Having written A in this way, then the equation in (4.3) for $u$ can be written schematically as

$$\Delta_\diamond u + [\nabla_{A_\diamond j} \mathcal{A}_j, u] + \yen(u) = -e^{-s} \mathfrak{X}_\diamond \ ,$$

(5.50)

where $\Delta_\diamond$ is depicted in (4.47)-(4.48) and $\yen(u)$ is short-hand for the sum

$$2[\mathcal{A}_j, \nabla_{A_\diamond j} u] + [\mathcal{A}_j, [\mathcal{A}_j, u]] + [\mathfrak{a}_a, [(\mathfrak{a}_a - \mathfrak{a}_{\diamond a}), u]] + [(\mathfrak{a}_a - \mathfrak{a}_{\diamond a}), [\mathfrak{a}_q, u]] + [(\mathfrak{a}_a - \mathfrak{a}_{\diamond a}), [(\mathfrak{a}_a - \mathfrak{a}_{\diamond a}), u]] \ ,$$

(5.51)

with it implicit that the index j is summed over the set {t, 1, 2} and the index a is summed over the set {1, 2, 3}.

To conform to the notation in Section 4d, let B now denote the ball of $2\eth(q)$ centered at q, and, supposing that p is in the concentric, radius $\eth(q)$ ball, let $G_{B,p}$ again denote the Dirichelet Green's function for B with pole at p. The inequalities in (4.49) again provide the relevant information with regards to $G_{B,p}$. Also as in Section 4e, let $\varpi$ denote the bump function $\chi(2\frac{|(\cdot) - q|}{\eth(q)} - 1)$.

Suppose now that $\mathfrak{P}$ is a given section of ad(P) on B. If $v$ is a section of ad(P) on B that obeys the equation $\Delta_\diamond v = -\mathfrak{P}$ on B, then $v$ at a point p with distance less than $\eth(q)$ from q can be written using $G_{B,p}$ and $\varpi$ as done in (4.50). In the relevant instance,

$$\mathfrak{P} = [\nabla_{A_\diamond j} A_j, u] + \yen(u) + e^{-s}\mathfrak{X}_\diamond .$$

(5.52)

The depiction of $v$ in (4.50) leads to the $\theta_0$-covariant derivative bound in (4.51).

*Part 3*: In the present circumstances, the integral on the right hand side of (4.51) has two terms with the right most term being this:

$$\frac{1}{\eth(q)} \int_{|p' - (\cdot)| < \frac{1}{4}\eth(q)} \frac{1}{|p' - (\cdot)|^2} |\nabla u| .$$

(5.53)

To get a useful bound for that integral, use the fact that $|A_\diamond| \le c_0 \frac{1}{\eth(q)}$ to bound it first by

$$c_0 \frac{1}{\eth(q)} \int_{|p' - (\cdot)| < \frac{1}{4}\eth(q)} \frac{1}{|p' - (\cdot)|^2} |\nabla_{A_\diamond} u| + c_\diamond \frac{1}{\eth(q)} \sup_B |u| ,$$

(5.54)

and then by

$$c_0 \frac{1}{\sqrt{\eth(q)}} \left( \int_{|p' - (\cdot)| < \frac{1}{4}\eth(q)} \frac{1}{|p' - (\cdot)|^2} |\nabla_{A_\diamond} u|^2 \right)^{1/2} + c_\diamond \frac{1}{\eth(q)} \sup_B |u| .$$

(5.55)

For the integral here, invoke Hardy's inequality in (5.48) with $f = |\nabla_{A_\diamond} u|$ to bound the preceding integral by integrals of $|\nabla_{A_\diamond} \nabla_{A_\diamond} u|^2$ and $\frac{1}{\eth(q)^2} |\nabla_{A_\diamond} u|^2$; and then use (5.42) for the former and the definition of $\mathfrak{X}_q$ for the latter. Doing this results in this bound:

$$\frac{1}{\eth(q)} \int_{|p' - (\cdot)| < \frac{1}{4}\eth(q)} \frac{1}{|p' - (\cdot)|^2} |\nabla u| \le c_\diamond \frac{1}{t_q}(1 + \mu_q^{-c_0})e^{-s}\mathfrak{X}_q ,$$

(5.56)

which is the desired bound for (5.53).

*Part 4*: The left most integral on the right hand side in (4.51) is the one with $\mathfrak{P}$ which in this case (and with p = q) is the following integral:

$$\int_{|q-(\cdot)|<\mathfrak{d}(q)} \frac{1}{|q-(\cdot)|^2}(|\nabla_{A_\Diamond j}A_j||u| + |¥(u)| + e^{-s}|\mathfrak{X}_\Diamond|)$$

(5.57)

The three terms in this integral are bounded separately in the subsequent paragraphs.

The contribution to (5.57) from the $\mathfrak{X}_\Diamond$ term has the bound:

$$\int_{|q-(\cdot)|<\mathfrak{d}(q)} \frac{1}{|q-(\cdot)|^2} e^{-s}|\mathfrak{X}_\Diamond| \le c_0 e^{-s} \frac{1}{\mu_q^{2/3} t_q^{2/3}} \frac{1}{(1+t_q^2)^{1/3}} (\int_B (1+t^2)|\mathfrak{X}_\Diamond|^2)^{1/3}.$$

(5.58)

This is derived by first fixing a number $r \in (0, 1]$ and then breaking the integral on the left in (5.58) into a sum of two integrals, the part where the distance to q is less than $r\mathfrak{d}(q)$ and the part where the distance is greater than $r\mathfrak{d}(q)$. The former is bounded by $c_0 r \mathfrak{d}(q)^{-1}$ because $|\mathfrak{X}_\Diamond|$ is at most $c_0\mathfrak{d}(q)^{-2}$. The latter is bounded by $c_0$ times the square root of the product of two integrals with the first integral being that of $|(\cdot)-q|^{-4}$ over the domain in question (it is at most $c_0 r^{-1}\mathfrak{d}(q)^{-1}$), and the second being the integral of $|\mathfrak{X}_\Diamond|^2$ over the domain. An appropriate choice for $r$ in the preceding decomposition of the left hand integral in (5.58) gives the bound on the right hand side of (5.58).

The contribution to (5.57) from $¥(u)$ is dealt with by first bounding $|¥(u)|$ by

$$|¥(u)| \le c_0(|A||\nabla_{A_\Diamond} u| + (|A|^2 + |\mathfrak{a}|^2 + |\mathfrak{a} - \mathfrak{a}_\Diamond|^2)\sup_B |u|^2,$$

(5.59)

and then invoking Hardy's inequality in (5.48) with $f = |A|$ using (5.37), with $f = |\nabla_{A_\Diamond} u|$ using (5.42), and then with $f = |\mathfrak{a} - \mathfrak{a}_\Diamond|$ using (5.37) again. Doing that leads to this bound:

$$\int_{|q-(\cdot)|<\mathfrak{d}(q)} \frac{1}{|q-(\cdot)|^2}|¥(u)| \le c_\Diamond \frac{1}{t_q}(1+\mu_q^{-c_0})e^{-s}\mathcal{X}_q(1+\mathcal{X}_q^{c_0}).$$

(5.60)

The left most and final term in (5.57) is the term with $\nabla_{A_\Diamond j}A_j$. To deal with this term, first use (4.1) to see that

$$\tfrac{\partial}{\partial s}(\nabla_{A_\Diamond j}A_j) = [B_3, u] - [A_1, \nabla_{A2} u] + [A_2, \nabla_{A1} u] - [\nabla_{At}\mathfrak{a}_3, u] + [A_t, [\mathfrak{a}_3, u]] - [\mathfrak{a}_3, \nabla_{At} u].$$

(5.61)

which can be written using the definition of $\mathfrak{X}$ (which is $e^{-s}\mathfrak{X}_\Diamond$) as

$$\tfrac{\partial}{\partial s}(\nabla_{A_\Diamond j}A_j) = [-e^{-s}\mathfrak{X}_\Diamond, u] - [\tfrac{i}{2}[\varphi, \varphi^*], u] - [A_1, \nabla_{A2} u] + [A_2, \nabla_{A1} u] + [A_t, [\mathfrak{a}_3, u]] - [\mathfrak{a}_3, \nabla_{At} u].$$

(5.62)

Integrate this formula invoking Hardy's inequality in (5.48) with $f = |A|$ using (5.37), with $f = |\nabla_A u|$ using (5.31), and then with $f = |\mathfrak{a} - \mathfrak{a}_\diamond|$ using (5.37) again. Doing so (and remembering from Section 5b that $|\varphi| \leq c_0 |\varphi_\diamond|$) leads to this bound

$$\int_{|q-(\cdot)|<\vartheta(q)} \frac{1}{|q-(\cdot)|^2} |\nabla_{A_\diamond j} A_j| |u| \leq c_\diamond \frac{1}{t_q} (1 + \mu_q^{-c_0}) e^{-s} \mathcal{X}_q.$$

(5.63)

Taken together, the bounds derived above for (5.57) and in Part 3 for (5.53) from the right hand side of the $w = u$ version of (4.51) lead to the asserted bound in (5.46) for the $A_\diamond$-covariant derivatives of $u$.

## 6. Integrals of $|B_{A3}|^2$ and $|E_{A1}|^2$ and $|E_{A2}|^2$

To summarize where things stand after Sections 4 and 5: Let $(A_\diamond, \varphi_\diamond, \mathfrak{a}_{\diamond 3})$ denote a data set that obeys the first three bullets of (1.3) and with the $(0, \infty) \times \mathbb{R}^2$ integral of $(1+t^2)|\mathfrak{X}_\diamond|^2$ being finite. By virtue of the analysis in Sections 4 and 5, there exists a smooth map $s \to \mathfrak{g}|_s$ from the half-line $[0, \infty)$ to the space of sections of $P \times_{\mathrm{Ad}(SU(2))} \mathrm{Sl}(2; \mathbb{C})$ that is described by (4.2)–(4.4) with $\mathfrak{g}|_{s=0}$ being the identity section. Moreover, this family when reparametrized by the interval $[0,1]$ by writing $s = -\ln(1-r)$ extends smoothly to the $r = 1$ boundary point as a map from $[0,1]$ into the space of smooth sections of $P \times_{\mathrm{Ad}(SU(2))} \mathrm{Sl}(2; \mathbb{C})$. The $r = 1$ limit is denoted both by $\mathfrak{g}|_{s=\infty}$ and $\mathfrak{g}_\infty$. The small t, large t and small $\frac{t}{|z|}$ asymptotics of any $s \in [0, \infty]$ version of $\mathfrak{g}|_s$ is described by Proposition 5.1.

Given $\{\mathfrak{g}|_s\}_{s \in [0, \infty]}$, there is the corresponding data set family $\{(A|_s, \varphi|_s, \mathfrak{a}_3|_s)\}_{s \in [0, \infty]}$ that is defined from $\{\mathfrak{g}|_s\}_{s \in [0, \infty]}$ using (4.1). The $s = \infty$ member of this family is also denoted by $(A_\infty, \varphi_\infty, \mathfrak{a}_{\infty 3})$. This end member $(A_\infty, \varphi_\infty, \mathfrak{a}_{\infty 3})$ necessarily obeys all four bullets of (1.3) because its version of $\mathfrak{X}$ is identically zero. With regards to these data set: If the function $\mu_{(\cdot)}$ from (5.1) is bounded, then the respective small t, large t and small $\frac{t}{|z|}$ asymptotics of each member of the family $\{(A|_s, \varphi|_s, \mathfrak{a}_3|_s)\}_{s \in [0, \infty]}$ is described by Proposition 5.2 which says in effect that these asymptotics are the same as those for $(A_\diamond, \varphi_\diamond, \mathfrak{a}_{\diamond 3})$. In particular, this asymptotic convergence in Proposition 5.2 implies that each $s \in [0, \infty]$ version of $(A|_s, \varphi|_s, \mathfrak{a}_3|_s)$ is described by (1.10) if $(A_\diamond, \varphi_\diamond, \mathfrak{a}_{\diamond 3})$ is described by (1.10). This implies in turn that each $s \in [0, \infty]$ version of $(A|_s, \varphi|_s, \mathfrak{a}_3|_s)$ is described by the first and second bullets of (1.9).

With an eye towards the third bullet in (1.9), this section is concerned with the $A = A_\infty$ version of the function on $(0, \infty)$ whose value at any given R is the integral

$$\int_{(0,\infty)\times\{z\in\mathbb{R}^2:|z|>R\}} (|B_{A3}|^2 + |E_{A1}|^2 + |E_{A2}|^2) \ .$$

(6.1)

The following proposition is the central observation in this section about this function.

**Proposition 6.1**: *Suppose that* $(A_\diamond, \varphi_\diamond, \mathfrak{a}_{\diamond 3})$ *obeys the bullets in (1.3) and that the* $(0,\infty)\times\mathbb{R}^2$ *integral of* $(1+t^2)|\mathfrak{X}_\diamond|^2$ *is finite. Assume also that* $(A_\diamond, \varphi_\diamond, \mathfrak{a}_{\diamond 3})$ *is described by (1.10) and that the corresponding function* $\mu_{(\cdot)}$ *is bounded. Finally, make the assumptions in the two bullets that follow.*
- *Supposing that* $R > 1$, *then*

$$\int_{(0,\infty)\times\{z\in\mathbb{R}^2:|z|>R\}} (|B_{A_\diamond 3}|^2 + |E_{A_\diamond 1}|^2 + |E_{A_\diamond 2}|^2) \le \tfrac{M}{R} \quad \text{and} \quad \int_{(0,\infty)\times\{z\in\mathbb{R}^2:|z|>R\}} |\mathfrak{X}_\diamond|^2 < \tfrac{M}{R}$$

  *with* M *being independent of* R.
- *There exists* $t_0 \in (0, 1]$ *such that if* $t \in (0, t_0]$, *then the following hold:*
  a) *The* $\{t\}\times\mathbb{R}^2$ *integral of* $|\mathfrak{X}_\diamond|^2$ *is finite with a* t *independent upper bound.*
  b) *The integral of* $|\mathfrak{X}_\diamond|^2$ *over the* $|z| > R$ *part of* $\{t\}\times\mathbb{R}^2$ *is bounded by* $\tfrac{M}{R}$ *when* $R > 1$
     *with* M *being independent of* t *and* R.

*Granted these assumptions, construct the 1-parameter family of data set* $\{\mathfrak{g}|_s\}_{s\in[0,\infty]}$ *as done in Propositions 5.1, and then define the corresponding 1-parameter family of data sets* $\{(A|_s, \varphi|_s, \mathfrak{a}_3|_s)\}_{s\in[0,\infty]}$ *via (4.1). There exists* $M' > 0$ *such that each* $R > 1$ *version of the integral in (6.1) is bounded by* $\tfrac{M'}{R}$ *when* $(A, \mathfrak{a})$ *comes from any* $s \in [0,\infty]$ *version of the data set* $(A|_s, \varphi|_s, \mathfrak{a}_3|_s)$.

The proof of this proposition occupies the remaining subsections of Section 6.

By way of an example and application, consider the data set from Section 3 when defined from a p-tuple of complex numbers $(a_1, \ldots, a_p)$ with $a_p \ne 0$ and a positive number $\delta$ chosen to be less than $\kappa_\diamond$ (which is the version of $\kappa$ from Lemma 3.3). Lemmas 3.4 and 3.5 make assertions to the effect that any such Section 3 version of $(A_\diamond, \varphi_\diamond, \mathfrak{a}_{\diamond 3})$ has a corresponding version of M for use as input to the proposition.

### a) Formula for the s-derivative of $B_{A3}$, $E_{A1}$ and $E_{A2}$

Supposing that the s-derivatives of the mapping $s \to (A|_s, \varphi|_s, \mathfrak{a}_3|_s)$ is described by (4.6), then the corresponding s-derivatives of $B_{A3}$, $E_{A1}$ and $E_{A2}$ at any given value of s where (4.6) holds are given by the rule

- $\tfrac{\partial}{\partial s} B_{A3} = -(\nabla_{A1}\nabla_{A1} + \nabla_{A2}\nabla_{A2})u$ .
- $\tfrac{\partial}{\partial s} (E_{A1} + iE_{A2}) = -i[u, E_{A1} + iE_{A2}] - i(\nabla_{At} + i[\mathfrak{a}_3, \cdot])(\nabla_{A1}u + i\nabla_{A2}u)$

(6.2)

The following inequality is a consequence:

$$\left|\tfrac{d}{ds}(|B_{A3}|^2 + |E_{A1}|^2 + |E_{A2}|^2)\right| \leq c_0(|\nabla_A \nabla_A u|^2 + |\mathfrak{a}_3|^2 |\nabla_A u|^2) \;.$$

(6.3)

With the preceding understood, suppose now that ô is a smooth, non-negative function on $(0,\infty) \times \mathbb{R}^2$ with compact support. Then (6.3) leads to the bound:

$$\left|\tfrac{d}{ds}\left(\int_{(0,\infty)\times\mathbb{R}^2} \hat{o}^2(|B_{A3}|^2 + |E_{A1}|^2 + |E_{A2}|^2)\right)^{1/2}\right| \leq c_0 \left(\int_{(0,\infty)\times\mathbb{R}^2} \hat{o}^2 |\nabla_A \nabla_A u|^2\right)^{1/2}$$
$$+ c_0 \left(\int_{(0,\infty)\times\mathbb{R}^2} \hat{o}^2 |\mathfrak{a}_3|^2 |\nabla_A u|^2\right)^{1/2}$$

(6.4)

The plan is prove Proposition 6.1 via (6.4) using a sequence of ô's that equal 1 on the $|z| > R$ part of a corresponding sequence of nested, open sets in $(0, \infty) \times \mathbb{R}^2$ with compact closure that exhaust $(0, \infty) \times \mathbb{R}^2$. Bounds for the corresponding sequence of integrals on the right hand side of (6.4) are derived that imply the proposition.

b) **Preliminary bounds for the integrals of $\hat{o}^2|\nabla_A\nabla_A u|^2$ and $\hat{o}^2|\mathfrak{a}_3|^2|\nabla_A u|^2$**

The best that can be said at this point about the $\hat{o}^2|\mathfrak{a}_3|^2|\nabla_A u|^2$ integral in (6.4) is that

$$\int_{(0,\infty)\times\mathbb{R}^2} \hat{o}^2 |\mathfrak{a}_3|^2 |\nabla_A u|^2 \leq c_\diamond \int_{(0,\infty)\times\mathbb{R}^2} \hat{o}^2 \tfrac{1}{t^2} |\nabla_A u|^2$$

(6.5)

with $c_\diamond$ denoting here and subsequently a number that is greater than 1 which depends only on the initial data set $(A_\diamond, \varphi_\diamond, \mathfrak{a}_{\diamond 3})$. It's value can be assumed to increase between successive appearances. This bound follows from the bound in Proposition 5.2 and the fact that $(A_\diamond, \varphi_\diamond, \mathfrak{a}_{\diamond 3})$ is described by (1.10). The $\tfrac{1}{t^2}\hat{o}^2|\nabla_A u|^2$ integral on the right hand side of (6.5) is the subject of the upcoming Lemma 6.3.

The analysis for the $\hat{o}^2|\nabla_A\nabla_A u|^2$ integral in (6.4) starts with the $\mathfrak{X} \equiv e^{-s}\mathfrak{X}_\diamond$ version of equation in (4.3) which leads immediately to the inequality below (it holds for any compactly supported function ô):

$$\int_{(0,\infty)\times\mathbb{R}^2} \hat{o}^2 |\nabla_A^\dagger \nabla_A u|^2 \leq c_0 \int_{(0,\infty)\times\mathbb{R}^2} \hat{o}^2 |\mathfrak{a}_3|^2 |\nabla_A u|^2 + c_0 e^{-2s} \int_{(0,\infty)\times\mathbb{R}^2} \hat{o}^2 |\mathfrak{X}_\diamond|^2 \;.$$

(6.6)

(The notation has $\nabla_A^\dagger$ denoting the formal $L^2$ adjoint of the A-covariant derivative $\nabla_A$.) For now, (6.5) is all that can be said about the $\hat{o}^2|\mathfrak{a}_3|^2|\nabla_A u|^2$ integral on the right hand side here. Meanwhile, the assumptions in the proposition will be brought to bear later for the

ô²|𝔛₀|² integral on the right hand side of (6.6). To proceed with the task at hand, integrate by parts on the left hand side of (6.6) (multiple times) and commute the covariant derivatives (multiple times) to change $\nabla_A^\dagger \nabla_A$ into $\nabla_A \nabla_A$. Doing that leads to this:

$$\int_{(0,\infty)\times\mathbb{R}^2} ô^2|\nabla_A\nabla_A u|^2 \leq c_0 \int_{(0,\infty)\times\mathbb{R}^2} ô^2|\nabla_A^\dagger\nabla_A u|^2 + c_0 \int_{(0,\infty)\times\mathbb{R}^2} ô^2|F_A||\nabla_A u|^2 + c_0 \int_{(0,\infty)\times\mathbb{R}^2} ô^2|F_A|^2|u|^2$$
$$+ c_0 \int_{(0,\infty)\times\mathbb{R}^2} (ô|\nabla\nabla ô| + |\nabla ô|^2)|\nabla_A u|^2 .$$

(6.7)

(The notation here has $F_A$ denoting the curvature 2-form of the connection A; thus $|F_A|$ denotes the square root of $|B_{A3}|^2 + |E_{A1}|^2 + |E_{A2}|^2$.)

There are two terms with $|F_A|$ in (6.7). A suitable bound for the right most one is obtained using (5.15) and (5.16):

$$\int_{(0,\infty)\times\mathbb{R}^2} ô^2|F_A|^2|u|^2 \leq c_◊ e^{-2s} \int_{(0,\infty)\times\mathbb{R}^2} ô^2|F_A|^2 .$$

(6.8)

The upcoming Lemma 6.2 leads directly to the bound below for the $ô^2|F_A||\nabla_A u|^2$ integral on the right hand side of (6.7).

$$\int_{(0,\infty)\times\mathbb{R}^2} ô^2|F_A||\nabla_A u|^2 \leq c_◊ e^{-2s} \int_{(0,\infty)\times\mathbb{R}^2} ô^2|F_A|^2 + c_◊ \int_{(0,\infty)\times\mathbb{R}^2} ô^2 \frac{1}{t^2}|\nabla_A u|^2 .$$

(6.9)

Here is the promised lemma.

**Lemma 6.2**: *Granted the assumptions in Proposition 6.1, there exists $\kappa > 1$ which is such that $|\nabla_A u|_s| \leq \frac{\kappa}{t} e^{-s}$ at any given $(t,z) \in (0,\infty)\times\mathbb{R}^2$ and any given $s \in [0,\infty]$.*

*Proof of Lemma 6.2*: This follows from (5.46) given the pointwise $c_◊ \frac{1}{t}$ bound in Proposition 5.2 for $A|_s$ and the sup-norm bound $c_◊ e^{-s}$ for $|u|_s|$.

Introduce by way of notation $\mathfrak{F}$ to denote the positive number whose square is:

$$\mathfrak{F}^2 \equiv \int_{(0,\infty)\times\mathbb{R}^2} ô^2(|B_{A3}|^2 + |E_{A1}|^2 + |E_{A2}|^2)$$

(6.10)

To summarize where things stand with regard to the s-derivative of $\mathfrak{F}$ from (6.4):

$$|\tfrac{d}{ds}\mathfrak{F}| - c_\diamond e^{-s}\mathfrak{F} \le c_\diamond \Big(\int_{(0,\infty)\times\mathbb{R}^2} \hat{o}^2 \tfrac{1}{t^2}|\nabla_A u|^2\Big)^{1/2} + c_\diamond e^{-s}\Big(\int_{(0,\infty)\times\mathbb{R}^2} \hat{o}^2|\mathfrak{X}_\diamond|^2\Big)^{1/2} + c_\diamond \mathfrak{Y}_{\hat{o}}$$

(6.11)

where $\mathfrak{Y}_{\hat{o}}$ denotes square root of the integral in (6.7) with the derivatives on $\hat{o}$.

To proceed from here requires a lemma concerning integrals of $\tfrac{1}{t^2}|\nabla_A u|^2$.

**Lemma 6.3**: *Granted the assumptions in Proposition 6.1, there exists $\kappa > 100$ and $\kappa_* > \kappa$ such that for all $s \in [0,\infty)$ and for all $t_* \in (0, 1/\kappa]$, the corresponding $u|_s$ obeys*

- $\int_{(0,t_*]\times\mathbb{R}^2} \tfrac{1}{t^2}|\nabla_A u|^2 \le \kappa_* t_*^{1-1/\kappa} e^{-2s}$ ,

- $\int_{(0,\infty)\times\mathbb{R}^2} \tfrac{1}{t^2}|\nabla_A u|^2 \le \kappa_* e^{-2s}$ ;

*and if $R > 2$, then*

- $\int_{(0,t_*]\times\{z\in\mathbb{R}^2:|z|>\tfrac{1}{2}R\}} \tfrac{1}{t^2}|\nabla_A u|^2 \le \kappa_* t_*^{1-1/\kappa} \tfrac{1}{R} e^{-2s}$ ,

- $\int_{(0,\infty)\times\{z\in\mathbb{R}^2:|z|>\tfrac{1}{2}R\}} \tfrac{1}{t^2}|\nabla_A u|^2 \le \kappa_* \tfrac{1}{R} e^{-2s}$ .

This lemma is proved in the Section 4d from another lemma that describes the small t behavior of the $\{t\}\times\mathbb{R}^2$ integrals of $|u|^2$. The preceding lemma is used directly in the next subsection to finish the proof of Proposition 6.1.

### c) The proof of Proposition 6.1

With regards to the choice of $\hat{o}$: Fix for the moment $\delta \in (0,1]$ and $T \in (1,\infty)$ and $r > R$, and then set

$$\hat{o}(t,z) \equiv \chi(2(1-\tfrac{t}{\delta}))\chi(\tfrac{t}{T}-1)\chi(\tfrac{|z|}{r}-1)\hat{i}_R(t,z).$$

(6.12)

where $\hat{i}_R(\cdot) \equiv 1$ in the case when $R \le 1$, and where $\hat{i}_R(t,z) \equiv \chi(2(1-\tfrac{|z|}{R}))$ when $R > 1$. To be sure: This version of $\hat{o}$ is equal to 1 where both $t \in (\delta, T)$ and $|z| \in [R, r]$ and it is equal to zero where any of the following conditions hold: Either $t < \tfrac{1}{2}\delta$ or $t > 2T$ or $|z| > 2r$. In the case when $R > 1$, it is also equal to zero where $|z| < \tfrac{1}{2}R$.

When $\hat{o}$ is given by (6.12), then the integral in (6.7) with derivatives of $\hat{o}$ is bounded by $c_0$ times a sum of three integrals when $R \le 1$ and by the sum of the same three and a fourth when R is greater than 1. The three that appear in both cases are

$$\tfrac{1}{\delta^2} \int_{(\tfrac{1}{2}\delta,\, \delta)\times \mathbb{R}^2} |\nabla_A u|^2 \quad \text{and} \quad \tfrac{1}{T^2} \int_{(T,\, 2T)\times \mathbb{R}^2} |\nabla_A u|^2 \quad \text{and} \quad \tfrac{1}{r^2} \int_{(0,\, \infty)\times \{z\in \mathbb{R}^2:\, r<|z|<2r\}} |\nabla_A u|^2 \;.$$

(6.13)

The fourth integral (which requires R > 1) is this:

$$\tfrac{1}{R^2} \int_{(0,\, \infty)\times \{z\in \mathbb{R}^2:\, \tfrac{1}{2}R<|z|<R\}} |\nabla_A u|^2 \;.$$

(6.14)

With regards to the first three: The right most integral is bounded by $c_\diamond e^{-2s} \tfrac{1}{r^2}$ (invoke the $v = u$, $\mathfrak{P}_0 \equiv e^{-s} \mathfrak{X}_\diamond$ version of (4.23)) so it has limit 0 as $r \to \infty$ as does the s-integral of its square root over the domain $[0, \infty)$. The middle integral in (6.13), the one with T, is bounded by $c_\diamond e^{-2s} \tfrac{1}{T^2}$ (also courtesy of the $v = u$, $\mathfrak{P}_0 \equiv e^{-s} \mathfrak{X}_\diamond$ version of (4.23)) so it has limit zero as $T \to \infty$ as does the s-integral of its square root over the domain $[0, \infty)$. The left most integral in (6.13) is bounded by $c_\diamond \delta^{1/2} e^{-2s}$ courtesy of the first bullet in Lemma 6.3, and so it has limit zero as $\delta \to 0$ and so does the s-integral of its square root over the domain $[0, \infty)$.

Meanwhile: The integral in (6.14) is bounded by $c_\diamond \tfrac{1}{R^2} e^{-2s}$ (use that same version of (4.23), so the s-integral of its square root over the domain $[0, \infty)$ is bounded by $c_\diamond \tfrac{1}{R}$.

Given the preceding, and given the proposition's assumption to the effect that there is a $c_\diamond \tfrac{1}{R}$ bound for the integral of $|\mathfrak{X}_\diamond|^2$ over the $|z| > R$ part of $(0, \infty) \times \mathbb{R}^2$ when $R \geq 1$, then the first two bullets in Lemma 6.3 when $R \leq 2$, or the third and fourth bullets in Lemma 6.3 when $R > 2$ lead directly to the following bound:

$$|\tfrac{d}{ds} \mathfrak{F}| - c_\diamond e^{-s} \mathfrak{F} \leq c_\diamond e^{-s} (\tfrac{1}{r} + \tfrac{1}{T} + \delta^{1/2} + \min(1, \tfrac{1}{\sqrt{R}}))$$

(6.15)

To elaborate briefly: The first or third bullet in Lemma 6.3 are used to bound the $t < 1$ part of the $\tfrac{1}{t^2} |\nabla_A u|^2$ integral in (6.11) whereas the second and fourth are used to bound the $t \geq 1$ part of that integral.)

The inequality in (6.15) can be integrated from $s = 0$ to any given value of s from the set $[0, \infty]$, and doing so results in the bound

$$\mathfrak{F}|_s \leq c_\diamond (\mathfrak{F}_\diamond + \tfrac{1}{r} + \tfrac{1}{T} + \delta^{1/2} + \min(1, \tfrac{1}{\sqrt{R}}))$$

(6.16)

with $\mathfrak{F}_\diamond$ denoting the $(A_\diamond, \varphi_\diamond, \mathfrak{a}_{\diamond 3})$ version of $\mathfrak{F}$ which is $\mathfrak{F}|_{s=0}$.

The inequality in (6.16) implies first that the $\delta \to 0$ and $r, T \to \infty$ limits of $\mathfrak{F}|_s$ exist for any $s \in [0, \infty]$. Then, granted that, it implies (given the proposition's assumptions about $\mathfrak{F}_\diamond$) that each $s \in [0, \infty]$ version of (6.1)'s integral is bounded by $c_\diamond \tfrac{1}{R}$.

## d) Integrals of $|u|^2$ and the proof of Lemma 6.3

The crucial input to the proof of Lemma 6.3 is the following lemma about integrals of $|u|^2$ over certain domains in $(0,\infty)\times\mathbb{R}^2$.

**Lemma 6.4**: *Granted the assumptions in Proposition 6.1, there exists $\kappa > 1$ and $\kappa_* > \kappa$ such that for all $s \in [0,\infty)$, the corresponding $u|_s$ obeys*

- $\displaystyle\int_{\{t_*\}\times\mathbb{R}^2} |u|^2 \leq \kappa_* t_*^{4-1/\kappa} e^{-2s}$ *when* $t_* \in (0, \tfrac{1}{\kappa}]$ ;

- $\displaystyle\int_{(0,\infty)\times\mathbb{R}^2} \tfrac{1}{t^4}|u|^2 \leq \kappa_* e^{-2s}$

*and if* $R > 2$, *then*

- $\displaystyle\int_{(0,\infty)\times\{z\in\mathbb{R}^2:|z|\geq\frac{1}{2}R\}} \tfrac{1}{t^4}|u|^2 \leq \kappa_* e^{-2s}\tfrac{1}{R}$ .

- $\displaystyle\int_{\{t_*\}\times\{z\in\mathbb{R}^2:|z|>\frac{1}{2}R\}} |u|^2 \leq \kappa_* t_*^{4-1/\kappa} e^{-2s}\tfrac{1}{R}$ *when* $t_* \in (0, \tfrac{1}{\kappa}]$ .

This lemma is proved momentarily. Coming first is the proof of Lemma 6.3.

*Proof of Lemma 6.3*: The R-independent assertions are proved first. To start, note that the lemma's second bullet follows from the first bullet (which bounds the $t < t_*$ part of the integral) and the $v = u$, $\mathfrak{P}_0 = e^{-s}\mathfrak{X}_\diamond$ version of (4.23) (which gives a $c_\diamond$ bound for the $t > t_*$ part). To see about the first bullet: Because $u$ obeys (4.3), it obeys the $v = u$, $\mathfrak{P}_0 = e^{-s}\mathfrak{X}_\diamond$ version of (4.25). With that identity in mind, fix some small but positive number $\delta$ and multiply both sides of the $v = u$, $\mathfrak{P}_0 = e^{-s}\mathfrak{X}_\diamond$ version of (4.25) by $\tfrac{1}{t^2}\chi(2(1-\tfrac{t}{\delta}))\chi(\tfrac{t}{t_*}-1)$; then integrate the result over $(0,\infty)\times\mathbb{R}^2$. Two instances of integration by parts on the left hand side of that integral identity leads to this:

$$\int_{[\delta,t_*]\times\mathbb{R}^2} \tfrac{1}{t^2}|\nabla_A u|^2 \leq c_0 \tfrac{1}{\delta^4}\int_{[\frac{1}{2}\delta,\delta]\times\mathbb{R}^2} |u|^2 + c_0\tfrac{1}{t_*^4}\int_{[t_*,2t_*]\times\mathbb{R}^2} |u|^2 + c_0 \int_{(0,2t_*]\times\mathbb{R}^2} \tfrac{1}{t^4}|u|^2$$
$$+ c_0 e^{-2s}\left(\int_{(0,2t_*]\times\mathbb{R}^2} \tfrac{1}{t^4}|u|^2\right)^{1/2}\left(\int_{(0,2t_*]\times\mathbb{R}^2} |\mathfrak{X}_\diamond|^2\right)^{1/2}.$$

(6.17)

Now invoke the first bullet in Lemma 6.4 and the assumption about the $\{t\}\times\mathbb{R}^2$ integrals of $|\mathfrak{X}_\diamond|^2$ being bounded to see that

$$\int_{[\delta,t_*]\times\mathbb{R}^2} \tfrac{1}{t^2}|\nabla_A u|^2 \le c_\lozenge e^{-2s}(\delta^{1-1/c_0} + t_*^{1-1/c_0}) \ .$$

(6.18)

if $t_* < c_\lozenge^{-1}$. The latter bound implies first that the $\delta \to 0$ limit of the left hand side in (6.18) is finite; and it then implies that there is a bound for that limit which has the form that is asserted by the first bullet Lemma 6.3.

Turn next to Lemma 6.3's third bullet inequality. No generality is lost here by assuming that $R > 4$ because the first of the lemma's inequalities can be used if $R < 4$. The argument for the third bullet is much like the argument for the first bullet with the differences being as follows: Multiply both sides of the $v = u$, $\mathfrak{P}_0 = e^{-s}\mathfrak{X}_\lozenge$ version (4.25) by the R-dependent product $\tfrac{1}{t^2}\chi(2(1-\tfrac{t}{\delta}))\chi(\tfrac{t}{t_*}-1)\chi(2(1-\tfrac{|z|}{R}))$. Then integrate by parts as before. The result of doing so is the analog of (6.17) which has the domain of all integrals on the right hand side being $(0, 2t_*]\times\{z\in\mathbb{R}^2: |z| > \tfrac{1}{2}R\}$ and which has a new term on the right hand side:

$$c_0 \tfrac{1}{R^2} \int_{(0,2t_*]\times\{z\in\mathbb{R}^2:\tfrac{1}{2}R<|z|<R\}} \tfrac{1}{t^2}|u|^2 \ .$$

(6.19)

Granted all of these terms, then the third bullet's inequality follows from the third and fourth bullets in Lemma 6.4 and from the assumption to the effect that the integral of $|\mathfrak{X}_\lozenge|^2$ on the $|z| > R$ part of $(0, c_0^{-1}]\times\mathbb{R}^2$ is at most an fixed multiple of $\tfrac{1}{R}$.

To derive the inequality in the lemma's fourth bullet, it is sufficient to consider the $t \ge t_*$ part of the integration domain because the lemma's third bullet deals with the remaining part. To deal with the $t \ge t_*$ part, multiply both sides of the $v = u$, $\mathfrak{P}_0 = e^{-s}\mathfrak{X}_\lozenge$ version (4.25) by the function $\tfrac{1}{t^2}\chi(2(1-\tfrac{t}{t_*}))\chi(2(1-\tfrac{|z|}{R}))$. Then integrate over $(0,\infty)\times\mathbb{R}^2$ and integrate by parts as before. One again obtains an analog of (6.17) which is this:

$$\int_{(t_*,\infty)\times\{z\in\mathbb{R}^2:|z|>R\}} \tfrac{1}{t^2}|\nabla_A u|^2 \le c_0 \tfrac{1}{t_*^4} \int_{[\tfrac{1}{2}t_*,t_*]\times\{z\in\mathbb{R}^2:|z|>\tfrac{1}{2}R\}} |u|^2 + c_0 \int_{[\tfrac{1}{2}t_*,\infty)\times\{z\in\mathbb{R}^2:|z|\ge\tfrac{1}{2}R\}} \tfrac{1}{t^4}|u|^2$$
$$+ c_0 \tfrac{1}{R^2} \int_{[\tfrac{1}{2}t_*,\infty)\times\{z\in\mathbb{R}^2:\tfrac{1}{2}R<|z|<R\}} \tfrac{1}{t^2}|u|^2 + c_0 e^{-2s} \int_{[\tfrac{1}{2}t_*,\infty)\times\{z\in\mathbb{R}^2:|z|\ge\tfrac{1}{2}R\}} |\mathfrak{X}_\lozenge|^2 \ .$$

(6.20)

The left most two terms that appears on the right hand side of (6.20) are no larger than $c_\lozenge e^{-2s}\tfrac{1}{R}$ when $R > 1$ courtesy of the fourth bullet in Lemma 6.4 for the left most term and the third bullet in Lemma 6.4 for the second to the leftmost term. Meanwhile, the third term from the left on the right hand side of (6.20) is bounced by $c_\lozenge e^{-2s}\tfrac{1}{R^2}$ courtesy of the $v=u$, $\mathfrak{P}_0 = e^{-s}\mathfrak{X}_\lozenge$ version of (6.23) and then Hardy's inequality from (4.15). (Because R is assumed to be greater than 1, this bound is sufficient for the purposes at hand.) The right most term on the right hand side of (6.20) is bounded by $c_0 e^{-2s} \frac{1}{R}$ by virtue of the assumptions of in Proposition 6.1.

***Proof of Lemma 6.4***: The proof has six parts. The first two address the first and second bullet's inequalities, the third addresses the inequality in the third bullet, and the remaining part addresses the one in the fourth bullet.

*Part 1*: The second bullet of the lemma follows directly from the first and the fact that $\frac{1}{t^4}|u|^2 \le \frac{1}{t_*^2}\frac{1}{t^2}|u|^2$ where $t \ge t_*$ and the fact that the $(0,\infty) \times \mathbb{R}^2$ integral of $\frac{1}{t^2}|u|^2$ is finite and bounded by $c_0 e^{-2s}$ (use the $v = u$, $\mathfrak{P}_\diamond = e^{-s}\mathfrak{X}_\diamond$ version of (4.23) with the $f = |u|$ version of (4.15).)

To prove the inequality in the first bullet, start with the $v = u$, $\mathfrak{P}_\diamond = e^{-s}\mathfrak{X}_\diamond$ version of (4.25) and integrate both sides over a given $t \in (0, t_0]$ slice $\{t\} \times \mathbb{R}^2$. That integral will be finite for all t. Doing that and integrating by parts twice to deal with the $\mathbb{R}^2$ derivatives leads to

$$-\frac{d^2}{dt^2}\left(\int_{\{t\}\times\mathbb{R}^2}|u|^2\right)^{1/2} + \left(\int_{\{t\}\times\mathbb{R}^2}|u|^2\right)^{-1/2}\int_{\{t\}\times\mathbb{R}^2}\|[\mathfrak{a},u]\|^2 \le c_0 e^{-s}\left(\int_{\{t_*\}\times\mathbb{R}^2}|\mathfrak{X}_\diamond|^2\right)^{1/2}.$$

(6.21)

The integration by parts used here can be justified by first using very large $r$ versions of the function $\chi(\frac{|z|}{r} - 1)$ to restrict integrals to the $|z| < 2r$ part of $\{t\} \times \mathbb{R}^2$ and then taking $r$ ever larger to see from the fact that $|u| \to 0$ as $\frac{t}{|z|} \to 0$ that the 'boundary' terms limit to zero as $r \to \infty$.

To continue: Remember that $(A_\diamond, \varphi_\diamond, \mathfrak{a}_{\diamond 3})$ obeys (1.10) and that it's version of the function $\mu_{(\cdot)}$ is bounded. This implies (via Proposition 5.2) the following: Given $\varepsilon > 0$, there exists $r_\varepsilon > 1$ with the following significance: The three components of any $s \in [0, \infty]$ version of $\mathfrak{a}|_s$ at any given point in the $|z| > r_\varepsilon t$ part of $(0, \infty) \times \mathbb{R}^2$ can be written as $-\frac{1}{2t}(\hat{\sigma}_3, \hat{\sigma}_1, \hat{\sigma}_2) + \mathfrak{e}$ with $\{\hat{\sigma}_a\}_{a=1,2,3}$ denoting an orthonormal frame for ad(P) at that point and with $\mathfrak{e}$ obeying $|\mathfrak{e}| \le \frac{\varepsilon}{t}$. Therefore, (6.21) leads in turn to

$$-\frac{d^2}{dt^2}\left(\int_{\{t\}\times\mathbb{R}^2}|u|^2\right)^{1/2} + \frac{2-c_0\varepsilon}{t^2}\left(\int_{\{t\}\times\{z\in\mathbb{R}^2:|z|>r_\varepsilon t\}}|u|^2\right)^{1/2} \le c_0 e^{-s}\left(\int_{\{t_*\}\times\mathbb{R}^2}|\mathfrak{X}_\diamond|^2\right)^{1/2}.$$

(6.22)

Because the norm of $|u|$ is bounded by $c_0\sqrt{t}\, e^{-s}$ where $t < 1$ (see (4.64), the integral of $|u|^2$ over the $|z| < r_\varepsilon t$ part of $\{t\} \times \mathbb{R}^2$ is at most $c_0 t^3 e^{-s}$. This can then be added to both sides of (6.22) to obtain:

$$-\tfrac{d^2}{dt^2}\Big(\int_{\{t\}\times\mathbb{R}^2}|u|^2\Big)^{1/2}+\tfrac{2-c_0\varepsilon}{t^2}\Big(\int_{\{t\}\times\mathbb{R}^2}|u|^2\Big)^{1/2}\leq c_\diamond(1+tr_\varepsilon^2)e^{-s}\;.$$

(6.23)

(The assumed a priori bound on the $\{t\}\times\mathbb{R}^2$ integral of $|\mathfrak{X}_\diamond|^2$ has been invoked here too.)

*Part 2*: Fix a positive number to be denoted by $c$ and let $\mathfrak{f}$ denote the function on $[0, t_0]$ that is defined by the rule

$$\mathfrak{f}(t)=\Big(\int_{\{t\}\times\mathbb{R}^2}|u|^2\Big)^{1/2}-ct^2(|\ln t|+t)\,e^{-s}\;.$$

(6.24)

The inequality in (6.23) says in effect that

$$-\tfrac{d^2}{dt^2}\mathfrak{f}+\tfrac{2-c_0\varepsilon}{t^2}\mathfrak{f}\leq 0$$

(6.25)

if $c>c_\diamond$. Assume henceforth that $c$ satifies this bound.

With (6.25) in mind, fix for the moment a positive number to be denoted by $\lambda$ so as to consider the kernel of the operator

$$-\tfrac{d^2}{dt^2}+\tfrac{\lambda}{t^2}$$

(6.26)

on the space of functions on $[0,\infty)$. In particular, the kernel is spanned by the elements

$$t^{(1+\varpi)/2}\quad\text{and}\quad t^{(1-\varpi)/2}$$

(6.27)

with $\varpi$ denoting here $(1+4\lambda)^{1/2}$. In the case at hand, $\lambda=2-c_0\varepsilon$ and so $\varpi=3-c_0\varepsilon$.

With (6.27) understood, and keeping in mind that the $\tfrac{1}{t^2}\mathfrak{f}$ has finite integral on the interval $(0,t_0]$ and that $\mathfrak{f}$ is bounded at $t=t_0$, the comparison principal can be invoked for (6.25) in the case at hand to see that $\mathfrak{f}\leq c_\diamond e^{-s}t^{(1+\varpi)/2}$ which is to say that $\mathfrak{f}\leq c_\diamond e^{-s}t^{2-c_0\varepsilon}$. The latter bound implies in turn that

$$\int_{\{t\}\times\mathbb{R}^2}|u|^2\leq c_\diamond(1+r_\varepsilon^2)^2\,e^{-2s}\,t^{4-c_0\varepsilon}\;.$$

(6.28)

Any $\varepsilon<c_\diamond^{-1}$ version of (6.28) gives the first inequality in Lemma 6.4.

*Part 3*: To prove the inequality in the third bullet: Fix some very small, positive number to be denoted by $\delta$ and multiply both sides of the $v = u$, $\mathfrak{P}_9 = e^{-s}\mathfrak{X}_\diamond$ version of (4.25) by the function $\frac{1}{t^2}\chi(2(1-\frac{t}{\delta}))^2\chi(2(1-\frac{|z|}{R}))^2$. (To be sure: This function is equal to $\frac{1}{t^2}$ where by $t > \delta$ and $|z| > R$, and it is equal to zero where either $t < \frac{1}{2}\delta$ or $|z| < \frac{1}{2}R$.) Having done that, integrate the result over the $(0,\infty) \times \mathbb{R}^2$ and then integrate by parts on the left side to remove derivatives from $|u|^2$. (There are no large t or large $|z|$ obstructions to doing that.) The result leads to an integral inequality for the function h given by the rule $h \equiv \chi(2(1-\frac{t}{\delta}))\chi(2(1-\frac{|z|}{R}))|u|$ that is given below in (6.29). By way of notation, (6.29) uses $\hat{u}$ to denote $\frac{u}{|u|}$.

$$\int_{(0,\infty)\times\mathbb{R}^2} \tfrac{1}{t^2}|\nabla h|^2 - 3\int_{(0,\infty)\times\mathbb{R}^2} \tfrac{1}{t^4}|h|^2 + \int_{(0,\infty)\times\mathbb{R}^2} \tfrac{1}{t^2}|[\mathfrak{a},\hat{u}]|^2 h^2 \leq c_0 \tfrac{1}{R^2}\int_{(0,\infty)\times\mathbb{R}^2} \tfrac{1}{t^2}|u|^2$$
$$+ c_0 \tfrac{1}{\delta^4}\int_{[\frac{1}{2}\delta,\delta]\times\mathbb{R}^2}|u|^2 + c_0 e^{-s}\Big(\int_{(0,\infty)\times\mathbb{R}^2}\tfrac{1}{t^4}|h|^2\Big)^{1/2}\Big(\int_{(0,\infty)\times\{z\in\mathbb{R}^2:\,|z|>\frac{1}{2}R\}}|\mathfrak{X}_\diamond|^2\Big)^{1/2}$$

(6.29)

Some points of note: First, the term on the right with the factor $\frac{1}{R^2}$ is no greater than $c_\diamond \frac{1}{R^2}e^{-2s}$, this being due to the $v = u$, $\mathfrak{P}_0 = e^{-s}\mathfrak{X}_\diamond$ version of (4.23) and Hardy's inequality in (4.15). Second, the term with the factor $\frac{1}{\delta^4}$ is bounded by $c_\diamond \delta^{1/c_0} e^{-2s}$, this being due to the first bullet of Lemma 6.4. Third: Having specified $\varepsilon \in (0,1]$, then the term with the integral of $\frac{1}{t^2}|[\mathfrak{a},\hat{u}]|^2 h^2$ is greater than $(2-\varepsilon)\frac{1}{t^4}h^2$ if $R > 2r_\varepsilon$ with $r_\varepsilon$ depending only on $\varepsilon$ and $(A_\diamond, \varphi_\diamond, \mathfrak{a}_{\diamond 3})$. (This is because of Proposition 5.2.) Therefore, having fixed $\varepsilon$, and supposing that $R > 2r_\varepsilon$, then (6.29) leads to the following:

$$\int_{(0,\infty)\times\mathbb{R}^2} \tfrac{1}{t^2}|\nabla h|^2 - (1+2\varepsilon)\int_{(0,\infty)\times\mathbb{R}^2} \tfrac{1}{t^4}|h|^2 \leq c_\diamond\big(\tfrac{1}{\varepsilon}\tfrac{1}{R} + \delta^{1/c_0}\big)e^{-2s}$$

(6.30)

after using the assumption about the integral of $|\mathfrak{X}_\diamond|^2$ on the $|z| > R$ part of $(0,\infty)\times\mathbb{R}^2$.

To proceed from here requires yet another version of Hardy's inequality which is this: Suppose that $f$ is a smooth function on $(0,\infty)$ that is zero near $t = 0$ and such that $|\frac{d}{dt}f|^2$ has finite integral on $(0,\infty)$. Then

$$\int_{(0,\infty)} \tfrac{1}{t^4}|f|^2 \leq \tfrac{4}{9}\int_{(0,\infty)} \tfrac{1}{t^2}|\tfrac{d}{dt}f|^2 \,.$$

(6.31)

(See the appendix for a proof.) Apply this inequality to each $z \in \mathbb{R}^2$ version of $h(\cdot, z)$ in (6.31) to see that

$$(1-2\varepsilon) \int_{(0,\infty)\times\mathbb{R}^2} \tfrac{1}{t^4}|h|^2 \le c_\diamond(\tfrac{1}{\varepsilon}\tfrac{1}{R}+ \delta^{1/c_0}) \, e^{-2s}.$$

(6.32)

Take $\varepsilon = \tfrac{1}{4}$ and then take $\delta \to 0$ to obtain assertion of Lemma 6.4's third bullet when R is larger than the $\varepsilon = \tfrac{1}{4}$ version of $r_\varepsilon$. The assertion of the lemma's third bullet follows from the lemma's second bullet when R is less that this same version of $r_\varepsilon$.

*Part 4*: Given the third bullet, the proof of the inequality in the fourth bullet is almost the same as that in the first but for the very beginning. In this case, the proof starts by multipling both sides of the $v = u$, $\mathfrak{P}_\mathfrak{g} = e^{-s}\mathfrak{X}_\diamond$ version of (4.25) by the function $\chi(2(1-\tfrac{|z|}{R}))^4$ which is equal to one where $|z| > R$ and equal to zero where $|z| < \tfrac{1}{2}R$. After this multiplication, integrate both sides over a given $\{t\}\times\mathbb{R}^2$ slice for $t \in (0, t_*]$. The result of doing that and integrating by parts leads to the following analog of (6.21):

$$-\tfrac{d^2}{dt^2}(\int_{\{t\}\times\mathbb{R}^2} \chi_R^4 |u|^2)^{1/2} + \tfrac{2-c_0\varepsilon}{t^2}(\int_{\{t\}\times\{z\in\mathbb{R}^2:|z|>r_\varepsilon t\}} \chi_R^4|u|^2)^{1/2} \le c_0 e^{-s}(\int_{\{t_*\}\times\mathbb{R}^2} \chi_R^4|\mathfrak{X}_\diamond|^2)^{1/2}$$
$$+ c_0 \tfrac{1}{R^2}(\int_{\{t\}\times\mathbb{R}^2} |u|^2)^{1/2}.$$

(6.33)

Note that the right most term on the right hand side is bounded by $c_\diamond \tfrac{1}{R^2} e^{-s}$ by virtue of the first bullet of the lemma. Meanwhile, the left most term on the right hand side of (6.33) is no greater than $c_\diamond \tfrac{1}{\sqrt{R}} e^{-s}$ because of the assumptions for Proposition 6.1.

Therefore, supposing that $R > 4r_\varepsilon$, the inequality in (6.33) implies that

$$\mathfrak{f}_R(t) \equiv (\int_{\{t\}\times\mathbb{R}^2} \chi_R^4|u|^2)^{1/2} - ct^2|\ln t| \, e^{-s}\tfrac{1}{\sqrt{R}}$$

(6.34)

obeys the inequality in (6.25) if $c > c_\diamond$. Granted such a choice for $c$, fix $s \in [\tfrac{1}{2}t_*, t_*]$ and, supposing that $t < \tfrac{1}{2}t_*$, use the maximum principle with (6.25) to see that

$$(\int_{\{t\}\times\mathbb{R}^2} \chi_R^4|u|^2)^{1/2} \le c_\diamond \, t^{2-c_0\varepsilon} \tfrac{1}{\sqrt{R}} e^{-s} + (\tfrac{t}{s})^{2-c_0\varepsilon}(\int_{\{s\}\times\mathbb{R}^2} \chi_R^4|u|^2)^{1/2}.$$

(6.35)

Averaging this over values of $s$ in the interval $[\tfrac{1}{2}t_*, t_*]$ and then invoking the third bullet's inequality leads directly to the fourth bullets inequality if R is greater than a $c_\diamond^{-1}$ version of $r_\varepsilon$. If R is less than that, then the fourth bullet's inequality follows from the first bullet of the lemma.

## 7. Proof of Theorem 1

The proof of Theorem 1 from Section 1 has two parts, the first is to verify that there exists at least one instanton solution for each non-negative integer $m$ and positive integer $p$. The second part is to verify the claim to the effect that there is a family of pairwise Aut(P)-inequivalent instanton solutions that is parametrized by $\mathbb{C}^{p-1} \times (\mathbb{C}-0)$. These respective parts of the proof are dealt with in Sections 7a and 7b. There is also a Section 7c concerning the possibility of a larger moduli space.

### a) Existence

Fix a p-tuple of complex numbers $(a_1, \ldots, a_p)$ with $a_p \neq 0$ and then fix a positive number $\delta$ which is small enough so that the pair $(A, \mathfrak{a})$ that is constructed in Section 3 is described by Lemmas 3.2-3.5. This pair is denoted by $(A_\diamond, \mathfrak{a}_\diamond)$ and, equivalently, $(A_\diamond, \varphi_\diamond, \mathfrak{a}_{\diamond 3})$. The construction is such that $(A_\diamond, \varphi_\diamond, \mathfrak{a}_{\diamond 3})$ obeys (1.3) and (1.10); and because it is described by Lemmas 3.2-3.5, it meets the requirements for Propositions 4.1, 5.1, 5.2 and 6.1. As a consequence, these propositions can be invoked to see that there is a 1-parameter family $s \to \mathfrak{g}|_s$ of smooth sections of $P \times_{\mathrm{Ad}(\mathrm{SU}(2))} \mathrm{Sl}(2; \mathbb{C})$ parameterized by $[0, \infty]$ that is described by (4.1)-(4.6). More to the point, the $s = \infty$ end member of corresponding family $(A|_s, \varphi|_s, \mathfrak{a}_3|_s)_{s \in [0, \infty]}$ from (4.1) obeys all of the bullets in (1.3) and it is described by (1.9) and (1.10). This end member is the sought after instanton solution to the equations in (1.3).

### b) The parameter space

This section explains why two different choices for the parameters $(a_1, \ldots, a_p)$ for use in Section 3 to construct an initial data set $(A_\diamond, \varphi_\diamond, \mathfrak{a}_{\diamond 3})$ lead via Propositions 4.1, 5.1,. 5.2 and 6.1 to corresponding instanton solutions to (1.3) that are not equivalent via the action of any automorphism of the bundle P. The proof is by contradiction: To this end, suppose for the sake of argument that parameter sets $\mathbf{a} \equiv (a_1, \ldots, a_p)$ and $\mathbf{a}' \equiv (a_1', \ldots, a_p')$ lead via Section 3's constructions to corresponding pairs $(A_\diamond, \mathfrak{a}_\diamond)$ and $(A_\diamond', \mathfrak{a}_\diamond')$ and then, via respective $s \to \infty$ limits (4.1)-(4.6), to corresponding instanton solutions to (1.3) that are Aut(P) equivalent. These corresponding instanton solution pairs are denoted by $(A, \mathfrak{a})$ and $(A', \mathfrak{a}')$. To say that they are Aut(P) equivalent means that there is a section (to be denoted by h) of the fiber bundle $P \times_{\mathrm{Ad}(\mathrm{SU}(2))} \mathrm{SU}(2)$ such that

$$A' = hAh^{-1} - (\nabla_A h)h^{-1} \quad \textit{and} \quad \mathfrak{a}' = h\mathfrak{a}h^{-1}.$$

(7.1)

The proof that this assumption leads to nonsense follows directly in two parts.

*Part 1*: There is an $(A_\diamond, \mathfrak{a}_\diamond)$ version of $\sigma$, $\varphi$ and $\lambda$ that obey (2.16) and (2.17) with $\lambda$ being $\sigma - 4q$ and with $q$ being a section of $\sigma$'s version of $\mathcal{L}^+$ which is smooth on the complement of the z=0 locus. This version is denoted by $\sigma_\diamond$, $\varphi_\diamond$ and $\lambda_\diamond$. These are described in Sections 3a-3c. Note in particular that $q_\diamond$ is smooth on the complement of the z = 0 locus where it has a pole, and thus so does $\lambda_\diamond$. The depiction of $q_\diamond$ is the $(A_\diamond, \varphi_\diamond, \mathfrak{a}_{\diamond 3})$ instance of the formula in (3.25).

For the present purposes, the salient features of $\lambda_\diamond$ are listed in the upcoming (7.3). The notation in (7.2) uses $\hbar$ to denote the meromorphic function

$$\hbar \equiv \tfrac{1}{4} \tfrac{1}{a_p} ( \tfrac{1}{z^{m+p}} + \mu_{p-1} \tfrac{1}{z^{m+p-1}} + \cdots + \mu_1 \tfrac{1}{z^{m+1}} ) .$$

(7.2)

that appears in (3.25). (It is constructed from the data set **a** using the polynomial $\mathfrak{P}$ that appears in (3.2).) What follows directly are the important features of $\lambda_\diamond$ (these are all straightforward conseqences of (3.25)'s depiction of $q_\diamond$):

- $|\lambda_\diamond| \to 1$ *as* $t \to 0$ *and uniformly so on sets where $|z|$ is bounded away from zero.*
- $\lambda_\diamond = -\hbar \varphi_\diamond + \mathfrak{e}_\diamond$ *near $z = 0$ with $\mathfrak{e}_\diamond$ extending continuously across the $z = 0$ locus.*

(7.3)

The pair $(A, \mathfrak{a})$ has a corresponding $\sigma$, $\varphi$ and $\lambda$ with $\varphi$ and $\lambda$ determined by the $s = \infty$ end member of the family $\{\mathfrak{g}|_s\}_{s \in [0, \infty]}$ according to the rule whereby

$$\varphi = \mathfrak{g}|_\infty \varphi_\diamond (\mathfrak{g}|_\infty)^{-1} \quad and \quad \lambda = \mathfrak{g}|_\infty \lambda_\diamond (\mathfrak{g}|_\infty)^{-1} .$$

(7.4)

Keep in mind that $\sigma$ is defined from $\varphi$ by writing $[\varphi, \varphi^*]$ as $2i|\varphi|^2\sigma$ and that the corresponding $q$ which is defined from $\lambda$ by setting $q \equiv -2i[\sigma, \lambda]$. Also keep in mind from Propositions 5.1 and 5.2 that $\mathfrak{g}|_\infty$ limits to the identity matrix as $t \to 0$ and as $\tfrac{t}{|z|} \to 0$ and as $t \to \infty$. A consequence of this is that the top bullet of (7.3) hold with $\lambda_\diamond$ replaced by $\lambda$ on the left hand side; and that the second bullet in (7.3) is replaced by the following:

$\lambda = -\hbar \varphi + \mathfrak{e}$ *near $z = 0$ with $\mathfrak{e}$ extending continuously across the $z = 0$ locus.*

(7.5)

Meanwhile, there are $(A_\diamond', \mathfrak{a}_\diamond')$ and $(A', \mathfrak{a}')$ versions of $(\sigma, \varphi, \lambda)$, these denoted respectively by $(\sigma_\diamond', \varphi_\diamond', \lambda_\diamond')$ and $(\sigma', \varphi', \lambda')$. And, there is a primed version of (7.2)–(7.5) with $\hbar$ replaced by a different meromorphic function which is denoted by $\hbar'$; it is defined from the data set $\mathbf{a}'$ by using $\mathbf{a}'$ to define a primed version of the polynomial that is depicted in (3.2).

*Part 2*: The right hand identity in (7.1) says in part that $\varphi' = h\varphi h^{-1}$ and this implies that $\sigma' = h\sigma h^{-1}$. However, this is not necessarily the case of $\lambda'$. What is true is that $\lambda'$ must have the form

$$\lambda' = h\lambda h^{-1} + f h\varphi h^{-1}$$

(7.6)

with $f$ being a t-independent, meromorphic function of z. (This is the only possibility that is consistent with the primed and unprimed versions of (2.16) and (2.17).) The next paragraphs use (7.3) and (7.5) and their primed analogs to constrain $f$.

The first constraint is that $f$ must be analytic on the complement of the $z = 0$ locus. This because both $\lambda$ and $\lambda'$ are are analytic there. The second constraint comes from the primed and unprimed versions of (7.5); it is that $f$ must have the form

$$f = \hbar' - \hbar + t$$

(7.7)

with $t$ being meromorphic with pole order along the $z = 0$ locus no greater than $m$. The third constraint comes from the primed and unprimed versions of the top bullet in (7.3), which is that $f \equiv 0$. This is because $t|\varphi| \to \frac{1}{\sqrt{2}}$ uniformly as $t \to 0$ on sets where $|z|$ has a positive lower bound. It is perhaps needless to say that this last constraint is not compatible with (7.7) unless $\hbar = \hbar'$. This incompatibility is the desired nonsense because $\hbar$ can't equal $\hbar'$ if **a** is not equal to **a**´.

### c) Other parameters?

An initial data set $(A_\diamond, \varphi_\diamond, \mathfrak{a}_{\diamond 3})$ from Section 3 required the choice of the complex parameters $(a_1, \ldots, a_p)$ and then a choice of a small real number, $\delta$. This section explains why two different choices of $\delta$ with the same $(a_1, \ldots, a_p)$ lead to Aut(P)-equivalent instanton solutions to (1.3). Two explanations are given, these being Parts 1 and 2 of what follows.

*Part 1*: Suppose that $(A, \mathfrak{a})$ and $(A', \mathfrak{a}')$ are two instanton solutions to (1.3), both described by (1.9) and (1.10). Then, by virtue of both obeying the first three bullets in (1.3), the connection $A'$ can be written as $A + \mathbb{A}$ and $\mathfrak{a}'$ as $\mathfrak{a} + \hat{\mathfrak{a}}$ with $\mathbb{A}$ and $\hat{\mathfrak{a}}$ given in terms of a section (denoted by $\mathfrak{g}$) of $P \times_{\text{Ad}(SU(2))} Sl(2; \mathbb{C})$:

- $\mathbb{A}_t - i\hat{\mathfrak{a}}_3 = -(\nabla_{At}\mathfrak{g})\mathfrak{g}^{-1} + i[\mathfrak{a}_3, \mathfrak{g}]\mathfrak{g}^{-1}$
- $\mathbb{A}_1 + i\mathbb{A}_2 = -((\nabla_{A1} + i\nabla_{A2})\mathfrak{g})\mathfrak{g}^{-1}$.
- $\hat{\mathfrak{a}}_1 - i\hat{\mathfrak{a}}_1 = -[\varphi, \mathfrak{g}]\mathfrak{g}^{-1}$

(7.8)

In the case at hand, where both $(A, \mathfrak{a})$ and $(A', \mathfrak{a}')$ come from a Section 3 data set, then the section $\mathfrak{g}$ can be changed by an automorphism of P so that it limits to the identity section $\mathbb{I}$ as $t \to 0$, as $t \to \infty$ and as $\frac{t}{|z|} \to 0$. Assuming only that (1.10) is obeyed by both $(A, \mathfrak{a})$ and $(A, \mathfrak{a}')$, then $\mathfrak{g}$ can be changed by an automorphism of P so that the $t \to 0$, $t \to \infty$ and $\frac{t}{|z|} \to 0$ limits of $t|A| + t|\mathfrak{a}|$ are zero. Granted this, then the third bullet in (7.8) implies directly that the $t \to 0$ and $t \to \infty$ limits of $\mathfrak{g}$ have the form

$$\mathfrak{g} \sim \mathbb{I} + h \varphi \tag{7.9}$$

with h being asymptotic to a t-independent, holomorphic function of z. Meanwhile $\mathfrak{g}$ as, $\frac{t}{|z|} \to 0$ should be of the form in (7.9) up to terms that grow slower as a function of z than a holomorphic function on $\mathbb{R}^2$. As a parenthetical remark in this regard: The three different limiting holomorphic functioins need not all be the same. In fact, if they are, then $\mathfrak{g}$ can be modified by subtracting $h \varphi$ so that (7.8) still holds and so that the new version of $\mathfrak{g}$ is asymptotic to $\mathbb{I}$ in all three limiting cases.

The crucial observation now is that the function $\tau \equiv \operatorname{trace}(\mathfrak{g}\mathfrak{g}^\dagger) - 2$ obeys the miraculous differential inequality

$$\Delta \tau > 0, \tag{7.10}$$

an inequality that was discovered by Simon Donaldson in the Hermitian Yang-Mills context [D] and used to good effect by He and Mazzeo in their Kapustin-Witten context in [HM1-3]. (It is a strenuous exercise to derive (7.10) from the fourth bullet of (1.3)!)

With (7.10) in mind: If the $t \to 0$, $t \to \infty$ and $\frac{t}{|z|} \to 0$ limits of $\mathfrak{g}$ are the identity section, then the corresponding limits of $\tau$ are zero; and in this case, $\tau$ is identically zero by virtue of (7.10) and the maximum principle. Indeed, since $\mathfrak{g}$ has determinant one, the trace of $\mathfrak{g}\mathfrak{g}^\dagger$ can't be less than 2 and it can equal 2 if and only if $\mathfrak{g}\mathfrak{g}^\dagger = \mathbb{I}$. If that is the case (which is saying that $\tau \equiv 0$), then $\mathfrak{g}$ is an automorphism of P (it is unitary) in which case (7.8) says that $(A', \mathfrak{a}')$ is Aut(P) equivalent to $(A, \mathfrak{a})$.

*Part 2*: What follows is a sketch leaving the details to the reader (the details use no significant new analysis; they use only slight extensions some from Sections 4-6.)

Letting $\delta$ vary for fixed $(a_1, \ldots, a_p)$ will give a $\delta$-parametrized family of limit data sets of the form $(A|_{s=\infty}, \varphi|_{s=\infty}, \mathfrak{a}_3|_{s=\infty})$ using Propositions 4.1, 5.1, 5.2 and 6.1. It follows from the manner of convergence as $s \to \infty$ that is described in Propositions 5.1 and 5.2 that this limit family will vary smoothly with respect to the parameter $\delta$. With that understood, suppose that $\delta_0$ and $\delta_1$ are two small choices for $\delta$. Because each end member of this family obeys the top three bullets in (1.3), there will be a smooth map $\delta \to \mathfrak{g}_{(\delta)}$ from the interval $[\delta_0, \delta_1]$ to the space of sections of $P \times_{\operatorname{Ad}(SU(2))} SL(2; \mathbb{C})$ which is asymptotic to $\mathbb{I}$ in the three relevant limits, which is equal to $\mathbb{I}$ at $\delta = \delta_0$, and such that the $s = \infty$ end member of any give $\delta \in [\delta_0, \delta_1]$ family $\{(A|_s, \varphi|_s, \mathfrak{a}_3|_s)\}_{s \in [0, \infty]}$ can be obtained from the $\delta = \delta_0$ end-member via (7.8) using $\mathfrak{g} = \mathfrak{g}_{(\delta)}$.

Now write the δ–derivative of $\mathfrak{g}_{(\cdot)}$ at a given value of δ as $(v+iu)\mathfrak{g}_{(\delta)}$ with $v$ and $u$ being sections of ad(P). These will necessarily limit to zero as $t \to 0$, $t \to \infty$ and $\frac{t}{|z|} \to 0$ because $\mathfrak{g}_{(\cdot)}$ limits to 1. The section $u$ in particular will obey (4.3) with $\mathfrak{X} \equiv 0$ because each $s = \infty$ end member data set obeys all four bullets of (1.3). (With out that fourth bullet, the right hand side of (4.3) would be $\frac{d}{d\delta}\mathfrak{X} + [v, \mathfrak{X}]$.) And, if $u$ obeys the $\mathfrak{X} = 0$ version of (4.3), its norm $|u|$ obeys the inequality $\Delta |u| \geq 0$ which implies via the maximum principal (and given $|u|$'s asymptotics) that $u \equiv 0$. The latter conclusion implies that $\mathfrak{g}_{(\cdot)}$ is unitary which implies in turn that the $s = \infty$ end member of any given $\delta \in [\delta_0, \delta_1]$ version of the family $\{(A|_s, \varphi|_s, \mathfrak{a}_3|_s)\}_{s \in [0,\infty]}$ is Aut(P) equivalent to the $s = \infty$ end member of the $\delta = \delta_0$ version of the family.

**Appendix on Hardy's inequality**

Three versions of Hardy's inequality were invoked in this paper, these were depicted respectively in (4.15), (5.48) and (6.31). They all stem from the following sequence of inequalities: Let λ denote a real number which is not the integer 1. Let $f$ denote a compactly supported function on $(0,\infty)$. Then

$$\int_{(0,\infty)} f^2 t^{\lambda-2}\, dt = \tfrac{1}{\lambda-1} \int_{(0,\infty)} f^2 \, d(t^{\lambda-1}) \,,$$

(A.1)

which can be written using integration by parts as

$$\int_{(0,\infty)} f^2 t^{\lambda-2}\, dt = -\tfrac{2}{\lambda-1} \int_{(0,\infty)} f(\tfrac{d}{dt} f)\, t^{\lambda-1}\, dt \,,$$

(A.2)

which leads via the Cauchy-Schwarz inequality to

$$\int_{(0,\infty)} f^2 t^{\lambda-2}\, dt \leq -\tfrac{4}{(\lambda-1)^2} \int_{(0,\infty)} |\tfrac{d}{dt} f|^2 \, t^{\lambda}\, dt \,,$$

(A.3)

The λ=0 version of this is (4.15) and the λ = -2 version is (6.31). The version in (5.48) would be the λ = 2 case if the integration domain were the whole of $\mathbb{R}^3$ in which case t would be the radial coordinate when using spherical coordinates for $\mathbb{R}^3$. When the integration domain for the λ = 2 case is $[0, \rho]$, then the integration by parts leading from (A.2) to (A.2) has a boundary term at ρ which changes (A.1) to this:

$$\int_{(0,\rho]} f^2 \, dt = 2 \int_{(0,\rho]} f(\tfrac{d}{dt}f) \, t \, dt + \rho \, f^2(\rho) \, .$$

(A.4)

This leads to the desired inequality

$$\int_{(0,\rho]} f^2 \, dt \leq c_0 \int_{(0,\rho]} (|\tfrac{d}{dt}f|^2 + \tfrac{1}{\rho^2} f^2) \, t^2 \, dt$$

(A.5)

by writing

$$f^2(\rho) = \int_{(0,\rho]} \tfrac{d}{dt}(\chi(2(1-\tfrac{t}{\rho})) f^2) \, dt$$

(A.6)

and using the fact that the integrand is zero where t is less than $\tfrac{1}{2}\rho$.